\documentclass{siamltex}
\usepackage{graphicx}
\usepackage{bm}
\usepackage{amsmath}
\usepackage{amssymb}
\usepackage[usenames, dvipsnames]{color}
\usepackage[breaklinks,colorlinks=true,allcolors=blue]{hyperref}
\font\msbm=msbm10

\newcommand{\E}{\hbox{{\msbm \char "45}}}
\newcommand{\R}{\hbox{{\msbm \char "52}}}
\newcommand{\bsc}[2]{#1^{(#2)}}
\newcommand{\bv}[1]{{\bf #1}}\usepackage{graphicx}
\usepackage{bm}
\usepackage{amsmath}
\usepackage{amssymb}

\newcommand{\iv}[2]{\bv{#1}^{(#2)}}

\newcommand{\vfld}[1]{\!\vec{\,#1}}

\newcommand{\ifiss}{{\sc ifiss}}
\newcommand{\vxi}{\boldsymbol \xi} 
\newcommand{\vphi}{\boldsymbol \phi}
\newcommand{\vdelta}{\boldsymbol \delta}
\newcommand{\vgamma}{\boldsymbol \gamma}

\newcommand{\calA}{\mathcal{A}}
\newcommand{\calC}{\mathcal{C}}
\newcommand{\calD}{\mathcal{D}}
\newcommand{\calE}{\mathcal{E}}
\newcommand{\calF}{\mathcal{F}}
\newcommand{\calH}{\mathcal{H}}
\newcommand{\calL}{\mathcal{L}}
\newcommand{\calM}{\mathcal{M}}
\newcommand{\calN}{\mathcal{N}}
\newcommand{\calV}{\mathcal{V}}
\newcommand{\Hone}{\calH^1(\calD)}
\newcommand{\uh}{\vfld{u}_h}
\newcommand{\vh}{\vfld{v}_h}
\newcommand{\deltah}{\vfld{\delta}_h}
\newcommand{\uhs}{\vfld{u}^{\,(s)}_h}
\newcommand{\Xh}{X_0^h}

\definecolor{otherblue}{rgb}{0,0.3,0.6}

   \def\rblack#1{{\textcolor{black}{#1}}}

\definecolor{gray}{rgb}{0.57,0.64,0.69}

\title{Collocation Methods for Exploring Perturbations in Linear Stability 
Analysis\thanks{This work was supported by the U.\, S.\, Department of 
Energy Office of Advanced Scientific Computing Research, Applied Mathematics 
program, under Award Number DE-SC0009301,
by the U.\, S.\, National Science Foundation under grant DMS1418754,
and by the U.\ K.\ EPSRC under grant EP/P013317.}}

\author{
Howard C. Elman\thanks{Department of Computer Science, University of Maryland,
College Park, MD 20742, USA (\tt{elman@cs.umd.edu})}
\and
David J. Silvester\thanks{
School of Mathematics, University of Manchester, UK (\tt{d.silvester@manchester.ac.uk})}.
}

\begin{document}

\maketitle

\begin{abstract}
Eigenvalue analysis is a well-established tool for stability analysis of
dynamical systems. 
However, there are situations where eigenvalues miss some important features 
of physical models.  
For example, in models of incompressible fluid dynamics, there are examples 
where linear stability analysis predicts stability but transient simulations exhibit
significant growth of infinitesimal perturbations.  
This behavior can be predicted by pseudo-spectral analysis.  
In this study, we show that an approach similar to pseudo-spectral analysis 
can be performed inexpensively using stochastic collocation methods and the 
results can be used to provide quantitative information about 
instability.  In addition, we demonstrate that the results of the perturbation analysis provide 
insight into the behavior of unsteady flow simulations.
\end{abstract}

\begin{keywords}
stability analysis, collocation, pseudospectra, flow simulation
\end{keywords}

\begin{AMS}
65L07, 65F15, 65P40
\end{AMS}

\pagestyle{myheadings}

\thispagestyle{plain}

\markboth{H.\ C.\ ELMAN and D.\ J.\ SILVESTER}
{Stochastic Collocation for Stability Analysis}

\section{Introduction} \label{sect-intro}
This study is concerned with a refined understanding of the classic problem of stability of 
dynamical systems.
Let 
\begin{equation} \label{dynamical}
\frac{\partial u}{\partial t} = f(u,t),
\end{equation}
represent a dynamical system, where 
$u:\R^d \times [0,T] \to \R$, $f:\R \times [0,T] \to \R$,
and let $u^{(s)}$ denote a steady solution to (\ref{dynamical}), i.e., 
$$
\frac{\partial \bsc{u}{s}}{\partial t} = f(\bsc{u}{s},t) = 0 \quad \mbox{for all } t. 
$$
Let $\gamma = \gamma(x,0)$ represent a small perturbation of $\bsc{u}{s}$.
Suppose the perturbed quantity $\hat{u}(x,0) := \bsc{u}{s}(x) + \gamma(x,0)$ is taken
as an initial condition for (\ref{dynamical}), for which integration leads to a solution
$\hat{u}(x,t) = \bsc{u}{s}(x) + \gamma(x,t)$.
If $\hat{u}(x,t)$ reverts to $\bsc{u}{s}(x)$ ($\gamma(x,t)\to 0$) as $t$ increases, then the 
steady solution is said to be stable; otherwise it is unstable.
In typical applications, $f(u,t) = f_{\alpha}(u,t)$ depends on a parameter $\alpha$, as 
does the resulting steady solution $\bsc{u}{s}_{\alpha}$, and we are interested in the set of 
values of  $\alpha$ for which $\bsc{u}{s}_{\alpha}$ is stable.

Spatial discretization of (\ref{dynamical}) leads to a 
discrete version of it, which has the form
\begin{equation} \label{discrete-dynamical}
M \frac{\partial \bv{u}}{\partial t} = \bv{f}(\bv{u},t), 
\end{equation}
where $\bv{u}$ and $\bv{f}(\bv{u},t)$ are finite-dimensional 
vectors of size $n_u$, the size of the spatial discretization.
For finite-element discretization, $M$ is a mass matrix.
As above, we wish to know if a steady solution $\iv{u}{s}$ to (\ref{discrete-dynamical}) 
is stable.

{\em Linear stability analysis} addresses this question by 
examining the eigenvalues of the algebraic system 
\begin{equation} \label{eig-linear-stability}
J \bv{v} = \lambda M \bv{v} ,
\end{equation}
where $J=\frac{\partial \bv{f}}{\partial \bv{u}}(\iv{u}{s})$ is the Jacobian matrix of $\bv{f}$ 
with respect to $\bv{u}$, evaluated at $\iv{u}{s}$; see, for example, \cite[Ch.\ 1]{Govaerts}.
A necessary condition for stability of $\iv{u}{s}$ is that all eigenvalues $\lambda$ of 
(\ref{eig-linear-stability}) have negative real part.
If any eigenvalue has positive real part, then there exists an arbitrary small perturbation 
$\vgamma$ such that if $\iv{u}{s}+\vgamma$ is used as an initial condition for 
(\ref{discrete-dynamical}), the integrated solution will not revert to $\iv{u}{s}$.

A problematic aspect of linear stability analysis is that it fails to account for transient 
effects that may take a long time to resolve.
In particular, it may happen that the solution of the system (\ref{dynamical}) with initial
condition $u(x,0) = \bsc{u}{s}(x) + \gamma(x)$, consisting of a small perturbation of a 
steady solution, exhibits large growth over a significant period of time even if $\bsc{u}{s}$ is
linearly stable.
This is discussed for models of flow in 
\cite[Sections 2.3,4.1]{Schmid-Henningson},  
\cite[Sections 20,22]{Trefethen-Embree}.
It can be explained using {\em pseudospectra}:
the $\epsilon$-{\em pseudospectrum} of the Jacobian matrix, defined for $M=I$ in 
(\ref{eig-linear-stability}), is the set of eigenvalues of $J+E$ for $\|E\|\le \epsilon$.
(A generalization to forms of $M$ considered in the present study is discussed in 
\cite{Embree-Keeler}.)\
Transient growth is exhibited when some elements of this set protrude into the 
right-half of the complex plane \cite{Trefethen-Embree}.

Our aim in this study is to develop and explore a simple procedure to study the sensitivity
of the eigenvalues of (\ref{eig-linear-stability}) when the dynamical system comes from 
models of incompressible flow.
As observed in \cite{Embree-Keeler}, it is not practical to compute pseudospectra for 
the large-scale systems  
\rblack{that arise in this setting.
This difficulty is addressed in \cite{Embree-Keeler} by projecting such systems
into invariant subspaces of (shifted versions of)  $J^{-1}M$, which have smaller 
dimension and for which computation of pseudospectra is feasible. 
It is shown in \cite{Embree-Keeler} that these pseudospectra estimates provide interior
bounds on pseudospectra of (\ref{eig-linear-stability}) as well as insight into transient
growth of solutions.
}

\rblack{In this work, we develop a complementary approach
to study the sensitivity of the eigenvalues of (\ref{eig-linear-stability}) 
for} models of incompressible flow.
The methodology derives from a two-fold procedure:
\begin{enumerate} \itemsep 0pt
\item
Introduce a simple way to construct perturbed versions of the eigenvalue problem 
(\ref{eig-linear-stability}) using spatial perturbations that depend on a finite number of 
parameters.
\item
Approximate the critical eigenvalues of the perturbed problem using a surrogate function
defined by interpolation.
\end{enumerate}
This requires the solution of a relatively small number of perturbed eigenvalue problems 
determined from a special set of parameter values, using sparse-grid  methods
\cite{Barthelmann-Novak-Ritter,Smolyak}.
The surrogate function interpolates the critical eigenvalues obtained from these eigenvalue
problems and provides a means of approximating the critical eigenvalues for an additional
set of perturbed problems.
The surrogate function is  very inexpensive to evaluate.
As a result, it is possible to generate many samples of (approximate) eigenvalues in order
to gain an understanding of the effects of perturbation.
We apply this technique to the eigenvalue problems arising from stability analysis of the 
incompressible Navier-Stokes equations.

An outline of the remainder of the paper is as follows.
In Section \ref{sect-approach}, we describe the collocation strategy and show in detail
how it is developed for the Navier-Stokes equations.
In Section \ref{sect-benchmark-and-eigs}, we describe two benchmark problems we use
to test the methodology and show how the perturbed eigenvalues behave with respect to 
Reynolds numbers and sizes of perturbation, and 
in Section \ref{sect-transient}, we demonstrate that the behavior of perturbed eigenvalues
predicts the behavior of transient solutions obtained from perturbed flow conditions.
Finally, in Section \ref{sect-conclusions}, we make some concluding remarks.

\section{Approach} \label{sect-approach}
In this section, we describe the methodology we will use to explore the sensitivity to
perturbation of the eigenvalue problem (\ref{eig-linear-stability}), which is based on
sampling.
We first outline the approach in general terms in Section \ref{sect-gen-approach},
and then we continue in Section \ref{sect-NS-details} with a more detailed statement of how
the ideas are applied to a specific benchmark problem, the incompressible Navier-Stokes 
equations.

\subsection{General approach} \label{sect-gen-approach}
Let $\iv{u}{s}$ be a steady solution to (\ref{discrete-dynamical}) and let $\vdelta$ be a small 
perturbation of $\iv{u}{s}$.
We will specify $\vdelta=\vdelta(\vxi)$ to depend on a vector of parameters
$\vxi := (\xi_1,\xi_2,\ldots,\xi_m)^T$ with $\vdelta(\bv{0})= \bv{0}$, and
we will explore a perturbed eigenvalue problem
\begin{equation} \label{perturbed-eig-linear-stability}
\hat{J}(\iv{u}{s},\vdelta(\vxi))\, \bv{v} = \hat{\lambda}(\vxi)\, M \bv{v},
\end{equation}
with the aim of understanding the impact of the perturbation $\vdelta$ on
the eigenvalues
$\{\hat{\lambda}\}$.
One way to define $\hat{J}$ is to evaluate the Jacobian at the perturbed
velocity, $\hat{J}(\iv{u}{s},\vdelta):=J(\iv{u}{s}+\vdelta)$.
In this study, which concerns the incompressible Navier-Stokes equations,
we will 
insist that the perturbation is not dissipative.
Details on the structure of the perturbation and its parameter dependence
are given in Section~\ref{sect-NS-details}. 

\vspace{.05in}

\noindent
{\bf Remark 2.1.}
We call attention here to an important aspect of the issue under study.
Classic linear stability analysis concerns the sensitivity of the steady solution $\iv{u}{s}$
to perturbation. Our (different) concern here, like that of  \cite{Trefethen-Embree}, 
is the sensitivity of the eigenvalues $\lambda$ to perturbation,
and in particular whether the conclusions reached from stability analysis predict behavior.
To highlight this distinction, we use different symbols for perturbation depending on context:
$\gamma$ is used for perturbations arising in linear stability analysis, and $\delta$ for perturbations
of eigenvalue problems as in (\ref{perturbed-eig-linear-stability}).

\vspace{.05in}

Given the eigenvalue problem (\ref{perturbed-eig-linear-stability}), let
\begin{equation} \label{stability-indicator-function}
g(\vxi):= \mbox{rightmost eigenvalue of (\ref{perturbed-eig-linear-stability})},
\end{equation}
where, if there is a complex conjugate pair of rightmost eigenvalues, $g(\vxi)$ can be taken 
to be the eigenvalue with positive imaginary part.
One way to explore the sensitivity of (\ref{eig-linear-stability}) is by sampling
$\vxi$, that is, to evaluate $g(\vxi)$ for a large set of sample values of $\vxi$.
If this function is very sensitive, that is, if small changes in $\vdelta(\vxi)$ lead to large 
changes in $g(\vxi)$, then linear stability analysis may not provide an accurate
assessment of stability; conversely, if $g$ is not sensitive to perturbation, then linear
stability analysis is likely to yield insight.

The point of view here is that the study of perturbation is done by sampling a large number 
of nearby problems.
A potential downside is that this approach requires the solution of 
many eigenvalue problems (\ref{perturbed-eig-linear-stability}), one for each choice of 
$\vxi$ and resulting $\vdelta(\vxi)$, which tends to incur a high computational cost.  
To reduce this expense, instead of evaluating the function of
(\ref{stability-indicator-function}) 
(by solving an eigenvalue problem), we will replace $g(\vxi)$ with an approximation, 
a {\em surrogate function} $g^{(I)}(\vxi)$, which is inexpensive to compute and therefore 
can be evaluated cheaply for many samples of $\vxi$.
For this, we will use the method of {\em collocation} designed to construct approximations
to functions on high-dimensional spaces \cite{Barthelmann-Novak-Ritter,Smolyak}.  
This entails evaluation of $g(\vxi)$ at a relatively small number of special points, 
$\{\iv{\vxi}{1},\iv{\vxi}{2},\ldots,\iv{\vxi}{n_{\vxi}}\}$.
The surrogate function is then taken to be the polynomial interpolant of $g$,
\begin{equation} \label{stability-indicator-surrogate}
\iv{g}{I}(\vxi):= \sum_{k=1}^{n_{\vxi}} g(\iv{\vxi}{k}) \, \ell_k(\vxi) ,
\end{equation}
where $\{\ell_k(\vxi)\}$ are multidimensional Lagrange interpolation polynomials,
$$
\ell_k(\iv{\vxi}{\ell})=\delta_{k\ell},\quad 1\le k,\ell \le m.
$$
For the interpolation points, we use sparse grids derived from the extrema of
one-dimensional Chebyshev polynomials \cite{Barthelmann-Novak-Ritter}.

\vspace{.05in}

\noindent
\rblack{{\bf Remark 2.2.}
It might happen that there are multiple eigenvalues of (\ref{eig-linear-stability}) with the same 
rightmost real part and different imaginary parts.
In this case, $g(\vxi)$ of (\ref{stability-indicator-function}) would also be multi-valued or nearly so,
and the ideas presented here would need to be applied to each of the rightmost eigenvalues.
As long as there are not too many such values, this would have minimal impact on costs.}

\subsection{Application to the Navier--Stokes equations} \label{sect-NS-details}
We will explore these ideas when the dynamical system (\ref{dynamical}) comes from 
the incompressible Navier--Stokes equations, and we now describe a way to specify a
perturbation $\vdelta(\vxi)$ for this benchmark problem for use in (\ref{perturbed-eig-linear-stability}). 
To this end, consider the  Navier--Stokes equations 
\begin{equation} \label{Navier-Stokes}
{\renewcommand{\arraystretch}{1.3}
\begin{array}{rcl}
 \vfld{u}_t - \nu \nabla^2  \vfld{u}  +  \vfld{u} \cdot \nabla  \vfld{u}  +
 \nabla p &=\ \vfld{0}, \\
-\nabla\cdot \vfld{u} &=\ 0, 
\end{array}
}
\end{equation}
posed on a domain $\calD\subset \R^d$, $d=2$ or $3$, with boundary conditions 
$$
\vfld{u}=\vfld{w}\enskip  \mbox{on } \partial\calD_D, \quad
\nu\frac{\partial \vfld{u}}{\partial n} -  \vfld{n} p =   \vfld{0}\enskip \mbox{on }
\partial\calD_N,
$$
for $\partial \calD=\partial \calD_D \cup \partial \calD_N$ consisting of the portions of the 
boundary of $\calD$ on which Dirichlet or Neumann boundary conditions hold.
In a typical scenario (see \cite[p.\,413]{ESW}), $\vfld{w}$ is 
a time-dependent inflow function that rapidly goes to a steady state, 
and the Neumann boundary condition is applied at  an outflow boundary.
Let  $\calH^1(\calD)$ be the Sobolev space of functions on $\calD$ with first derivatives
in $L_2(\calD)$,  and let
$$
\calH^1_E := \{ \vfld{u}\in \calH^1(\calD)^d\; | \,\vfld{u}=\vfld{w} \;\,
{\rm on}\;\,\partial\calD_D\}, \quad
\calH^1_{E_0} := \{ \vfld{v}\in \calH^1(\calD)^d\; | \,\vfld{v}=\vfld{0}\;\,
{\rm on}\;\,\partial\calD_D\}.
$$
{For fixed time $t\in (0,\infty)$,} the weak formulation of (\ref{Navier-Stokes}) is to find 
${\vfld{u}(\cdot,t)} \in \calH^1_E$, ${p(\cdot,t)} \in L_2({\calD})$ such that 
\begin{equation} \label{weak-Navier-Stokes}
{\renewcommand{\arraystretch}{1.7}
\begin{array}{l}
\int_{\calD} \vfld{u}_t \cdot \vfld{v} +
\nu\int_{\calD} \nabla \vfld{u} : \nabla \vfld{v} +
\int_{\calD} (\vfld{u} \cdot \nabla \vfld{u})\cdot  \vfld{v} -\int_{\calD} p \>(\nabla\cdot \vfld{v})
   =  \int_{\calD} {\vfld f} \cdot \vfld{v} \\
\hspace{3in}  \mbox{for all  }\vfld{v}\in \calH^1_{E_0} \\
\hspace{1in} -\int_{\calD}q \, (\nabla\cdot \vfld{u}) = 0
\quad \mbox{for all }q\in L_2(\calD) .
\end{array}
}
\end{equation}

Linear stability analysis uses a linearized form of the first (momentum)
equation of 
(\ref{Navier-Stokes})--(\ref{weak-Navier-Stokes}).
Given a steady velocity field $\vfld{u}$ (i.e., $\vfld{u}_t=0$), consider a
perturbation $\vfld{u}+\vfld{\gamma}$.
Substitution of this perturbed velocity into the quadratic term from
(\ref{Navier-Stokes}) gives
$$
(\vfld{u}+\vfld{\gamma})\cdot \nabla(\vfld{u}+\vfld{\gamma}) =
\vfld{u} \cdot \nabla \vfld{\gamma} +
\vfld{\gamma} \cdot \nabla \vfld{u} +
\vfld{\gamma} \cdot \nabla \vfld{\gamma}
\approx 
\vfld{u} \cdot \nabla \vfld{\gamma} +
\vfld{\gamma} \cdot \nabla \vfld{u},
$$
where the approximation on the right is made under the
assumption that $\vfld{\gamma}$ is small.
Addition of the diffusion operator and specification of a perturbed weak
formulation then leads to a trilinear form associated with the linearized
convection-diffusion operator,
\begin{equation} \label{trilinear}
a(\vfld{\gamma},\vfld{v};\vfld{u}) :=
\nu\int_{\calD} \nabla \vfld{\gamma} : \nabla \vfld{v} +
\int_{\calD} (\vfld{u} \cdot \nabla \vfld{\gamma})\cdot  \vfld{v}  +
\int_{\calD} (\vfld{\gamma} \cdot \nabla \vfld{u})\cdot  \vfld{v}.
\end{equation}

Mixed finite-element discretization of (\ref{weak-Navier-Stokes}) uses
finite-dimensional subspaces $X_0^h$ $\subset$ ${\calH}_{E_0}^1$ and
$Y^h \subset L_2(\calD)$ together with $X_E^h \subset \calH_E^1$ containing functions
that interpolate the Dirichlet boundary data at element nodes lying in
$\partial \calD_D$.
We will assume that this discretization is div-stable \cite[Sect. 2.2]{Gunzburger}.
The discrete weak formulation is to find $\vfld{u}_h \in X_E^h$ and $p_h\in Y^h$
such that
\begin{equation} \label{discrete-weak-Navier-Stokes}
{\renewcommand{\arraystretch}{1.7}
\begin{array}{l}
\int_{\calD} [\vfld{u}_h]_t \cdot \vfld{v}_h +
\nu\int_{\calD} \nabla \vfld{u}_h : \nabla \vfld{v}_h +
\int_{\calD} (\vfld{u}_h \cdot \nabla \vfld{u}_h)\cdot  \vfld{v}_h -
\int_{\calD} p_h \>(\nabla\cdot \vfld{v}_h)
   =  \int_{\calD} {\vfld f} \cdot \vfld{v}_h \\
\hspace{3in}  \mbox{for all  }\vfld{v}_h \in X_0^h, \\
\hspace{1in} -\int_{\calD}q_h \, (\nabla\cdot \vfld{u}_h) = 0
\quad \mbox{for all }q_h \in Y^h .
\end{array}
}
\end{equation}

Let $\uhs$ be a discrete steady solution to
(\ref{discrete-weak-Navier-Stokes}), i.e., $[\uhs]_t=0$.
The eigenvalue problem (\ref{eig-linear-stability}) is derived from a 
linearized discrete formulation associated with (\ref{discrete-weak-Navier-Stokes})  
where the aim is to find eigenvalues $\lambda_h$ and associated eigenfunctions satisfying
\begin{equation} \label{functional-eigenvalue-pb}
{\renewcommand{\arraystretch}{1.7} \renewcommand{\arraycolsep}{2pt}
\begin{array}{rcll}
a(\uh, \vh ; \uhs) - \int_{\calD} p_h (\nabla\cdot \vh)
&= &\lambda_h  \int_{\calD}  \uh  \cdot \vh \qquad &\hbox{for all} \;  
\vh \in \Xh, \\
\int_{\calD} q_h (\nabla\cdot \uh) &= &0 \qquad &\hbox{for all} \; q_h \in Y^h.
\end{array}
} 
\end{equation}
Here, we have linearized around  a steady flow velocity field $\uhs$ satisfying
(\ref{discrete-weak-Navier-Stokes}). 

\vspace{.05in}

\noindent
{\bf Remark 2.3.}   
A complete discussion of the development of the trilinear form
$a(\cdot,\cdot;\cdot)$ of (\ref{trilinear})  and the derivation of (\ref{functional-eigenvalue-pb}) 
is given in \cite[Sections 8.2--8.3]{ESW}.
This form also arises from use of Newton's method for solving the
nonlinear system of equations arising from implicit time discretization of
(\ref{weak-Navier-Stokes}).

\vspace{.05in}

Let the dimensions of $X_0^h$ and $Y^h$ be $n_u$ and $n_p$, respectively.
Let $\iv{u}{s}$ be the vector of coefficients of the steady finite-element 
solution $\vfld{u}_h^{(s)}$ appearing in (\ref{functional-eigenvalue-pb}).
Then the eigenvalue problem (\ref{eig-linear-stability}) has the structure
\begin{equation} \label{NS-eig-stability}
\left[
\begin{array}{cc}
F & B^T \\ B & 0
\end{array}
\right]
\left[
\begin{array}{c}
\bv{u} \\ \bv{p}
\end{array}
\right]
=
\lambda \,
\left[
\begin{array}{cc}
-Q & 0 \\ \ 0 & 0
\end{array}
\right]
\left[
\begin{array}{c}
\bv{u} \\ \bv{p}
\end{array}
\right].
\end{equation}
Here, $F=F(\iv{u}{s})$ is the matrix of order $n_u$ derived from the 
bilinear form $a(\cdot,\cdot;\vfld{u}_h^{(s)})$, $B$ and $B^T$ are matrix 
representations of negative-divergence and gradient operators, respectively
($B$ is of size $n_p \times n_u$),
and  $Q$ is a velocity mass matrix, also of order $n_u$.

\vspace{.05in}
\noindent
{\bf Remark 2.4.}  
The matrix on the right side of  (\ref{NS-eig-stability}) is singular, and the 
resulting infinite eigenvalue can lead to
instability in eigenvalue computations~\cite{Meerbergen-Spence}.
This can be avoided by replacing the matrix by 
$\left[\!\!
{\renewcommand{\arraycolsep}{2pt}
\begin{array}{cc} -Q & \alpha B^T \\ \ \alpha B & 0 \end{array}}
\!\!\right]$, 
which leaves the finite eigenvalues intact and maps the infinite eigenvalue to $1/\alpha$,
see~\cite{Cliffe-Garratt-Spence}.\footnote{We use this variant of the mass matrix with 
$\alpha=-1/10$  in all of our computations. 
The resulting mapped eigenvalue $\lambda=-10$ is far enough from 
the near-critical ones that it  does not affect the results.}

\vspace{.05in}
We want to explore the sensitivity of  our modified eigenvalue problem to perturbation.
For this, we add a small perturbation $c(\uh,\vh;\deltah)$ to
$a(\uh,\vh;\deltah)$ in (\ref{functional-eigenvalue-pb}).
The perturbation is defined in two steps.
First, we specify a discretely divergence-free vector field $\deltah$, that is, one satisfying 
\begin{equation} \label{discrete-div-free}
\int_{\calD} q_h \, (\nabla\cdot \deltah) = 0 
\quad \hbox{for all} \; q_h \in Y^h .
\end{equation}
\rblack{This ensures that the perturbed velocity $\uh+\deltah$ would 
be appropriate as an initial condition for testing the stability of $\uh$.}
Second, the field $\deltah$ is used to generate a nondissipative perturbation 
$c_h(\uh,\vh;\deltah)$ of $a_h(\uh,\vh;\deltah)$.
\rblack{This means that the perturbation does not introduce any damping
effects associated with numerical diffusion.}
To illustrate the construction, we will suppose that $T_h$ denotes a
subdivision of $\calD \subset \R^2$ into triangular or rectangular elements. 
The extension to three-dimensional problems is  perfectly straightforward.

\vspace{.05in}

\noindent
{\bf \rblack{Claim} 2.1.}        
Suppose that $\phi_h\in \Hone$ is a finite element function defined on $T_h$ and 
that $\deltah$ is defined locally on every element $k\in T_h$ ,
$\deltah \vert_{k} := \; \deltah^{\,(k)}$, via
\begin{equation} \label{local-deltah-def}
\deltah^{\,(k)} = \hbox{curl} \; \phi_h^{(k)} =
\left[ -{\partial \phi_h^{(k)}/ \partial x_2}, 
{\partial \phi_h^{(k)}/ \partial x_1} \right]^T,
\end{equation}
so that $\deltah^{\,(k)}$ is divergence-free on each element. 
Then $\deltah := \sum_{k \in T_h} \deltah^{\,(k)}$ satisfies (\ref{discrete-div-free}).

\vspace{.05in}

\noindent
{\em Proof.}
For any function $q_h \in Y^h$,
$$
\int_{\calD} q_h\, (\nabla\cdot \deltah) = \sum_{k \in T_h}
\int_{k} q_h^{(k)} (\underbrace{\nabla\cdot  \hbox{curl} \;
  \phi_h^{(k)}}_{= 0}) =0 .
\qquad \square
$$
Note that the local construction (\ref{local-deltah-def}) generates a 
discontinuous  velocity field so that in general $\deltah\not\in {\Hone}^d$. 

\vspace{.05in}

\noindent
{\bf \rblack{Claim} 2.2.}     
Let the perturbation operator on element $k\in T_h$ be given by
\begin{equation} \label{local-skew-adjoint}
c_h^{(k)} (u_h^{(k)}, v_h^{(k)} ;\deltah^{(k)} ) :=
 \int_{k} (\deltah^{(k)} \cdot  \nabla u_h^{(k)}) \, v_h^{(k)}  
 -\frac{1}{2} \int_{\partial k} u_h^{(k)} v_h^{(k)} \;
 \deltah^{(k)}  \cdot  {\vfld{n}},
\end{equation}
Then $c_h^{(k)}$ is skew-adjoint on $T_h$, 
and 
$$ 
c_h (u_h, v_h;\deltah)  :=  \sum_{k \in T_h} c_h^{(k)}
(u_h^{(k)}, v_h^{(k)} ;\deltah^{(k)})
$$
is skew-adjoint on $\calD$.\footnote{In this discussion, $u_h$ and $v_h$
are discrete scalar functions.  For vector-valued arguments, 
e.g., $\uh=([u_h]_1,[u_h]_2)^T$, 
$c_h(\uh, \vh; \deltah) = \sum_i c_h([u_h]_i,[v_h]_i; \deltah)$
is the sum of contributions from individual scalar components. 
}

\vspace{.075in}

\noindent
{\em Proof.}
Since  $\deltah\not\in {\Hone}^d$, we need to apply  Green's theorem on 
each element:   
$$
{\renewcommand{\arraystretch}{2.2}  \renewcommand{\arraycolsep}{1pt}
\begin{array}{lcl}
\int_{k}(\deltah^{(k)}\!\cdot\!\nabla u_h^{(k)})\,v_h^{(k)}  
&= &\int_{k}(v_h^{(k)}\deltah^{(k)}) \!\cdot \!\nabla u_h^{(k)} 
 = -\int_{k} \nabla \cdot (v_h^{(k)} \deltah^{(k)}\!) u_h^{(k)} 
+ \int_{\partial k} u_h^{(k)} v_h^{(k)} \; \deltah^{(k)} \!\!  \cdot \vfld{n} \\
&= &-\int_{k} (  v_h^{(k)} \nabla \cdot \deltah^{(k)}  +
\deltah^{(k)} \! \cdot \! \nabla v_h^{(k)} )\, u_h^{(k)} 
+ \int_{\partial k} u_h^{(k)} v_h^{(k)} \; \deltah^{(k)} \!\! \cdot \vfld{n} \\
 &= &-\int_{k}  (  \deltah^{(k)} \cdot  \nabla v_h^{(k)} ) u_h^{(k)} 
+ \int_{\partial k} u_h^{(k)} v_h^{(k)} \; \deltah^{(k)} \!\! \cdot \vfld{n} ,
\end{array}
}
$$
where the last equality follows from 
\rblack{the fact that $\deltah^{(k)}$ is divergence-free on each element}.
It follows that 
$$
\int_{k} \! (\deltah^{(k)} \cdot  \nabla u_h^{(k)}) \, v_h^{(k)}  
- \frac{1}{2} \! 
\int_{\partial k} \!\!\!u_h^{(k)} v_h^{(k)} \; \deltah^{(k)} 
\!\! \cdot  {\vfld{n}}  = \!
- \!\int_{k} \! (\deltah^{(k)} \cdot  \nabla v_h^{(k)}) \, u_h^{(k)}  
+ \frac{1}{2} \! 
\int_{\partial k} \!\!\!u_h^{(k)} v_h^{(k)} \; \deltah^{(k)} 
\!\! \cdot \vfld{n} ,
$$
that is, $c_h^{(k)}$ is skew-adjoint. 
Summation over all the elements establishes the same property for $c_h$.
$\qquad \square$


The perturbed variant of  (\ref{functional-eigenvalue-pb}) is
$$
{\renewcommand{\arraystretch}{1.7} \renewcommand{\arraycolsep}{2pt}
\begin{array}{rcll}
a(\uh, \vh ; \uhs) + c_h(\uh,\vh; \deltah) - \int_{\calD} p_h (\nabla\cdot \vh)
&= &\hat{\lambda}_h  \int_{\calD}  \uh  \cdot \vh \qquad &\hbox{for all} \;  
\vh \in \Xh, \\
\int_{\calD} q_h (\nabla\cdot \uh) &= &0 \qquad &\hbox{for all} \; q_h \in Y^h.
\end{array}
} 
$$
This leads to the perturbed matrix eigenvalue problem (\ref{perturbed-eig-linear-stability}) 
\begin{equation} \label{perturbed-NS-eig-stability}
\left[
\begin{array}{cc}
F  + N(\vxi) & B^T \\ B & 0
\end{array}
\right]
\left[
\begin{array}{c}
\bv{u} \\ \bv{p}
\end{array}
\right]
=
\hat{\lambda} \,
\left[
\begin{array}{cc}
-Q & \alpha B^T \\ \  \alpha B & 0
\end{array}
\right]
\left[
\begin{array}{c}
\bv{u} \\ \bv{p}
\end{array}
\right],
\end{equation}
\rblack{where the perturbation matrix $N=N(\vxi)$ 
is determined from
\begin{equation} \label{N-from-ch}
(\bv{u},N\bv{v}) = c_h(u_h,v_h,\vec{\delta}_h)
\end{equation}
so that in particular} $N$ is a skew-symmetric matrix, $N^T = -N$, for all 
parameter values $\vxi$ independent of the boundary conditions of the 
flow problem.

It remains to specify the finite element function $\phi_h$ used  
in \rblack{Claim} 
2.1 to define the vector field $\deltah$. 
Following \cite{Phillips-Elman-stochastic}, we take
$\phi_h(x,\vxi)\in \Hone$ 
to be a parameter-dependent scalar potential specified using a
covariance function $\calC(\bsc{x}{1},\bsc{x}{2})$ for $\bsc{x}{i}\in\calD$.
In particular, given $\calC$, let $C:= \calC(\bv{x},\bv{x})$ 
be the covariance matrix of order $n$ 
consisting of the vertices in the subdivision associated with
$\Xh$, so that $C_{ij}=\calC\left(x_i,x_j\right)$.
Now let $\vphi$ be an $n$-dimensional  
zero-mean stationary random process with covariance 
matrix $C$, i.e., $C=\E(\vphi\vphi^T)$, where ``$\E$'' refers to expected
value.
If $C=\sigma^2 V\Theta V^T$ is an eigenvalue--eigenvector decomposition 
(scaled by the variance), then
$\vphi$ can be defined using a discrete Karhunen--Lo\`eve (KL) expansion 
\begin{equation} \label{full-discrete-KL}
\vphi(\vxi) := \sigma V \Theta^{1/2} \vxi = 
\sigma \sum_{j=1}^{n} \sqrt{\theta_j}\,\bv{v}_j\, \xi_j,
\end{equation}
where the eigenvector $\bv{v}_j$ is the $j$th column of $V$ and
$\{\xi_j\}_{j=1}^{n}$
are uncorrelated random variables with zero mean and unit variance
\cite[Section 5.4]{LPS}.
It is often the case that many of the eigenvalues are small and some of the 
terms in (\ref{full-discrete-KL}) can be removed without significant loss of 
accuracy.
We will choose $m<n$ such that 
$\left(\sum_{j=m+1}^{n}\theta_j\right) 
\left/\left(\sum_{j=1}^{n}\theta_j\right) \right. \le 5/100$,
and, in the sequel, $\vxi:= (\xi_1,\ldots,\xi_m)^T$ will represent an
$m$-dimensional vector of parameters and $\vphi(\vxi)$ is defined using the
truncated KL-expansion
\begin{equation} \label{truncated-discrete-KL}
\vphi(\vxi) := \sigma \sum_{j=1}^m \sqrt{\theta_j}\, \bv{v}_j\, \xi_j.
\end{equation}
This $\vxi$-dependent coefficient vector (of length $n$) then characterizes a 
piecewise-defined  linear or bilinear  function $\phi_h(\vxi)$. 
For the computational results described in Section~\ref{sect-benchmark-and-eigs},  we 
take the smooth covariance function
\begin{equation} \label{covariance-function}
\calC(\bsc{x}{1},\bsc{x}{2}):= \sigma^2
\exp \left({\textstyle -\left[
    \left(\dfrac{\textstyle \bsc{x}{1}_1-\bsc{x}{2}_1}{\textstyle c_1}\right)^2
    + \left(\dfrac{\textstyle \bsc{x}{1}_2-\bsc{x}{2}_2}{\textstyle c_2}\right)^2
    \right]} \right) ,
\end{equation}
where $c_1$ and $c_2$ are correlation lengths.
We will also assume that $\{\xi_j\}$ in (\ref{truncated-discrete-KL}) 
are mutually independent, with each satisfying a 
truncated Gaussian distribution with range $[-3,3]$, so that $\xi_j$ has the density function
$$
\rho(\xi) = \left\{
{\renewcommand{\arraystretch}{1.7}
\begin{array}{ll}
 \frac{1} {\operatorname{erf} (3\sqrt{2})} \, \frac{1}{\sqrt{2\pi}}  \, 
 \exp\left(\frac{\xi^2}{2}\right) 
  & \mbox{for } |\xi| \le 3\\ 
  \  0 & \mbox{for } |\xi| > 3.
\end{array}
}
\right.
$$
A  {\sc Matlab} implementation of this distribution is given in
\cite{Botev-truncated-code} and described in \cite{Botev-truncated-ref}.

\rblack{We note two differences between this formulation of perturbations and traditional 
approaches based on pseudospectra.
First, the perturbed eigenvalue problem (\ref{perturbed-eig-linear-stability}), as specified
in (\ref{perturbed-NS-eig-stability}), is restricted to have a structure determined from
that of the original problem, so that the resulting perturbed eigenvalues are closer in form
to  {\em structured pseudoeigenvalues}  \cite[Ch.\ 50]{Trefethen-Embree}.
Moreover, the perturbation itself derives explicitly from the nonlinear term 
$\vfld{u} \cdot \nabla  \vfld{u}$ in (\ref{Navier-Stokes}).
Indeed, if the original problem (\ref{dynamical}) were linear so that the Jacobian $J$ did not 
depend on $\bv{u}$, then the eigenvalues of $J$ might still be sensitive to perturbation but 
the ideas discussed here would not give insight into this. 
Despite these limitations, as will be shown in Sections 
\ref{sect-benchmark-and-eigs}-\ref{sect-transient}, 
the behavior of the perturbed eigenvalues gives insight into transient growth
of solutions and other features of transient solvers.
}

\rblack{We conclude this section with an analytic result bounding the
  size of the eigenvalue perturbation in proportion to the perturbation
  size.
}
\rblack{
\begin{theorem} \label{eig-perturb-bound}
If $\delta_h$ is the perturbation defined using {\rm (\ref{local-deltah-def})} and $\lambda$ is the 
rightmost eigenvalue of {\rm (\ref{eig-linear-stability})}, then the eigenvalue $\hat \lambda$ 
of the perturbed problem {\rm (\ref{perturbed-eig-linear-stability})} closest to $\lambda$ satisfies 
$|\hat \lambda - \lambda| \le (c/h)\, \| \delta_h \|_{\infty}$.
\end{theorem}
}
\begin{proof}
\rblack{
Let $\calM$ denote the matrix on the right side of (\ref{perturbed-NS-eig-stability}), 
and let $\calF$ and $\widehat{\calF}=\calF + \calN$ denote the unperturbed and perturbed 
matrices on the left sides of (\ref{NS-eig-stability}) and (\ref{perturbed-NS-eig-stability}), 
respectively; here $\calN = \calN(\vxi) =
\left[\begin{array}{cc}N(\vxi) & 0 \\ 0 & 0 \end{array}\right]$. 
We are interested in the eigenvalue problems $\calF v = \lambda \calM v$ 
and $\widehat{\calF} \hat v = \hat \lambda \hat \calM v$.
$\calM$ can be factored as
\begin{equation} \label{calM-factor1}
\calM = 
\left[
\begin{array}{cc}
I & 0 \\ -\alpha B Q^{-1} & I
\end{array}
\right]
\left[
\begin{array}{cc}
-Q & 0 \\ 0 & \alpha^2 B Q^{-1} B^T 
\end{array}
\right]
\left[
\begin{array}{cc}
I & -\alpha Q^{-1} B^T \\ 0 & I
\end{array}
\right] .
\end{equation}
The velocity mass matrix $Q$ and Schur complement  $\alpha^2 BQ^{-1}B^T$ 
each admit a Cholesky decomposition, $Q=LL^T$, $\alpha^2 BQ^{-1}B^T=RR^T$, so that
(\ref{calM-factor1}) can be refined to
$$
\begin{array}{c}
\calM = 
\underbrace{\left[
\begin{array}{cc}
I & 0 \\ -\alpha B Q^{-1} & I
\end{array}
\right]
\left[
\begin{array}{cc}
L & 0 \\ 0 & R 
\end{array}
\right]
}
\ 
\underbrace{\left[
\begin{array}{cc}
-I & 0 \\ 0 & I
\end{array}
\right]
}
\ 
\underbrace{\left[
\begin{array}{cc}
L^T & 0 \\ 0 & R^T 
\end{array}
\right]
\left[
\begin{array}{cc}
I & -\alpha Q^{-1} B^T \\ 0 & I
\end{array}
\right]}. \\
\hspace{.17in} \calL \hspace{1.2in} \calD \hspace{1.2in} \calL^T
\end{array}
$$
Thus, we can consider the unperturbed and perturbed standard eigenvalue problems
\begin{equation} \label{standard-eig-pbs}
\calD^{-1}\calL^{-1}\calF\calL^{-T} w = \lambda w, \quad 
(\calD^{-1}\calL^{-1}\calF\calL^{-T} +
\calD^{-1}\calL^{-1}\calN\calL^{-T} )
\hat w = \hat \lambda \hat w.
\end{equation}
Let $\calA=\calD^{-1}\calL^{-1}\calF\calL^{-T}$ and $\calE=\calD^{-1}\calL^{-1}\calN\calL^{-T}$.
If $\calA =\calV \Lambda  \calV^{-1}$ is diagonalizable, then we can use the Bauer-Fike Theorem 
\cite[Ch.\ 7]{GvL} to explore the eigenvalue perturbation:
given an eigenvalue $\lambda$ of $\calA$, 
$$
\min_{\hat \lambda \in \sigma(\calA+\calE)} |\hat \lambda -\lambda| \le \kappa(\calV) \|\calE\|_2\,.
$$
We seek a bound on $\|\calE\|_2 = \|\tilde \calE\|_2$, where $\tilde \calE = \calL\calN\calL^{-T}$;
the equality here follows from the fact that $\calD$ is unitary. 
The structure of $\calL$ leads to
$$
\tilde \calE = \calL^{-1}\calN\calL^{-T} = 
\left[
\begin{array}{cc}
L^{-1}NL^{-T} & \alpha L^{-1}NQ^{-1}B^TR^{-T} \\
\alpha R^{-1}B Q^{-1}N L^{-T} & \alpha^2 R^{-1}BQ^{-1}NQ^{-1}B^TR^{-T}
\end{array}
\right],
$$
so that $\|\tilde \calE \|_2 \le \| L^{-1}NL^{-T}\|_2 + O(\alpha)$.
This holds for all $\alpha \ne 0$, and so it also holds in the limit as $\alpha \to 0$, giving 
$\|\tilde \calE \|_2 \le \| L^{-1}NL^{-T}\|_2$.
Since both $N$ and $L^{-1}NL^{-T}$ are skew-symmetric,
it follows that $\|\calE\|_2\le \rho(L^{-1}NL^{-T})$, the spectral radius.
}

\rblack{
Thus, we require a bound on the Rayleigh quotient $\frac{|(y,Ny)|}{(y,Qy)}$.
For this, we use (\ref{N-from-ch}) together with a standard bound on $|c_h|$ (see 
\cite[p.\ 243]{ESW}) to get
$$
\frac{|(y,Ny)|}{(y,Qy)} = \frac{|c_h(y_h,y_h;\delta_h)|}{\|y_h\|_{L_2(\calD)}^2} \le 
\frac{\|\delta_h\|_{\infty}\, \|\nabla y_h\|_{L_2(\calD)}\, \|y_h\|_{L_2(\calD)}} {\|y_h\|_{L_2(\calD)}^2} = 
\frac{\|\delta_h\|_{\infty}\, \|\nabla y_h\|_{L_2(\calD)}} {\|y_h\|_{L_2(\calD)}} \, .
$$
The assertion then follows from the inverse estimate $\|\nabla y_h\|_{L_2(\calD)}/\|y_h\|_{L_2(\calD)} \le O(1/h)$. 
}
\end{proof}

\section{Benchmark problems and structure of eigenvalues} \label{sect-benchmark-and-eigs}
We will illustrate these ideas for two benchmark problems.
In this section, we specify the problems and their features of interest,
the eigenvalues associated with linear stability analysis and the 
effect of perturbation of these eigenvalues.
Each of these is a model of flow through a channel for which there are inflow and 
outflow boundaries. 
The position of the outflow boundary is far enough downstream that the flow is fully developed.
The spatial approximation is done using $Q_2$--$P_{-1}$ (biquadratic velocity; discontinuous
linear pressure) mixed approximation \cite[Section 3.3.1]{ESW},
implemented in the \ifiss\ software package
\cite{IFISS-paper1,IFISS-paper2}.
\rblack{Unless otherwise specified, 
the discretization is done on a uniform grid with element width $h=1/32$,
which gives $64$ elements along the vertical interval $[-1,1]$.
This corresponds to ``grid level'' $\ell=6$ in \ifiss, with element width
$h=2/2^{\ell}$.}
\rblack{For both benchmark problems, we will explore the stability of solutions obtained
  for choices of the viscosity parameter near the critical values for linear stability.
  We determined these critical values experimentally, as described for example
  in \cite{Elman-Meerbergen-Spence-Wu}.
}

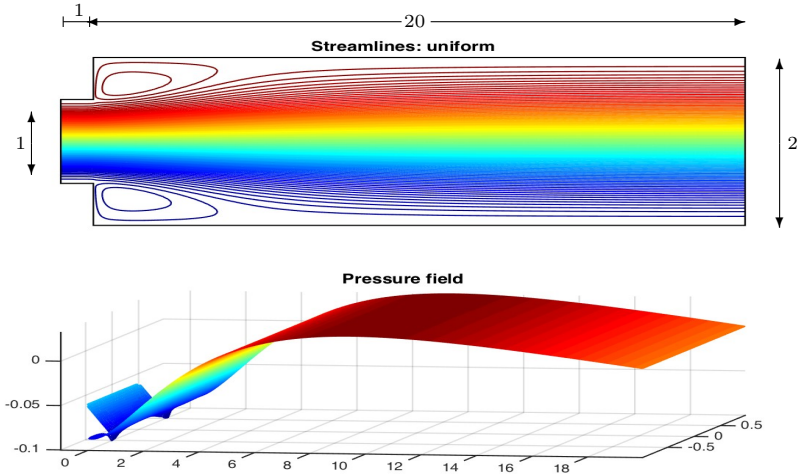
\begin{figure}
\begin{center}
\begin{picture}(290,190) (24,0) 
\put(0,0){\includegraphics[width=.9\textwidth, height=.5\textwidth]{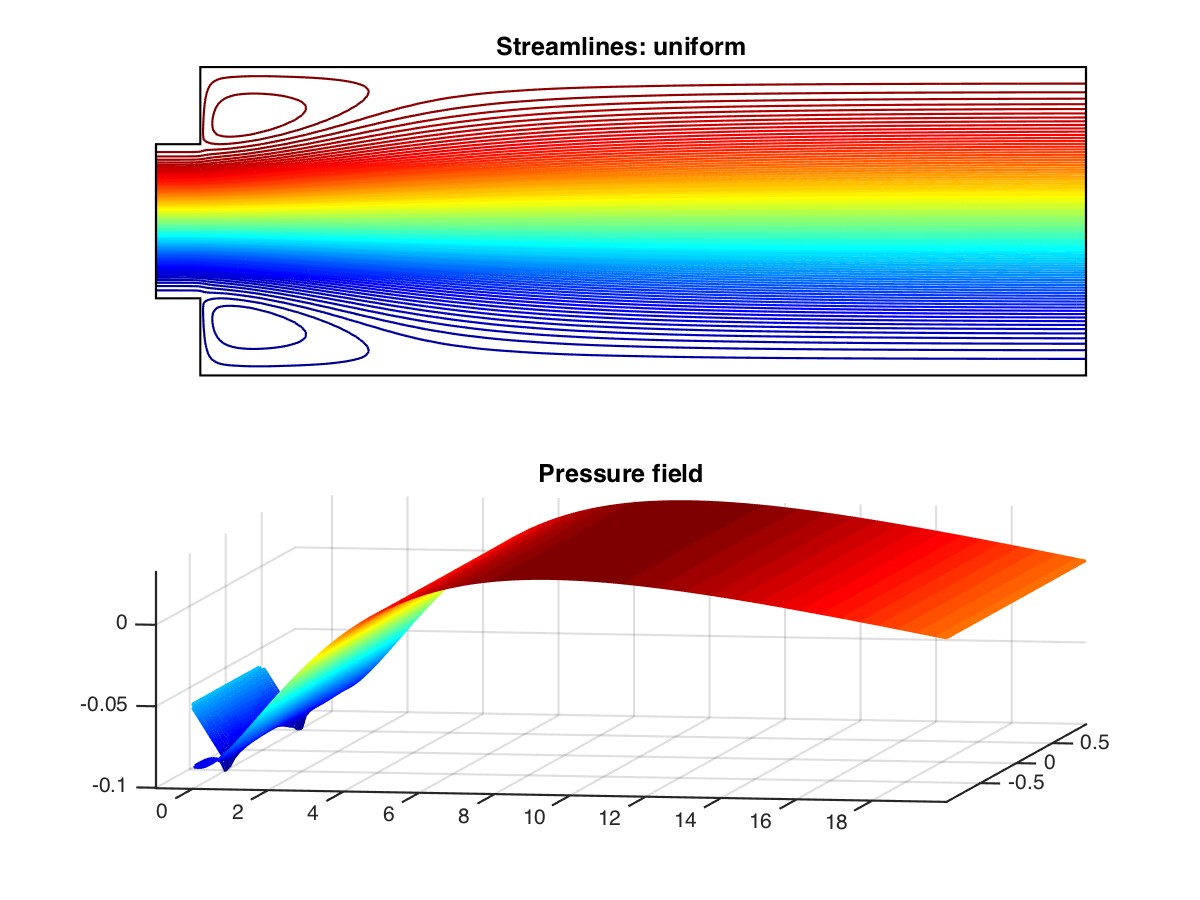}}
\put(32,139){\vector(0,1){12}}
\put(32,139){\vector(0,-1){12}}
\put(26,137){{\scriptsize 1}}
\put(44,185){\line(1,0){10}}
\put(44,183){\line(0,1){4}}
\put(54,183){\line(0,1){4}}
\put(48,187){{\scriptsize 1}}
\put(315,139){\vector(0,1){32}}
\put(315,139){\vector(0,-1){32}}
\put(318,137){{\scriptsize $2$}}
\put(183,185){\vector(1,0){119}}
\put(172,185){\vector(-1,0){119}}
\put(173,183){{\scriptsize $20$}}
\end{picture}
\end{center}
\vspace{-.2in}
\caption{Symmetric step domain and velocity/pressure solutions for $\nu=1/220$.} 
\label{step-figure}
\end{figure}

\subsection{Expansion flow around a symmetric step} \label{benchmark-step}
The domain $\calD$ is a rectangular duct with a symmetric expansion, with
boundary conditions
\vspace{.05in}
\begin{itemize}  \itemindent -.3in \itemsep 3pt
\item 
parabolic profile  $\vfld{u}(-1,y) = (1-4y^2,0)$ at the inflow boundary 
$(-1,y),\,|y|\le .5$
\item
natural conditions $\nu \frac{\partial u_x}{\partial x}=p$, 
$\frac{\partial u_y}{\partial x}=0$  at the outflow boundary  $(20,y),\,|y|\le 1$
\item
no-flow conditions $\vfld{u}=0$  along fixed walls 
 
 $(x,\pm 1), 0\le x \le 20$;
  $(x, \pm .5), -1 \le x \le 0$;   $(0,y)$, $.5 \le |y| \le 1$.
\end{itemize}
\vspace{.05in}

\noindent
Details of the domain and a sample solution are shown in Figure~\ref{step-figure}.
The discretization is defined on a uniform grid of square elements.
The key feature  of this solution is that it is  reflectionally symmetric
with respect to the centerline $y=0$, i.e., the stream function $\psi$ satisfies $\psi(x,y)=-\psi(x,-y)$.
It follows that for the velocity, 
\begin{equation} \label{step-flow-character}
u_x(x,y)=u_x(x,-y), \quad u_y(x,y)=-u_y(x,-y).
\end{equation}

This flow problem exhibits a pitchfork bifurcation~\cite{drikakis97}:
as the viscosity decreases through a critical value (approximately $\nu=1/220.5$),
the rightmost eigenvalue of (\ref{NS-eig-stability}), 
which is real, changes from negative (indicating linear stability) to positive (instability).
Figure \ref{eigs_step} shows the ten rightmost eigenvalues for three values of $\nu$ in this
range and the rightmost eigenvalues (detail in the inset) for each choice, whose
values are also identified on the right.

\begin{figure}
\begin{center}
\begin{picture}(200,190) (30,0) 
\put(0,0){\includegraphics[width=.7\textwidth,height=.5\textwidth]{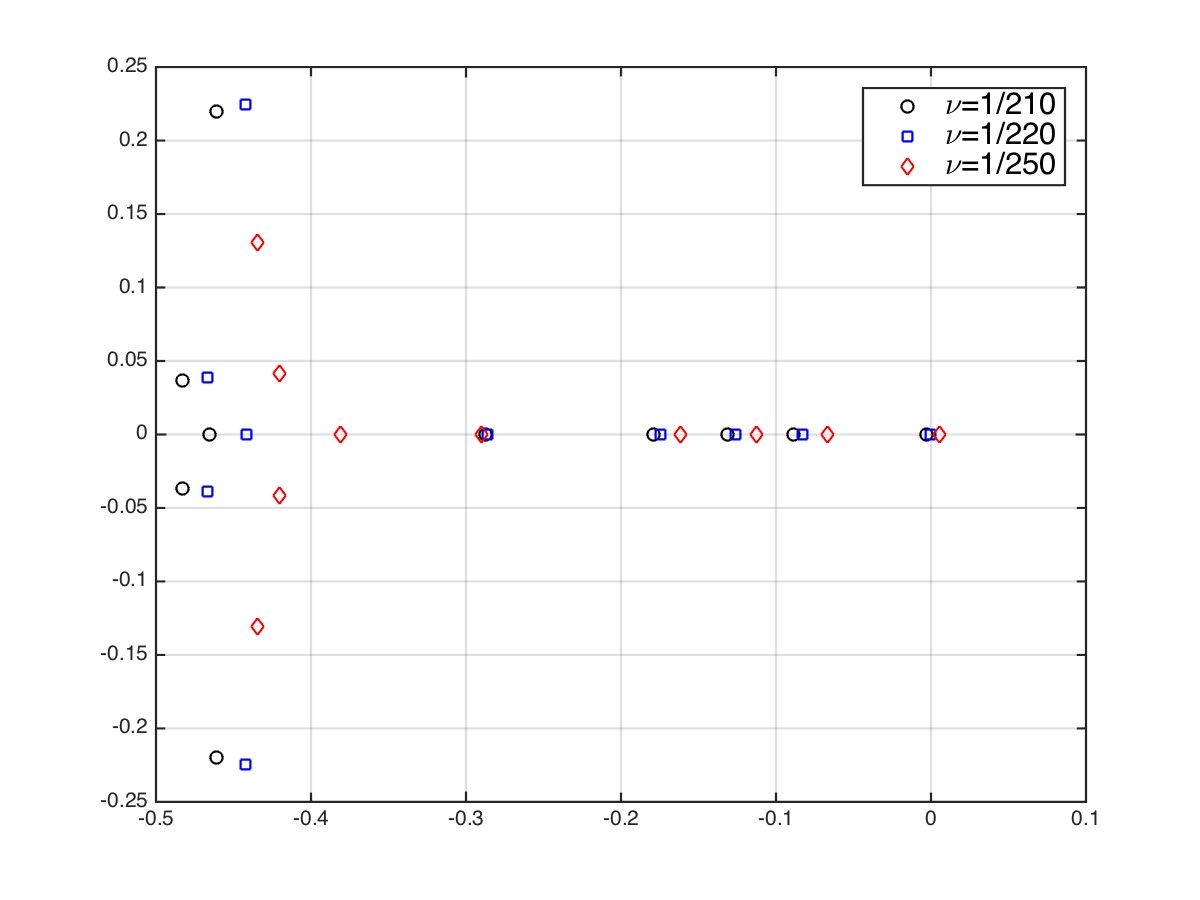}}
\put(160,30){\includegraphics[width=.15\textwidth,height=.09\textwidth]{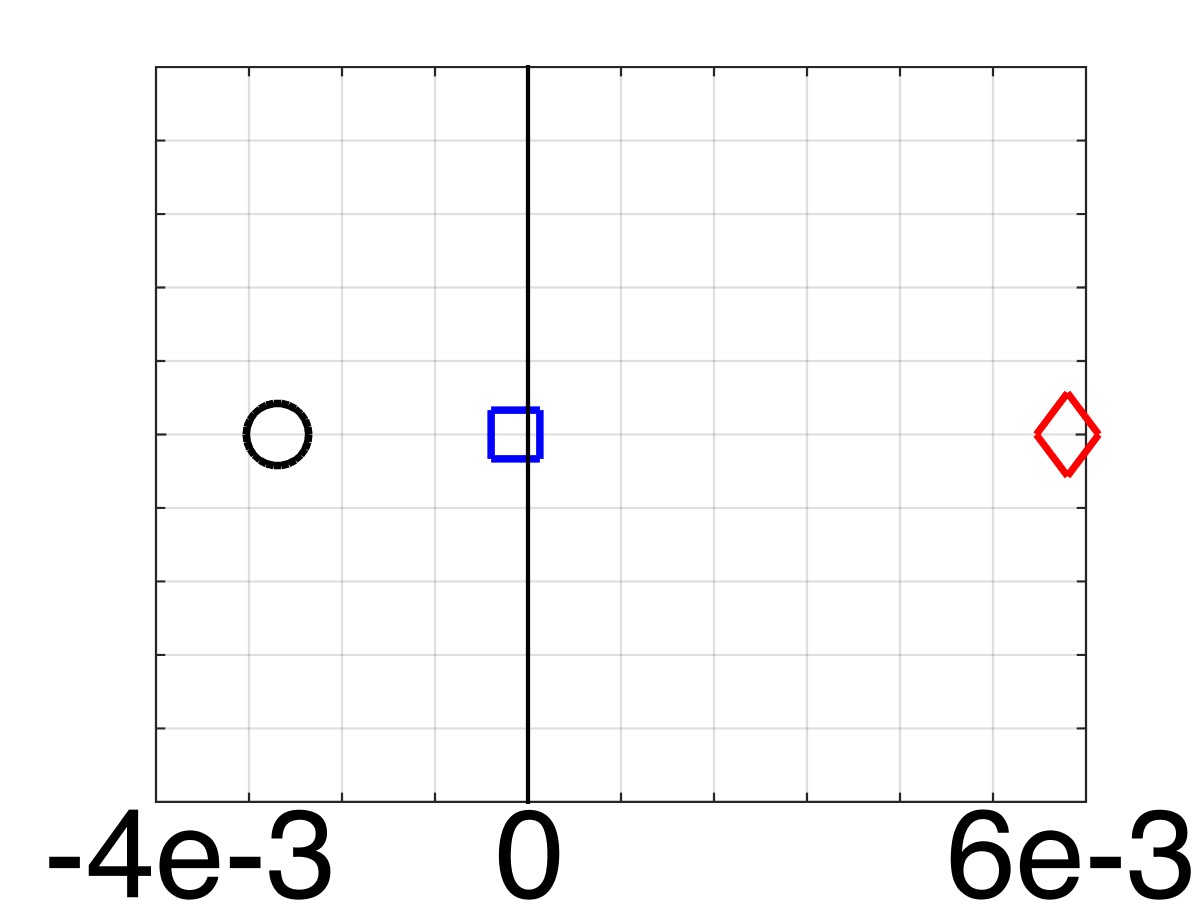}}
\put(200,90){\vector(-1,-2){12}}
\put(255,60){\scriptsize $\nu$}     \put(290,60){\scriptsize $\lambda$}
\put(245,50){\scriptsize $1/210$} \put(273,50){\scriptsize $-2.7\times 10^{-3}$}
\put(245,40){\scriptsize $1/220$} \put(273,40){\scriptsize $-1.4\times 10^{-4}$}
\put(245,30){\scriptsize $1/250$} \put(279,30){\scriptsize $5.8\times 10^{-3}$}
\end{picture}
\end{center}
\vspace{-.2in}
\caption{Eigenvalues for the symmetric step problem}.
\label{eigs_step}
\end{figure}
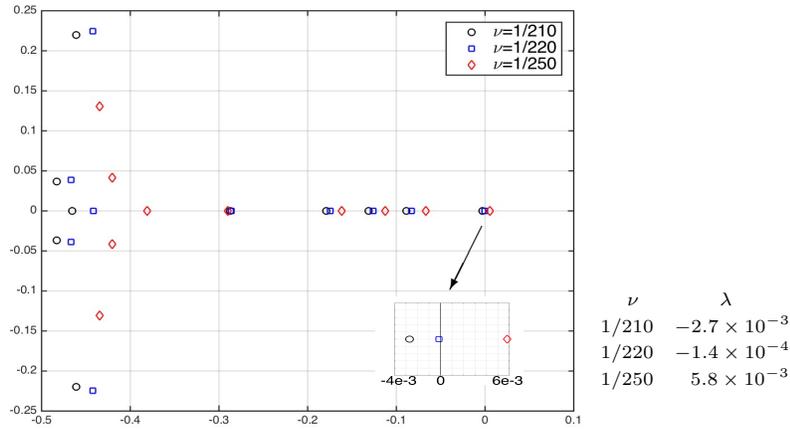

\begin{figure}
\begin{center}
\begin{picture}(300,155) (50,0) 
\put(0,0){\includegraphics[width=.55\textwidth,height=.4\textwidth]{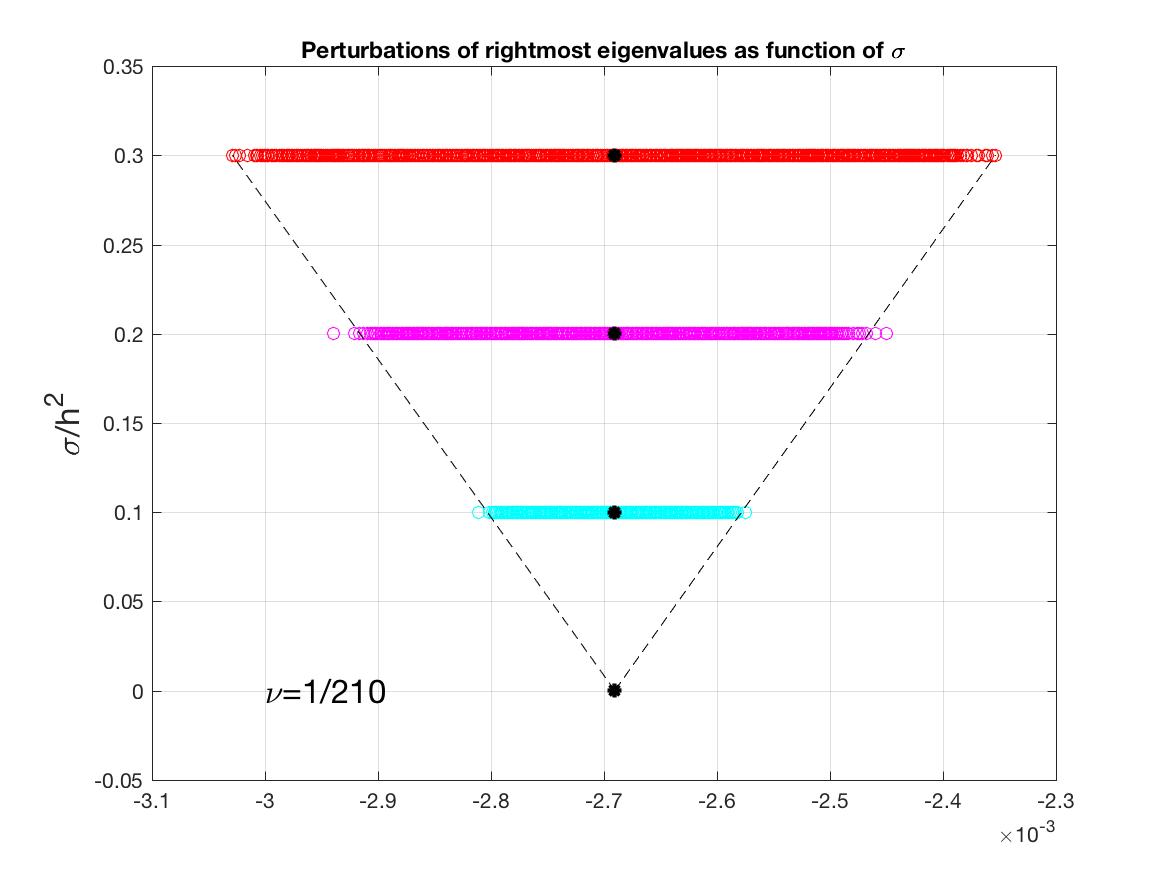}}
\put(195,-3){\includegraphics[width=.55\textwidth,height=.41\textwidth]{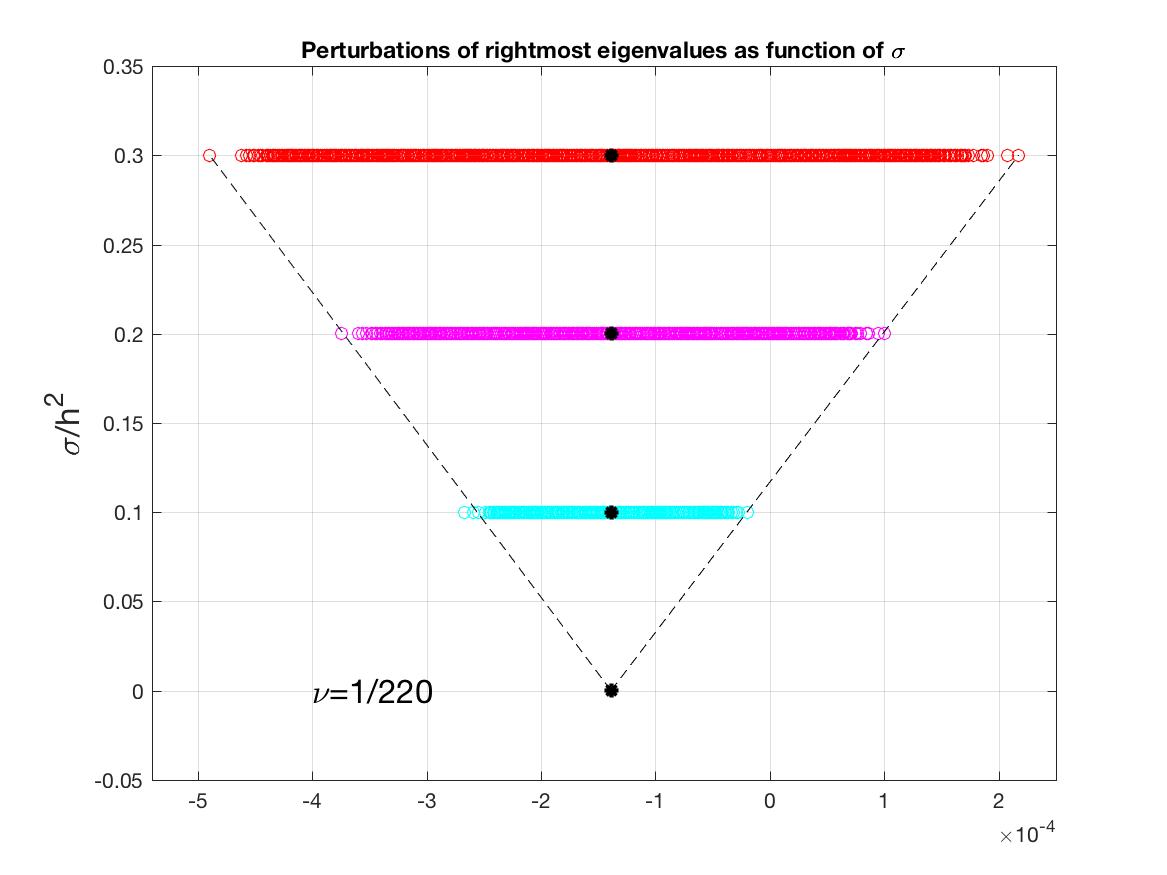}}
\end{picture}
\end{center}
\vspace{-.2in}
\caption{Surrogate perturbed rightmost eigenvalues for the symmetric step 
problem, for $\sigma=.1h^2$, $.2h^2$, and $.3h^2$, and $\nu=1/210$ (left) and 
$\nu=1/220$ (right).} 
\label{step_eigs_perturbations}
\end{figure}

We explored the sensitivity of the rightmost negative eigenvalues using the perturbed
eigenvalue problem (\ref{perturbed-NS-eig-stability}).
This was derived using correlation lengths $c_1=c_2=2$ in (\ref{covariance-function}), which
resulted in a finite expansion (\ref{truncated-discrete-KL}) with $m=19$ terms.
The surrogate function $g^{(I)}$  of (\ref{stability-indicator-surrogate}) used to estimate eigenvalues
was constructed from a two-level sparse grid on the $m$-dimensional parameter space,
which in turn resulted in $n_{\vxi}=761$ sparse grid nodes.
Thus, it is necessary to solve $760$ eigenvalue problems, that is, find
the rightmost eigenvalues of $760$ perturbed systems (\ref{perturbed-NS-eig-stability}), one
for each sparse-grid node.  
(One of the sparse-grid nodes is $\vxi = \bv{0}$, which corresponds to an unperturbed system.)\ 
Once these are available, the estimates of eigenvalues for other choices of $\vxi$ can be 
obtained by evaluating $g^{(I)}$.
We implemented the sparse-grid interpolation using the {\sc matlab} toolbox {\sc spinterp} 
\cite{kli07,kliwoh05}.

\rblack{
The dependence of the estimated perturbed eigenvalues on the perturbation
size is illustrated in Figure \ref{step_eigs_perturbations}.
We found that insight can be provided using small values of the standard 
deviation in (\ref{truncated-discrete-KL}), and in these tests we used
$\sigma=\beta h^2$ for $\beta=.1$, $.2$ and $.3$.
For each $\sigma$, the  figure shows the distribution of 
one million eigenvalue estimates, computed 
}
\rblack{using the interpolant (\ref{stability-indicator-surrogate}),
with results} 
for $\nu=1/210$ shown on the left and for $\nu=1/220$ on the right.
For both values of $\nu$, the rightmost unperturbed eigenvalue (center of
the set of perturbations) is negative, showing that
the associated steady solution is stable, but for the smaller, closer-to-critical value $\nu=1/220$,
some of the \rblack{estimated} perturbed eigenvalues are positive, whereas all the 
perturbations are negative for $\nu=1/210$.
The two figures have the same horizontal scaling, indicating that the magnitude of the 
perturbations does not depend on $\nu$.
The bounding dashed lines show that the magnitude of perturbations varies linearly with 
$\sigma$, 
\rblack{as the bound of Theorem \ref{eig-perturb-bound} suggests.}

\begin{figure}
\begin{center}
\begin{picture}(300,155) (50,0) 
\put(100,0){\includegraphics[width=.55\textwidth,height=.4\textwidth]{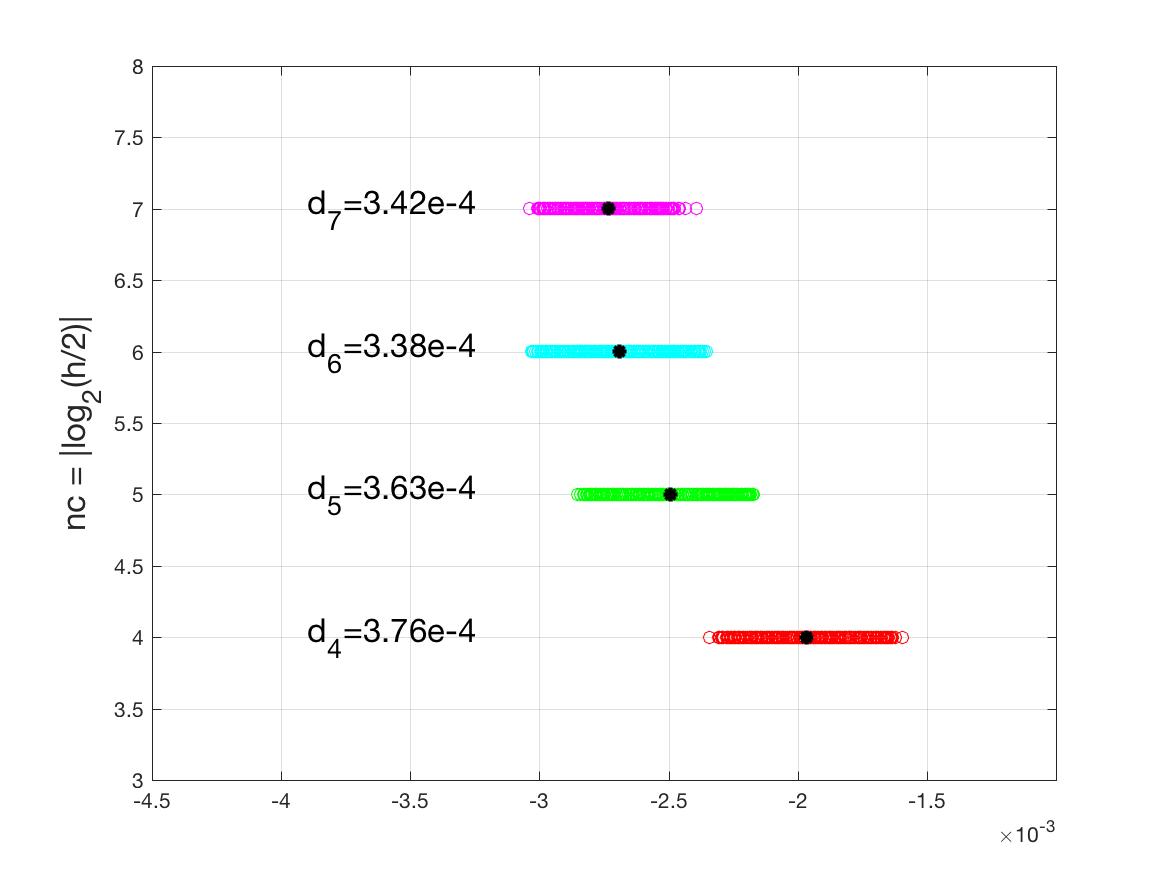}}
\end{picture}
\end{center}
\vspace{-.2in}
\caption{Rightmost eigenvalues together with surrogate perturbed rightmost eigenvalues, for the 
symmetric step problem with $\nu=1/210$, 
$\sigma=.3\times(1/32)^2$ and various mesh sizes.}
\label{eigs_perturbations_hdependence}
\end{figure}

\rblack{
Finally, Figure \ref{eigs_perturbations_hdependence} shows the behavior of 
the eigenvalue perturbations as well as the critical eigenvalues (highlighted 
in the middle of each set of perturbations) as the discretization mesh size 
varies.  
Results are shown for four mesh sizes, $h=1/8$, $1/16$, $1/32$ and $1/64$ 
(corresponding to grid levels $\ell=4$ through $7$). 
For each $\ell$, $d_{\ell}$ is the maximal difference 
$|\hat \lambda - \lambda|$ between the surrogate rightmost perturbed
eigenvalue and the rightmost true eigenvalue among all surrogate values.
The range of perturbations appears to be independent of the discretization,
which suggests that the $O(1/h)$ dependence of the bound in Theorem 
\ref{eig-perturb-bound} is pessimistic.
We attribute this to limitations on what can be obtained using the Bauer-Fike 
Theorem for the transformed problems (\ref{standard-eig-pbs}).
Also note that the critical eigenvalues move to the left with mesh refinement, indicating that stability of
the discrete systems is enhanced with mesh refinement, although it is also clear that a limiting
value is approached with refinement.
}

\subsection{Flow around a square obstacle} \label{benchmark-obstacle}
In this case, the domain is a rectangular duct containing a square obstacle, with boundary conditions
\vspace{.05in}
\begin{itemize}  \itemindent -.3in \itemsep 3pt
\item 
parabolic profile  $\vfld{u}(-1,y) = (1-4y^2,0)$ at the inflow boundary 
$(0,y),\,|y|\le 1$
\item
natural conditions $\nu \frac{\partial u_x}{\partial x}=p$, 
$\frac{\partial u_y}{\partial x}=0$  at the outflow boundary   $(8,y),\,|y|\le 1$
\item
no flow conditions $\vfld{u}=0$ along the top and bottom walls $(x,\pm 1), 0\le x \le 8$, 
and
  
on the obstacle, a square centered at $(2,0)$ with sides of length $0.5$.
\end{itemize}
\vspace{.05in}

\noindent
\rblack{For this example, we used a level-6 stretched grid with local refinement near the obstacle.}
A representative  steady solution that retains the reflectional symmetry is 
shown in Figure \ref{obstacle-figure}.
In this case, there is a symmetry-breaking Hopf bifurcation for $\nu\approx 1/186$; that is, 
for $\nu$ in this range 
there is a complex conjugate pair of rightmost eigenvalues whose real parts change from
negative to positive as $\nu$ is reduced.
Figure \ref{eigs_obstacle} shows the $100$ smallest eigenvalues for three values of $\nu$, two
near critical ($\nu=1/175$ and $1/185.6$ and one super-critical ($\nu=1/200$), as well as
a detail of the rightmost eigenvalues.

\begin{figure}
\begin{center}
\begin{picture}(260,190) (25,0) 
\put(0,0){\includegraphics[width=.84\textwidth, height=.48\textwidth]{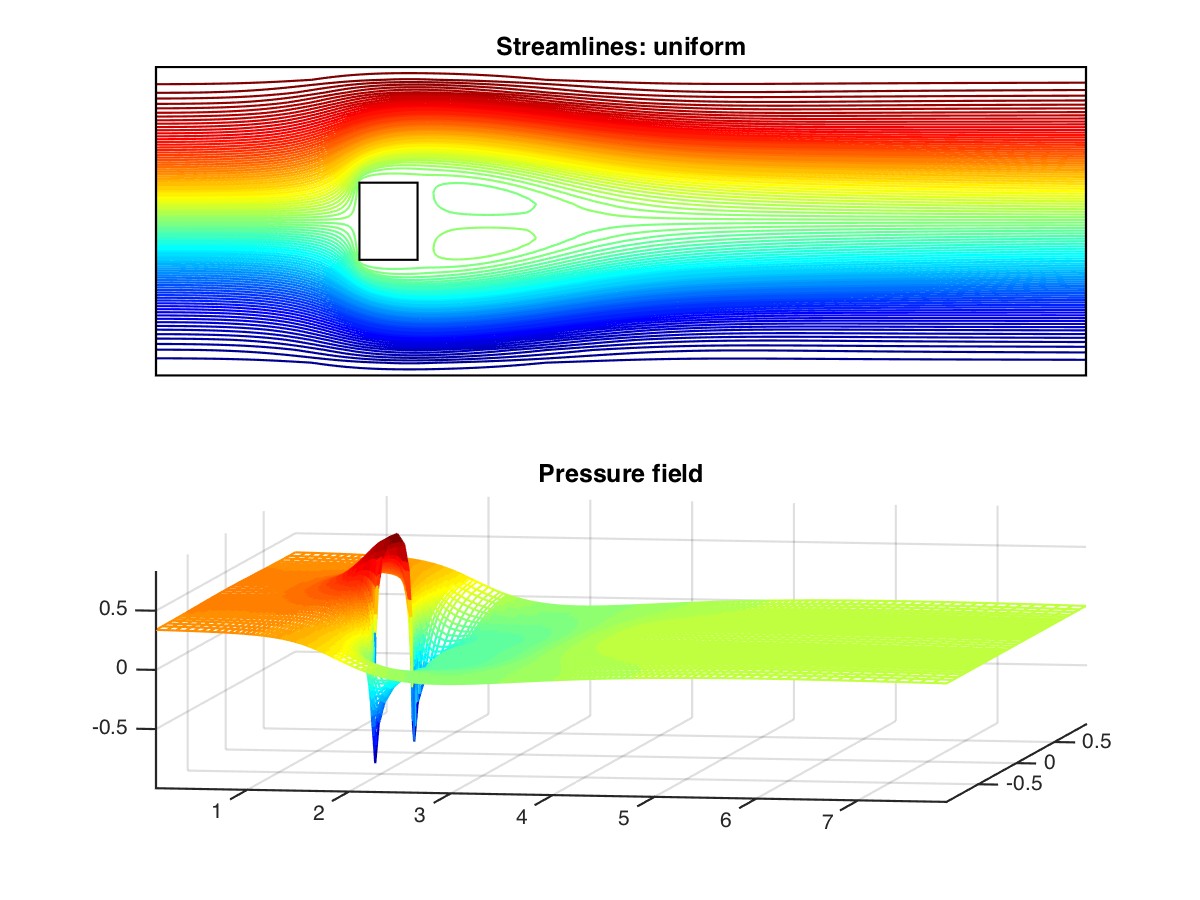}}
\put(151,178){\vector(-1,0){111}}
\put(170,178){\vector(1,0){111}}
\put(159,176){{\scriptsize $8$}}
\put(292,133){\vector(0,1){31}}
\put(292,133){\vector(0,-1){31}}
\put(295,132){{\scriptsize $2$}}
\put(93,150){\line(1,0){15}}
\put(93,148){\line(0,1){4}}
\put(108,148){\line(0,1){4}}
\put(96,144){{\scriptsize $.5$}}
\put(26,127){\line(0,1){15}}
\put(24,127){\line(1,0){4}}
\put(24,142){\line(1,0){4}}
\put(15,132){{\scriptsize $.5$}}
\end{picture}
\end{center}
\vspace{-.2in}
\caption{Obstacle domain and velocity/pressure solutions for $\nu=1/185.6$.} 
\label{obstacle-figure}
\end{figure}
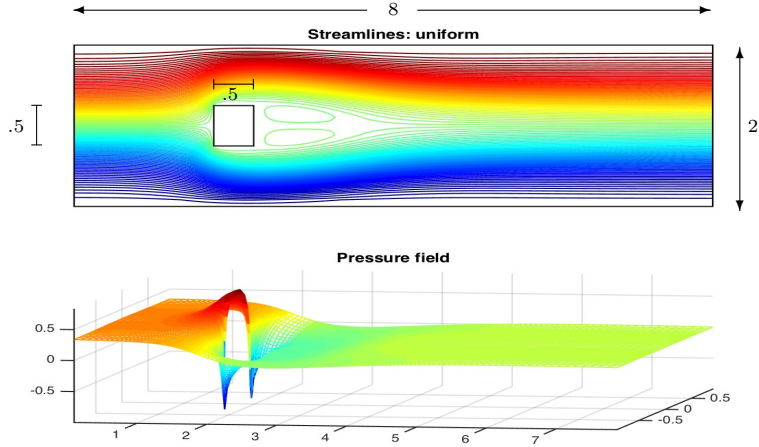
\begin{figure}
\begin{center}
\begin{picture}(200,190) (30,0) 
\put(0,0){\includegraphics[width=.7\textwidth,height=.5\textwidth]{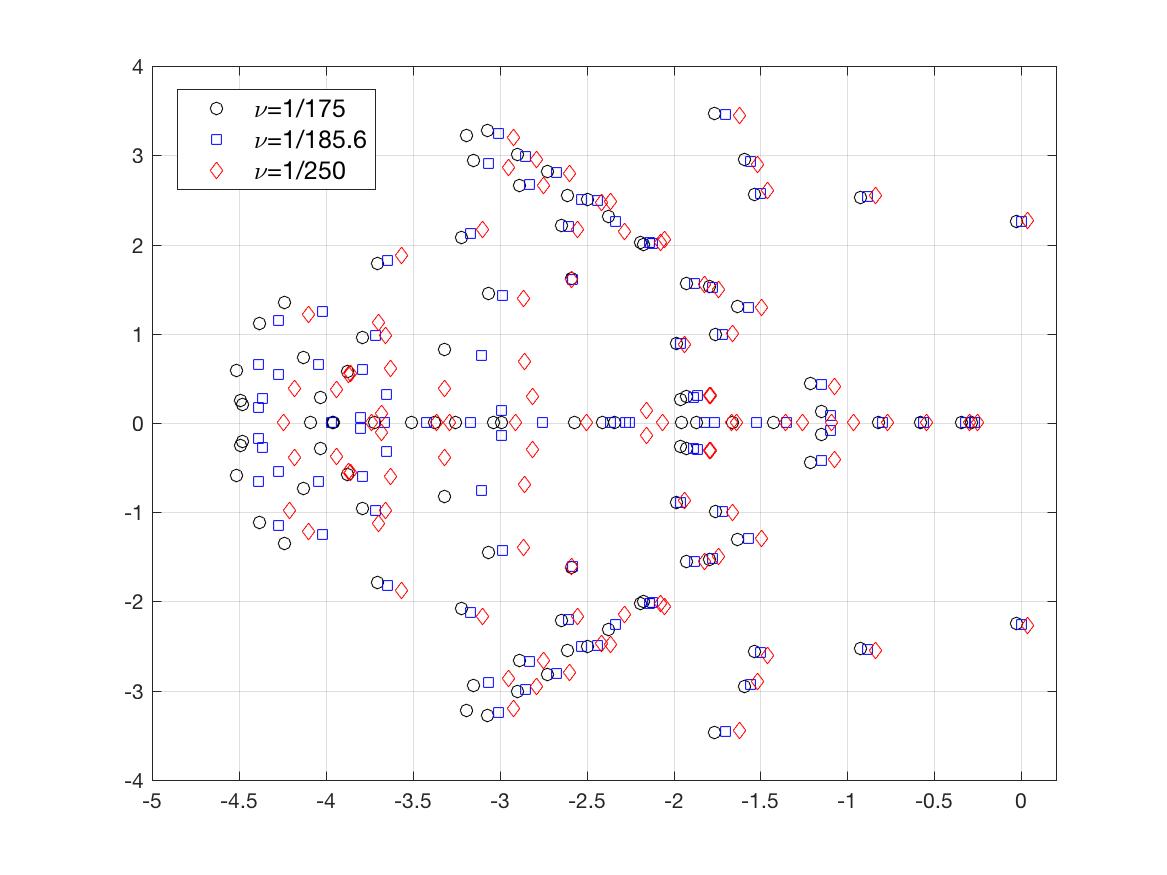}}
\put(250,135){\includegraphics[width=.15\textwidth,height=.09\textwidth]{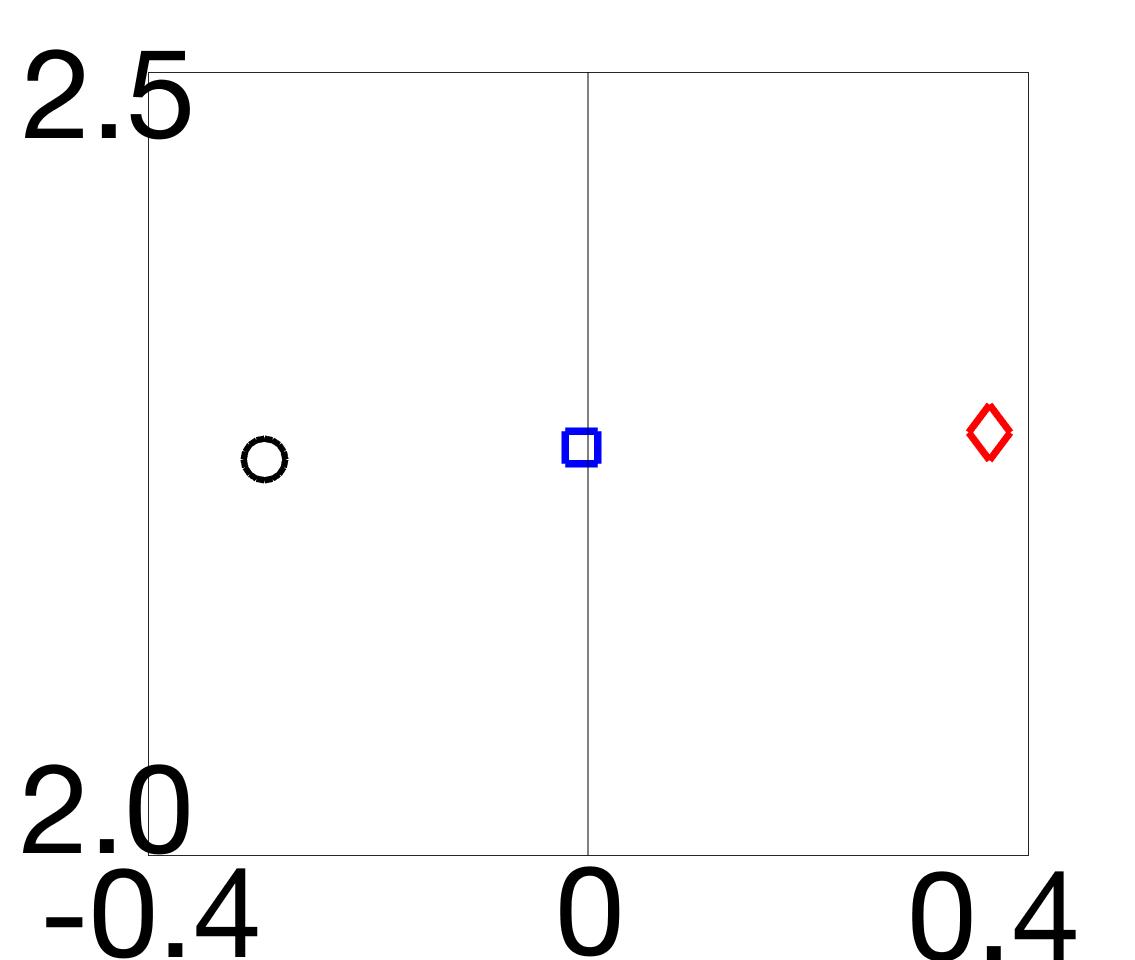}}
\put(232,140){\vector(2,1){18}}
\put(250,120){\scriptsize $\nu$}    \put(285,120){\scriptsize Re($\lambda$)}
\put(240,110){\scriptsize $1/175$} \put(273,110){\scriptsize $-2.9\times 10^{-2}$}
\put(240,100){\scriptsize $1/185.6$} \put(273,100){\scriptsize $-3.0\times 10^{-4}$}
\put(240,90){\scriptsize $1/200$} \put(279,90){\scriptsize $3.7\times 10^{-2}$}
\end{picture}
\end{center}
\vspace{-.2in}
\caption{Eigenvalues for the obstacle problem}.
\label{eigs_obstacle}
\end{figure}
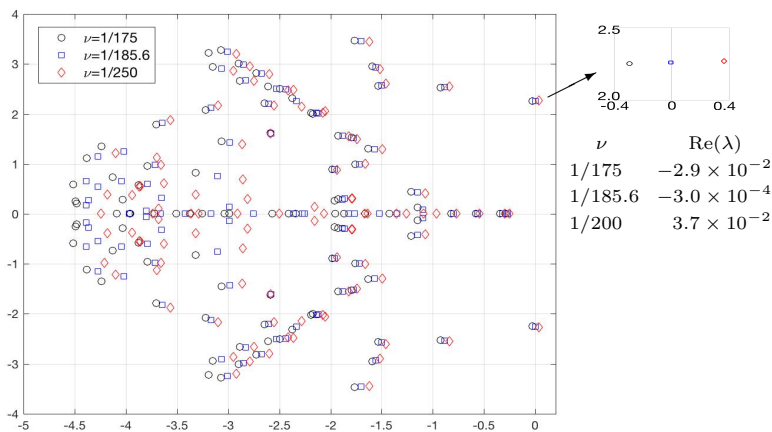

The behavior of the perturbed (estimated) eigenvalues is illustrated in 
Figure~\ref{obstacle_eigs_perturbations}.
\rblack{Once again,  we computed one million eigenvalue estimates
for  three values of the standard deviation in (\ref{truncated-discrete-KL}),
$\sigma=.1h^2$, $.2h^2$ and $.3h^2$ with $h=1/32$.}
These are shown in the figure, for $\nu=1/175$ on the left and $\nu=1/185.6$ on the right.
For both values of $\nu$, the rightmost unperturbed eigenvalue (center of
the set of perturbations) has negative real part, showing that
the associated steady solution is stable, but for the close-to-critical value $\nu=1/185.6$,
some of the perturbed eigenvalues have a positive real part.
As for the step problem,  it is readily seen that the magnitude of the  perturbations does
 not depend on $\nu$   
 and varies linearly with  $\sigma$.

\begin{figure}
\begin{center}
\begin{picture}(300,160) (50,0) 
\put(0,-3){\includegraphics[width=.55\textwidth,height=.45\textwidth]{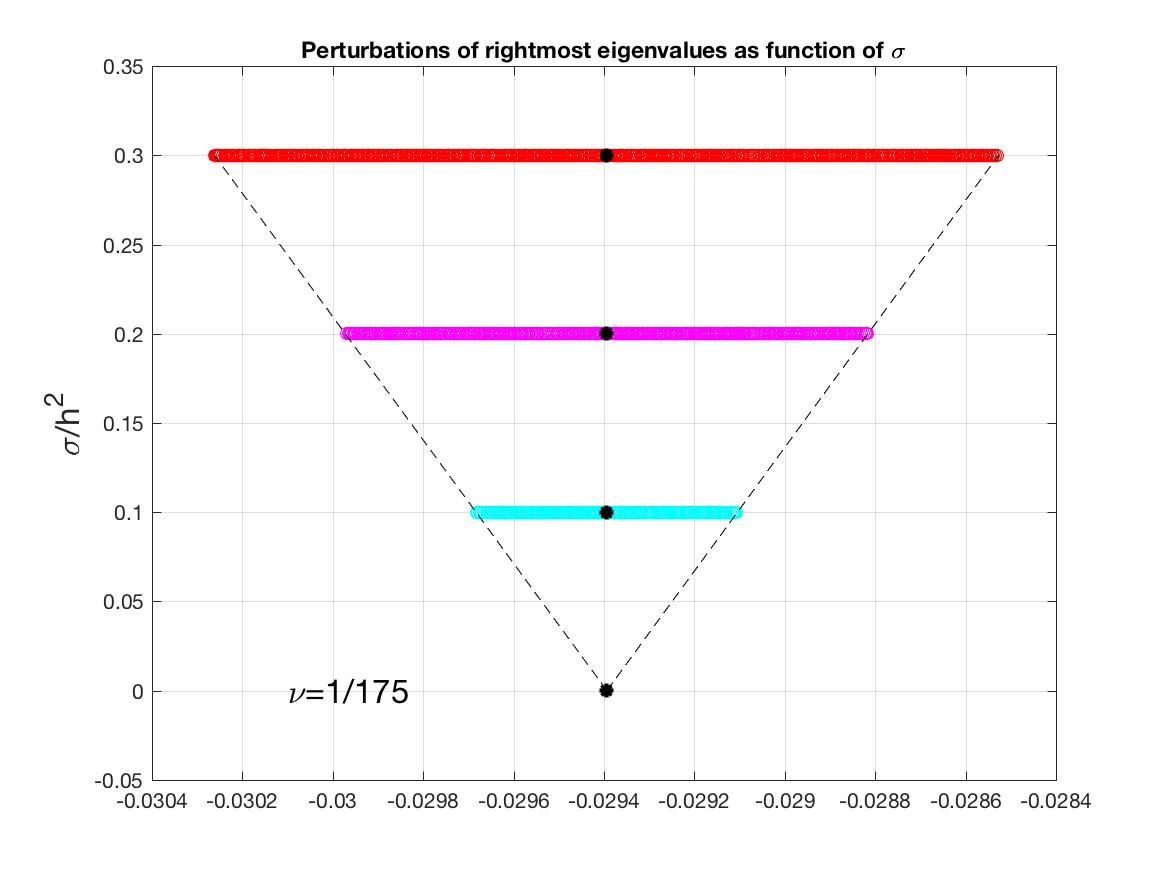}}
\put(195,-3){\includegraphics[width=.55\textwidth,height=.45\textwidth]{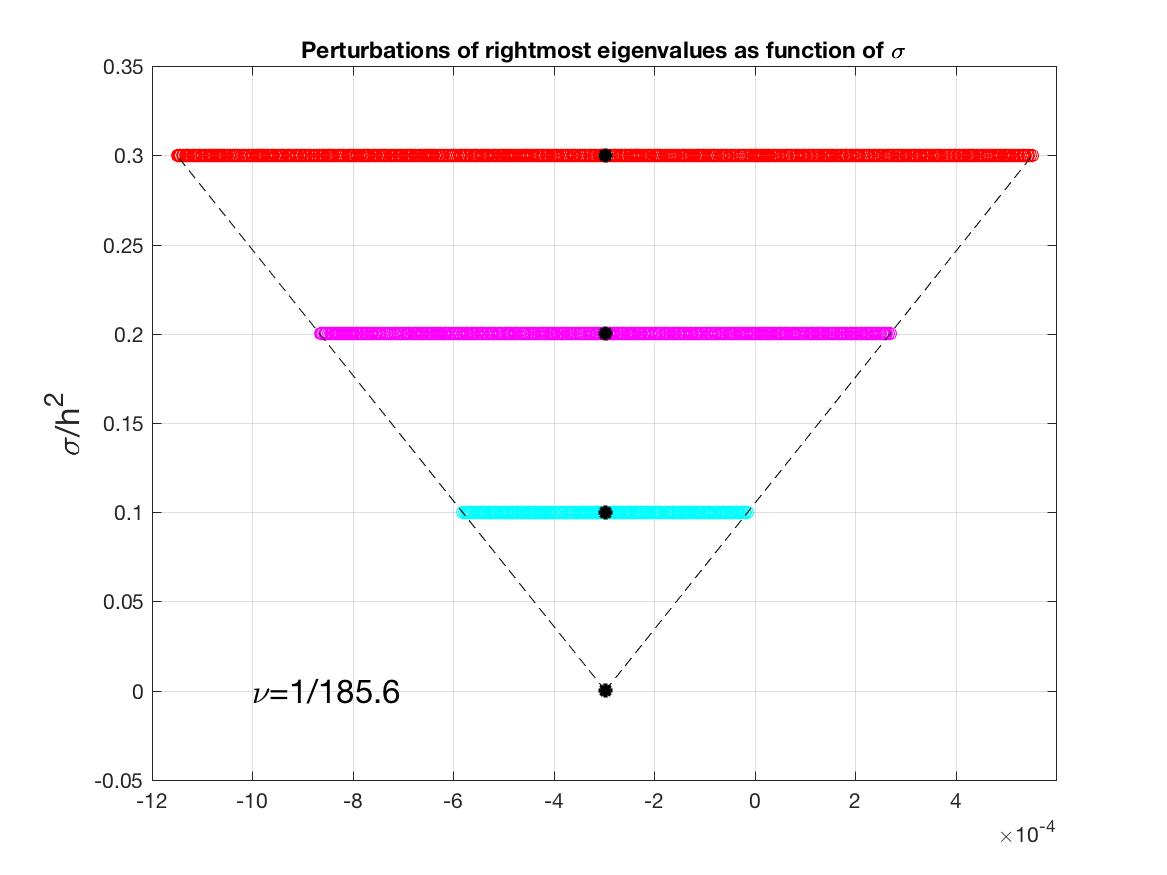}}
\end{picture}
\end{center}
\vspace{-.2in}
\caption{Real parts of surrogate perturbed rightmost eigenvalues for the obstacle problem, for 
$\sigma=.1h^2$, $.2h^2$, and $.3h^2$, and $\nu=1/175$ (left) and $\nu=1/185.6$ (right).} 
\label{obstacle_eigs_perturbations}
\end{figure}

\section{Unsteady flow simulations} \label{sect-transient}
In this section, we  explore the connection between time integration of 
 the Navier--Stokes equations and the eigenvalue perturbation results in the 
previous section. 
We will do this by  computing time-accurate solutions of the Navier--Stokes equations
using the adaptive (stabilized) Trapezoidal Rule (sTR)  time stepping methodology built 
into \ifiss.
The suitability of sTR for long-time integration is discussed  in \cite{GGS08}.
Full  details of the 
\ifiss\  implementation of  sTR can be found in section 10.2.3 of  \cite{ESW}.
(Stabilization is based on time step averaging, which prevents 
the ``ringing'' to which TR is susceptible for stiff systems.)
In what follows, we  present results obtained from a nonlinear version of the 
integrator, denoted (sTR$k$), where  a fixed number ($k=1$ or $k=2$) of  Picard corrections 
are performed at every time step.
We  present results for the benchmark problems of Sections 
\ref{benchmark-step} and \ref{benchmark-obstacle}.
Our objective is to test the sensitivity of the reference flow with respect
to instantaneous spatial perturbations,  loosely simulating a laboratory experiment 
where a reference steady flow is subject to an external disturbance, and the
flow is monitored to see if it returns  to the steady state.

\begin{figure}
\begin{center}
\begin{picture}(260,175) (51,0) 
\put(0,0){\includegraphics[width=1.0\textwidth, height=.5\textwidth]{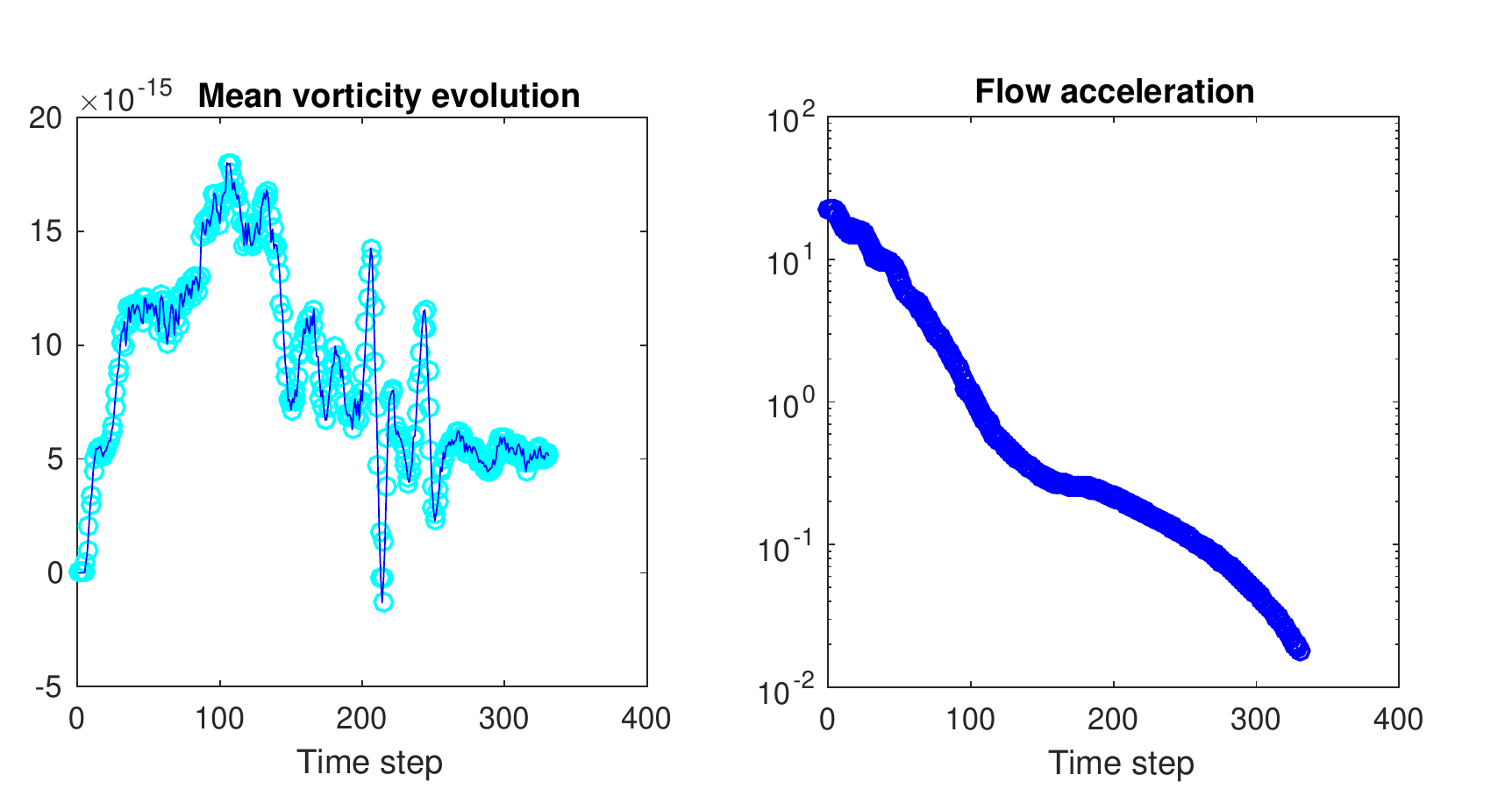}}
\end{picture}
\end{center}
\vspace{-.1in}
\caption{First phase for symmetric step flow with $\nu=1/250$.} 
\label{first-phase-250}
\end{figure}

\subsection{Evolution of expansion flow around a symmetric step}\label{evolution-step}  
Motivated by the  eigenvalue calculations shown in Section \ref{benchmark-step},
we consider the three  distinct values of the viscosity parameter 
$\nu=1/210$ (linearly stable),   $\nu=1/220$ (close to critical) and  
$\nu=1/250$ (unstable). 

We model the laboratory  scenario computationally via a two-stage process.
\begin{itemize} 
 \item[1.]
We start from a quiescent state and  a tiny time step ({\tt 1e-9}).
The inflow profile is smoothly ramped up to a fully developed flow using 
an exponential startup. 
The  sTR2 integration is then carried out  for  330 time steps with a relatively tight accuracy tolerance 
(i.e., a bound on an estimate of local truncation error), {\tt 3e-5}.
The number  of steps taken is arbitrary but needs to be chosen large enough so that the  
reference flow is visually steady.  
More precisely,  when this phase is complete, the instantaneous {\em acceleration} $a(t)$,
defined  in terms of  the  flow velocity $\uh(\cdot,t)$  at time $t$ by
$a(t) = \sqrt{\int_\calD \left( \frac{\partial \vec u_h}{\partial t}\right)^2}\,,$
should be around $10^{-2}$ or even less.

At the point  in time, $T$ say, where the first stage is completed, the integration is interrupted 
and a perturbation is added to the flow field $\vec{u}_h(\cdot,T)$.
The perturbation is of the form $\deltah$ specified in \rblack{Claim} 
2.1 where the associated scalar field $\phi_h$ derives from 
(\ref{truncated-discrete-KL})--(\ref{covariance-function}).
\rblack{We construct $\phi_h$ with $\sigma=.3(h^2)$}.\footnote{It is also 
necessary to  scale the perturbation so as not to ``shock'' the transient simulation --- the
perturbation field is thus scaled by a factor of  {\tt 1e-5}.
This ensures  that the magnitude of the perturbation is
comparable to the time accuracy used for the simulation.}
\item[2.] 
The time integration is then restarted without reducing the time step, using sTR1 
in place  of sTR2 (because it is marginally less dissipative). 
The restarted integration is continued for a fixed number (typically, 200 or 700) time 
steps, stopping prematurely only  if the time reaches  $T^*$={\tt 1e14} --- which we 
interpret  as reaching a ``computational steady-state" --- at which point 
\rblack{the adaptive time-stepping routine is taking very large time steps, see Figure 
\ref{timestep-histories-210-220}  below, and} the acceleration
 $a(T^*)$ will almost certainly be smaller than  unit roundoff.
\end{itemize}

The   {\it mean vorticity}  $\omega(t)$, or the average vertical velocity at the outflow,
 $$ 
 \omega(t) =  \int_\calD \nabla \times \uh(\cdot,t)  
 =  \int_{\partial\calD_N} \!\!\!\! u_y(\cdot,t) \, \rm{d}s  ,
 $$
provides a convenient way of assessing the degree of departure from 
the  reflectionally symmetric flow (for which $\omega$=0). 
At the conclusion of the first phase of the time integration, a pseudo-steady
flow is computed for each of the three values of $\nu$. 
The evolution of the mean vorticity  and the acceleration 
visualized  in Figure \ref{first-phase-250} 
shows that a symmetric flow is  established  for $\nu=1/250$ 
before the interruption is made after 330 time steps;  this corresponds to time $T \approx 32$.

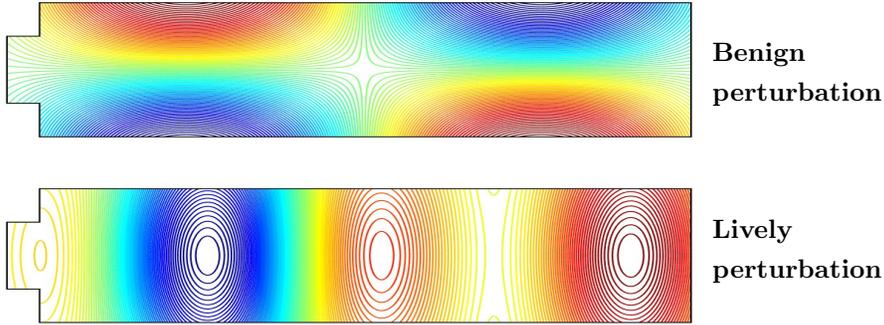
\begin{figure}
\begin{center}
\begin{picture}(260,155) (70,0) 
\put(0,0){\includegraphics[width=.9\textwidth, height=.4\textwidth]{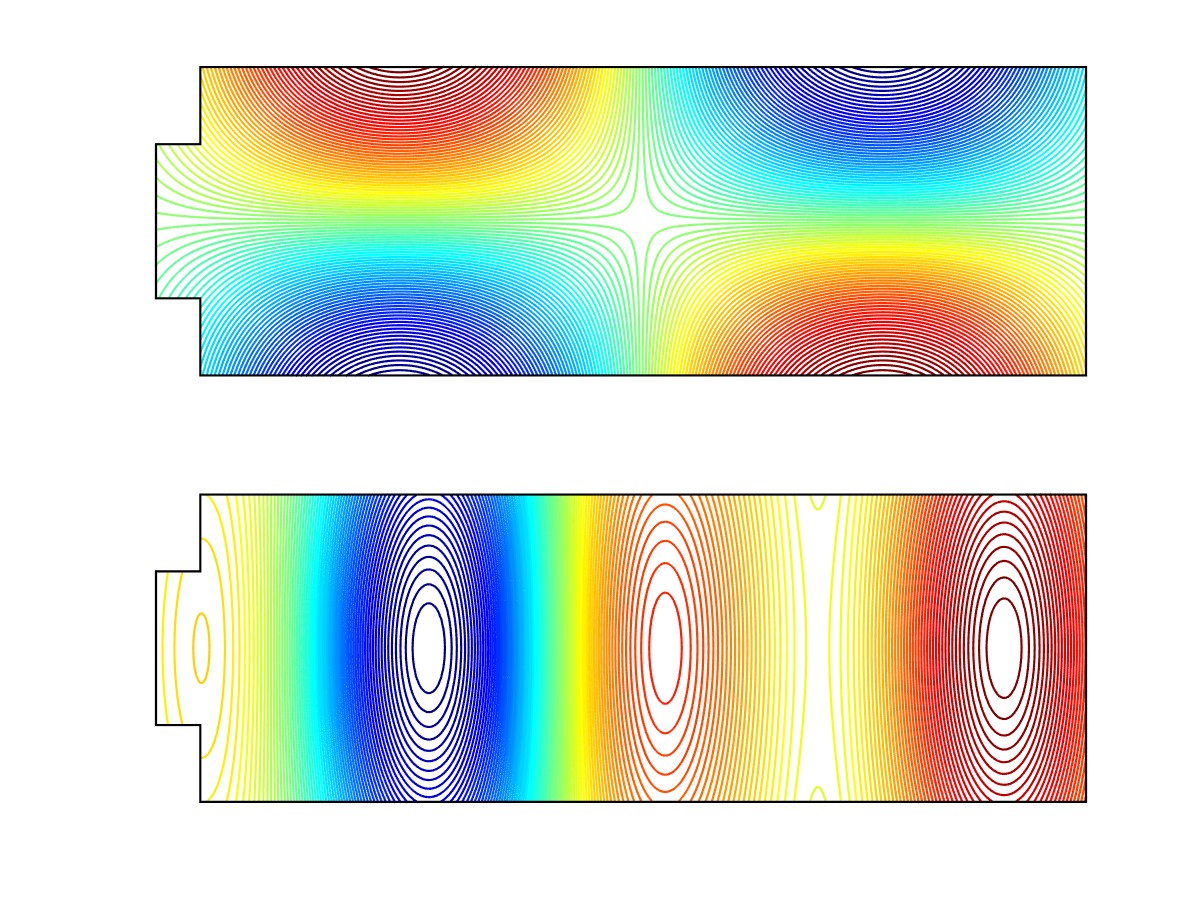}}
\put(310,115){\bf Benign} \put(310,100) {\bf perturbation}
\put(310,48){\bf Lively} \put(310,33) {\bf perturbation}
\end{picture}
\end{center}
\vspace{-.1in}
\caption{Structure of two scalar potentials $\phi_h$ used to generate velocity perturbations $\vfld{\delta}_h$}.
\label{phi-pictures}
\end{figure}

Moving on to the second stage,  we  show results 
for the subcritical cases $\nu=1/210$ and $1/220$, for
three representative flow perturbations, each of which derives from a particular
collocation point 
$\vxi^{(k)}$ used in (\ref{stability-indicator-surrogate}).
For the first of these, no perturbation is made  
(this corresponds to the point $\vxi\equiv \bv{0}$)
and the integration simply continues from the first-stage stopping point $T$. 
The other two are representatives of a ``benign'' perturbation and a ``lively''
perturbation and have the spatial structure shown  in Figure \ref{phi-pictures}.

The evolution of the flow after the restart for $\nu=1/210$  is
depicted in Figure~\ref{second-phase-210},
\rblack{where mean vorticity is shown in both actual and logarithmic scales.}
The unperturbed flow is perfectly stable; 
the sTR1 integrator reaches the end time ($T^*$={\tt 1e14}) \rblack{at time step 396,}
66 steps after the restart.
The distinctive jumps in the acceleration are associated with the stabilization
of the integrator, which has the effect of periodically injecting a small amount
of dissipation into the flow.  
The benign perturbation, which respects the reflectional symmetry, has no effect
on the long-term flow evolution. 
In contrast, the lively perturbation excites visible instability at about
time step 390, 60 steps after the restart. 
But (as seen in particular from the acceleration),
the size of the perturbation is not big 
enough to stop the long-term evolution to the symmetric
flow at the designated end time  $T^*$.
\rblack{The growth in vorticity for the stable examples toward the end of the
simulation (about step 385) is a roundoff effect caused by allowing the simulation
to proceed after changes in the steady solution are near machine precision.}

\begin{figure}
\begin{center}
\begin{picture}(270,360) (78,50) 
\put(198,153){\includegraphics[width=.61\textwidth, height=.88\textwidth]{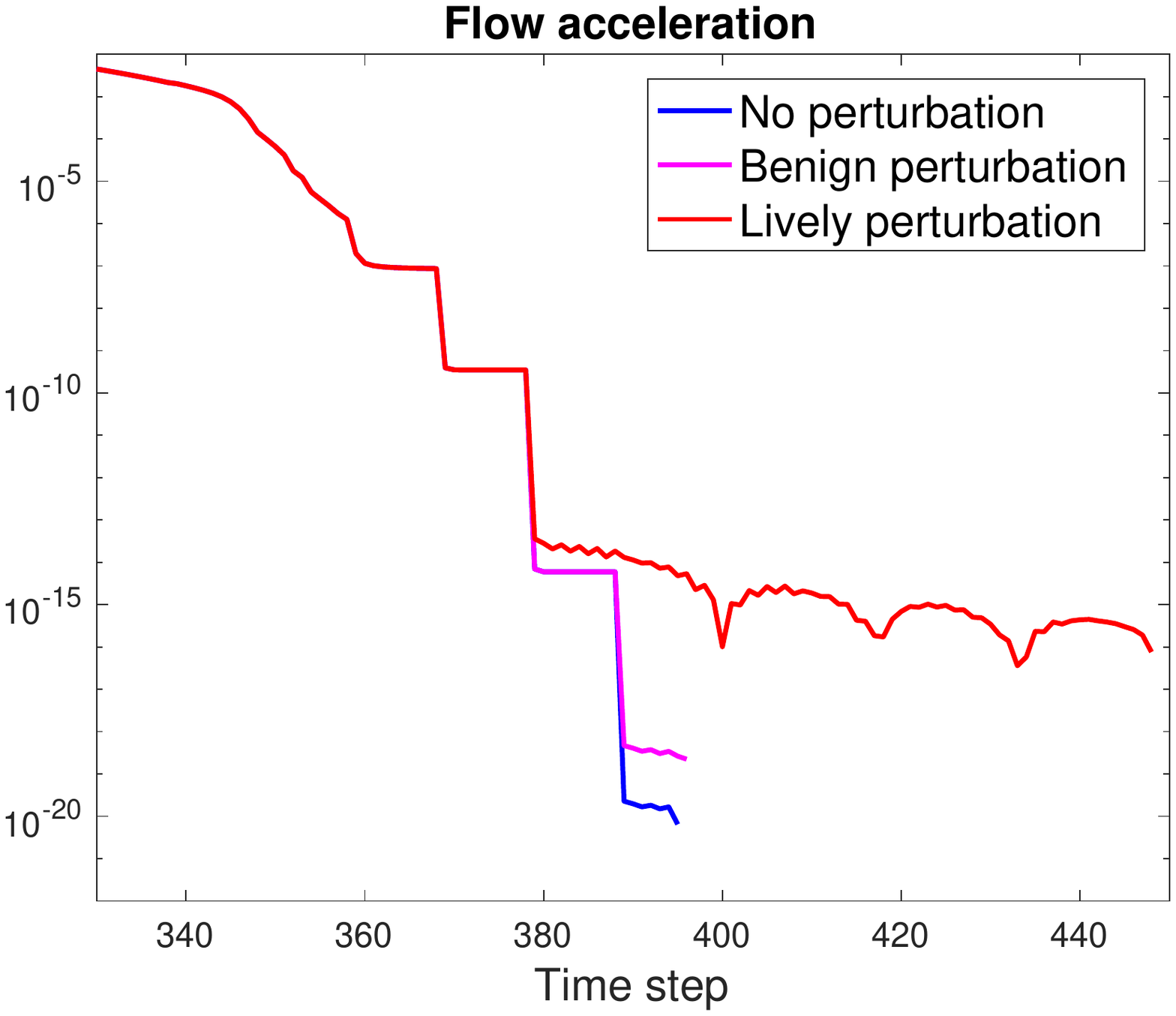}}
\put(4,153){\includegraphics[width=.61\textwidth, height=.88\textwidth]{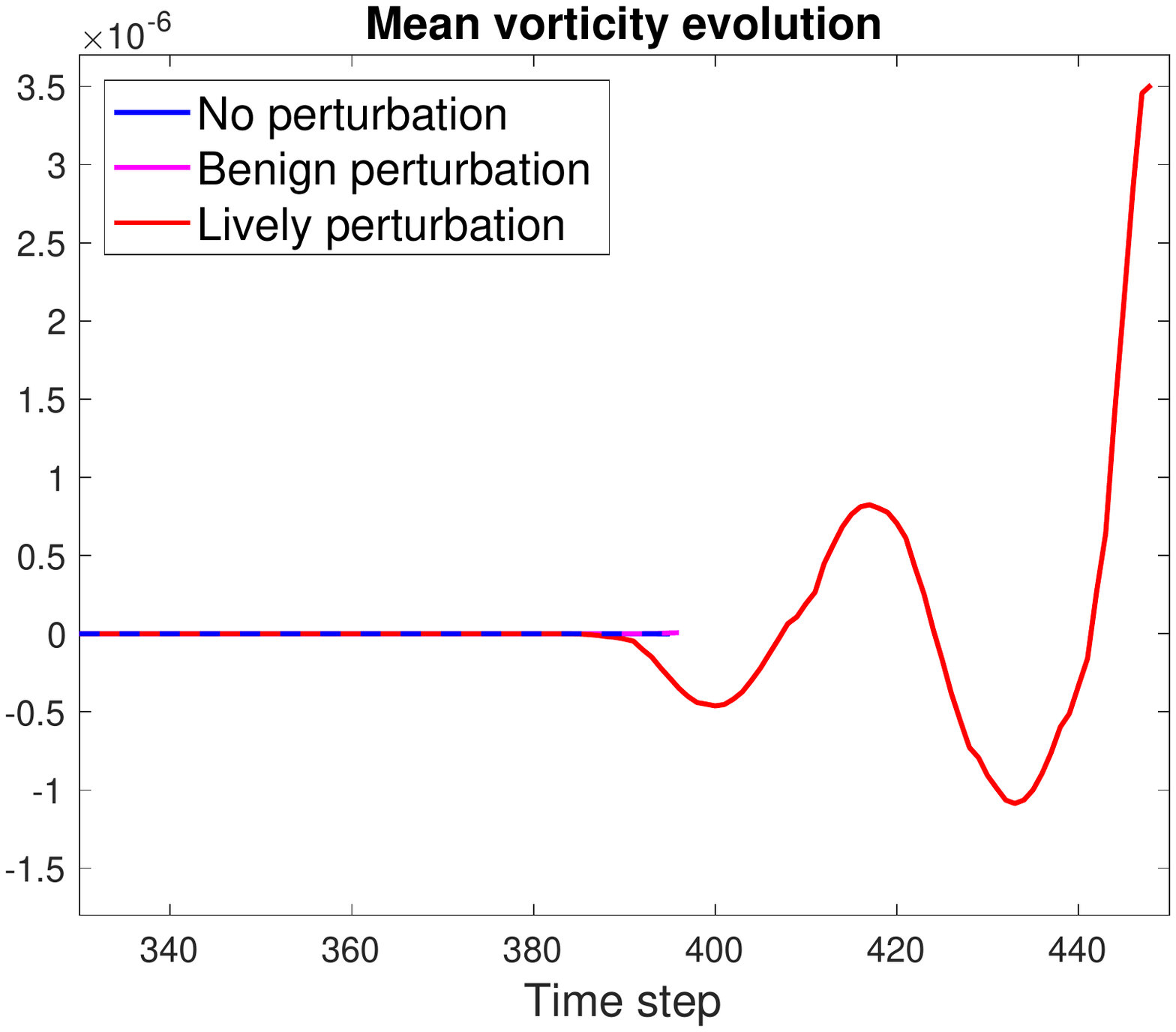}}
\put(102,-30){\includegraphics[width=.61\textwidth, height=.88\textwidth]{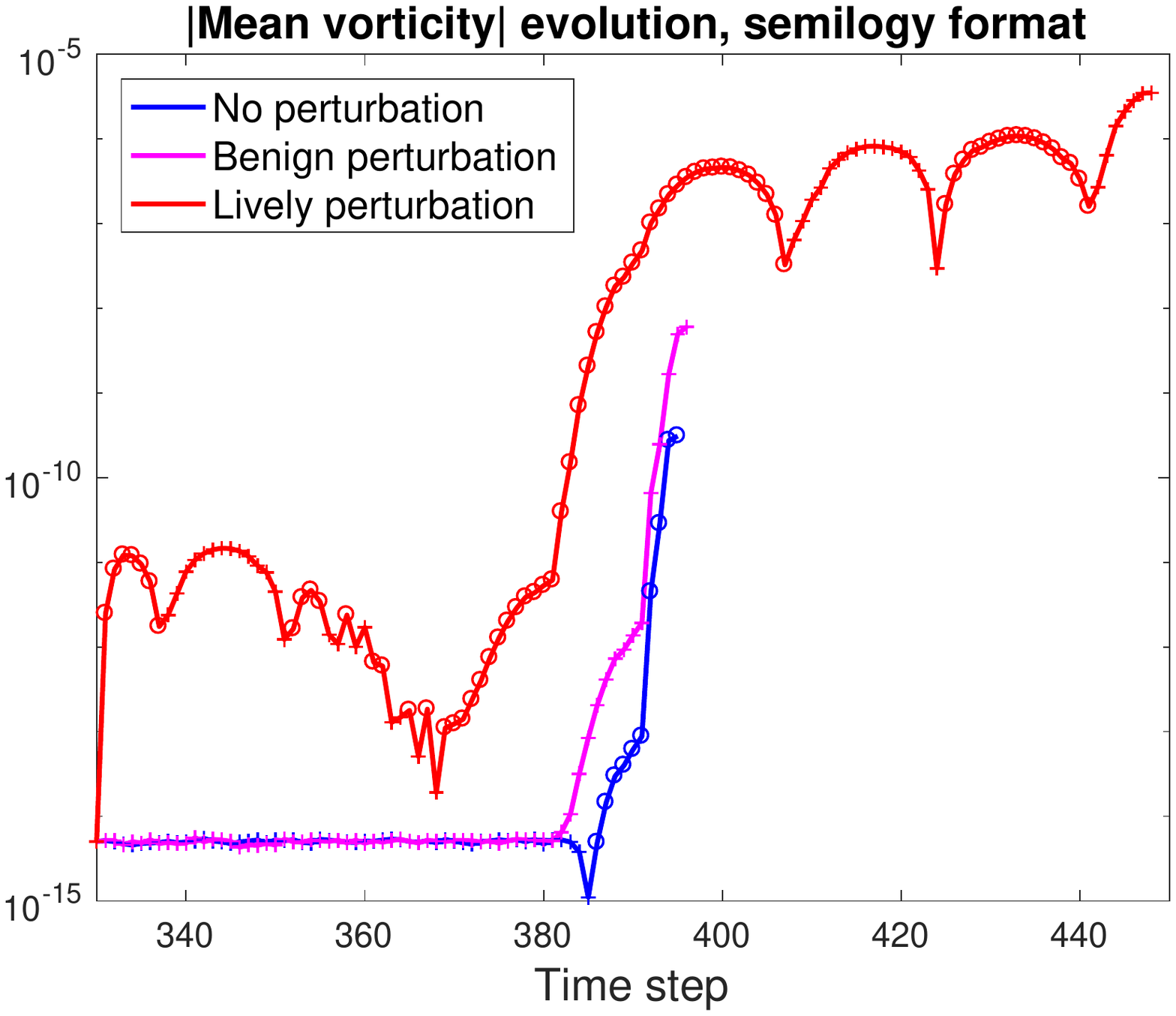}}
\end{picture}
\end{center}
\caption{Long-term evolution for different perturbations, symmetric step flow with $\nu=1/210$.} 
\label{second-phase-210}
\end{figure}

\begin{figure}
\begin{center}
\begin{picture}(270,360) (78,50) 
\put(198,153){\includegraphics[width=.61\textwidth, height=.88\textwidth]{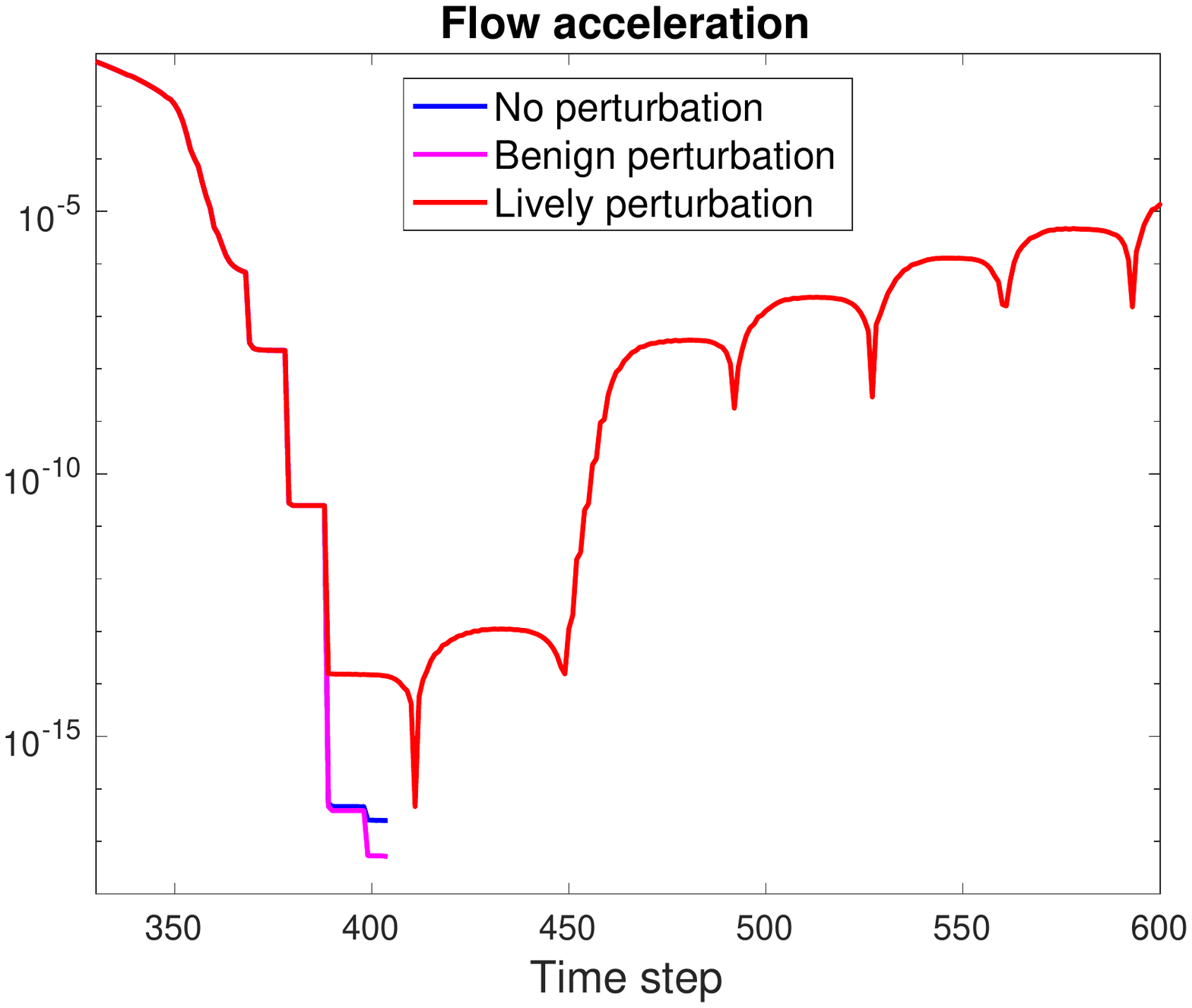}}
\put(4,153){\includegraphics[width=.61\textwidth, height=.88\textwidth]{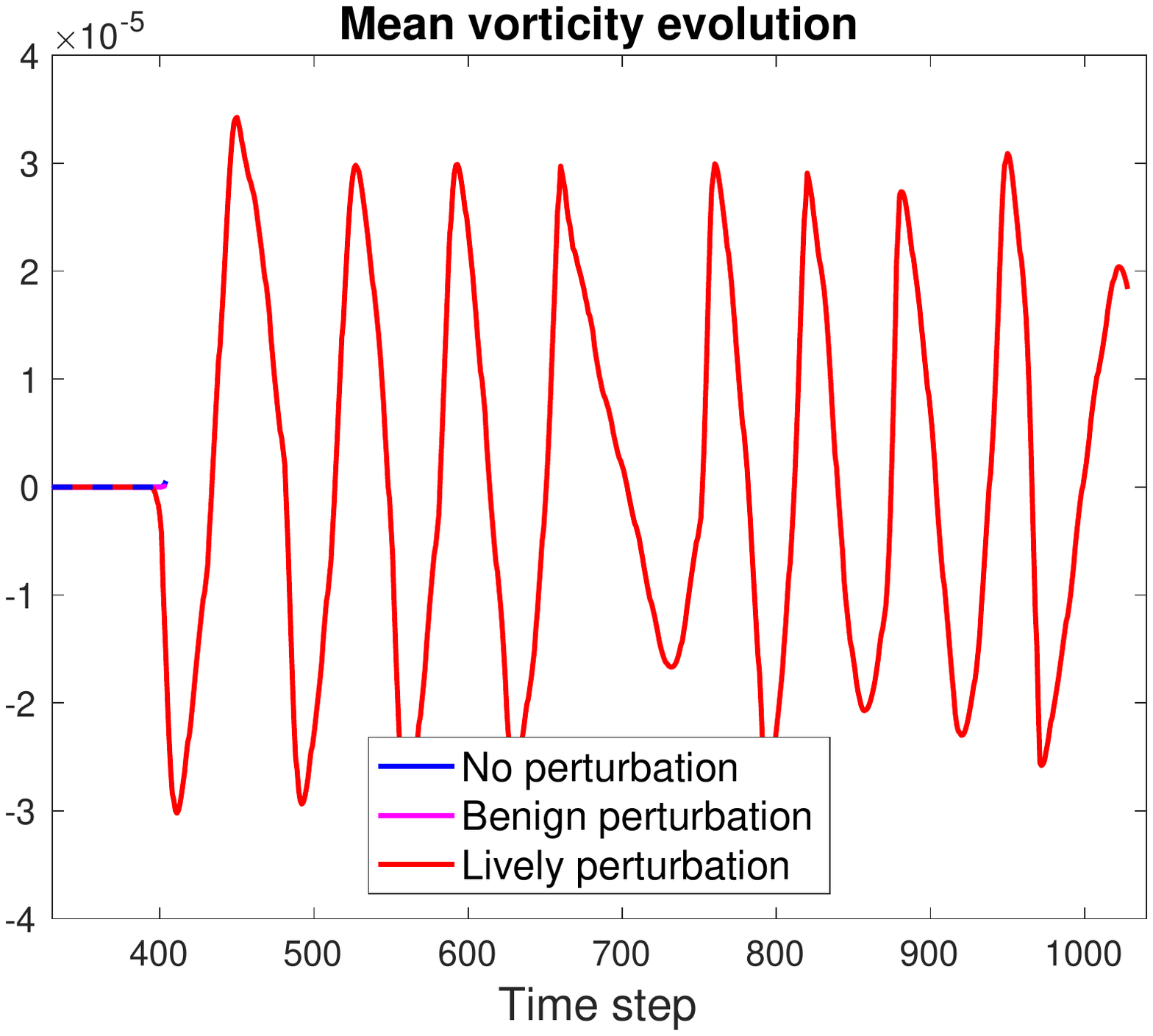}}
\put(102,-30){\includegraphics[width=.61\textwidth, height=.88\textwidth]{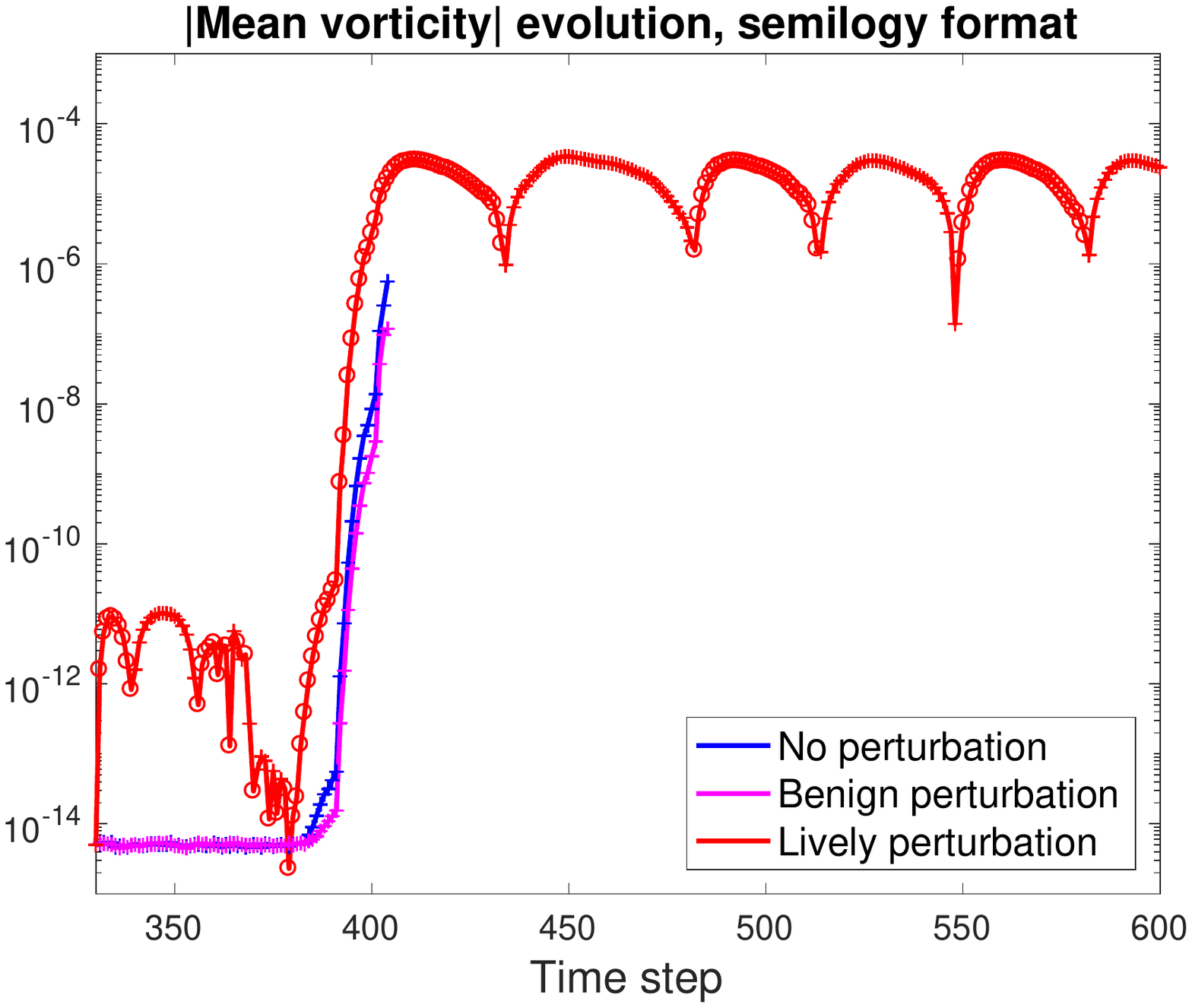}}
\end{picture}
\end{center}
\caption{Long-term evolution for different perturbations, symmetric step flow with $\nu=1/220$.} 
\label{second-phase-220}
\end{figure}

\begin{figure}
\begin{center}
\begin{picture}(270,620) (62,20) 
\put(0,370){\includegraphics[width=.6\textwidth, height=.90\textwidth]{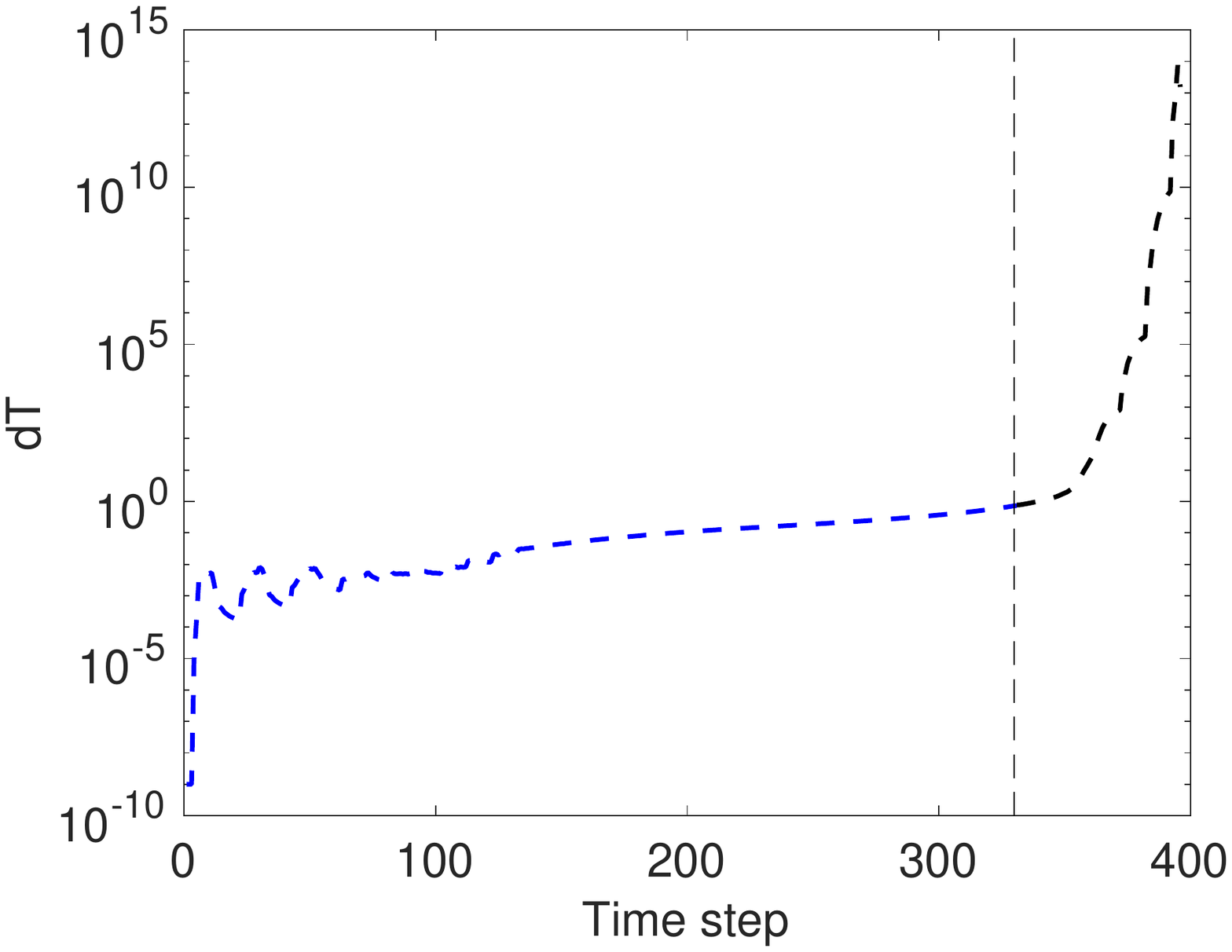}}
\put(190,370){\includegraphics[width=.6\textwidth, height=.90\textwidth]{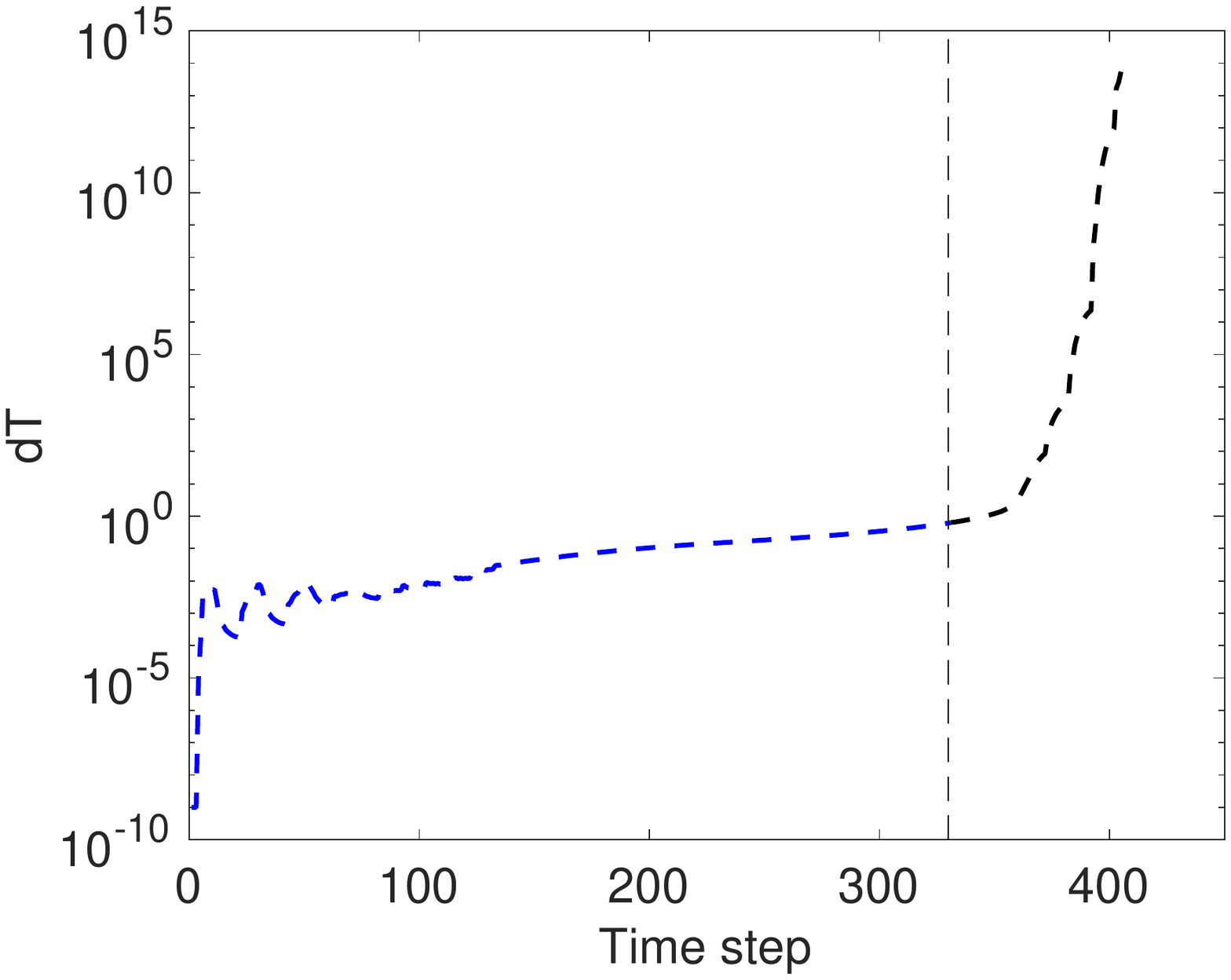}}
\put(170,440){\bf No perturbation}
\put(0,170){\includegraphics[width=.6\textwidth, height=.90\textwidth]{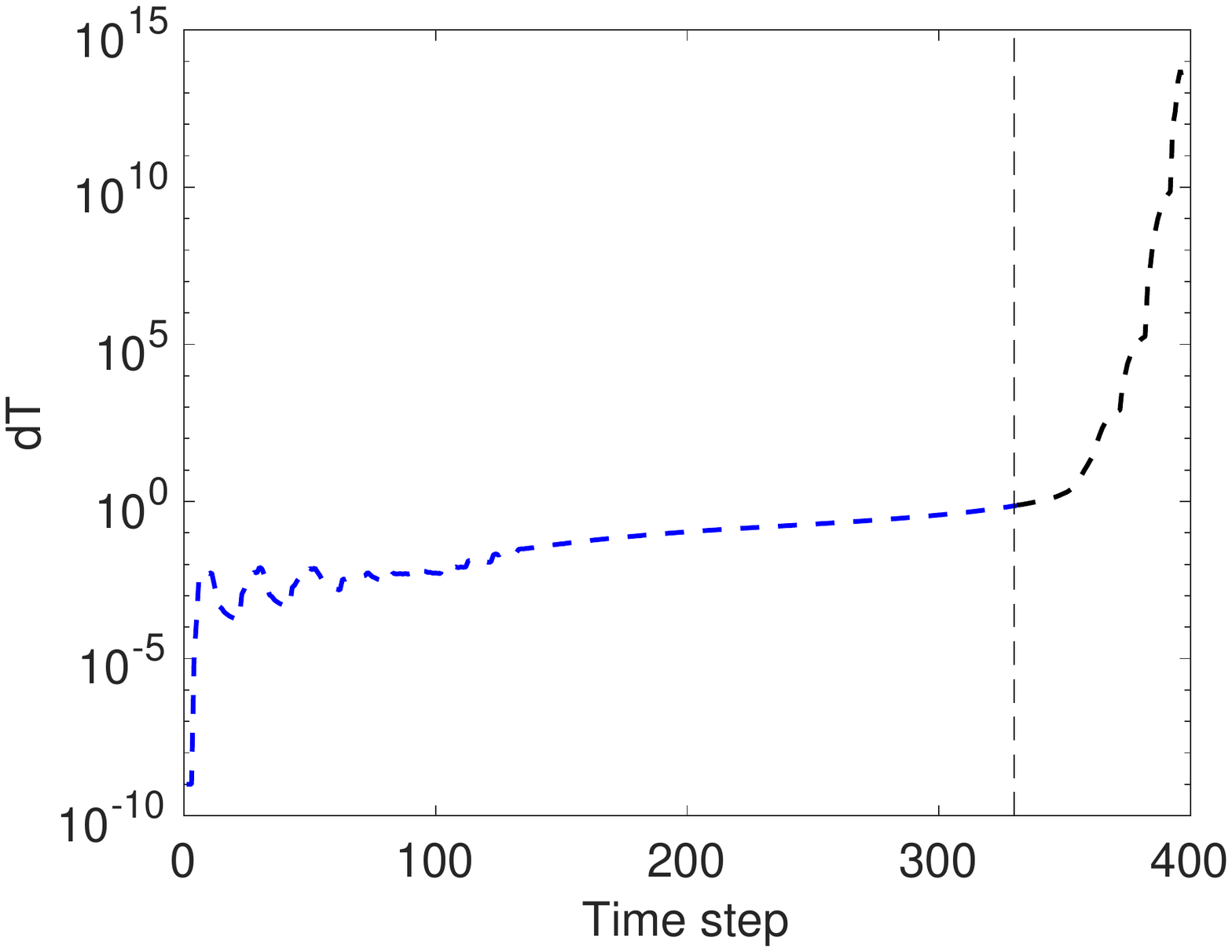}}
\put(190,170){\includegraphics[width=.6\textwidth, height=.90\textwidth]{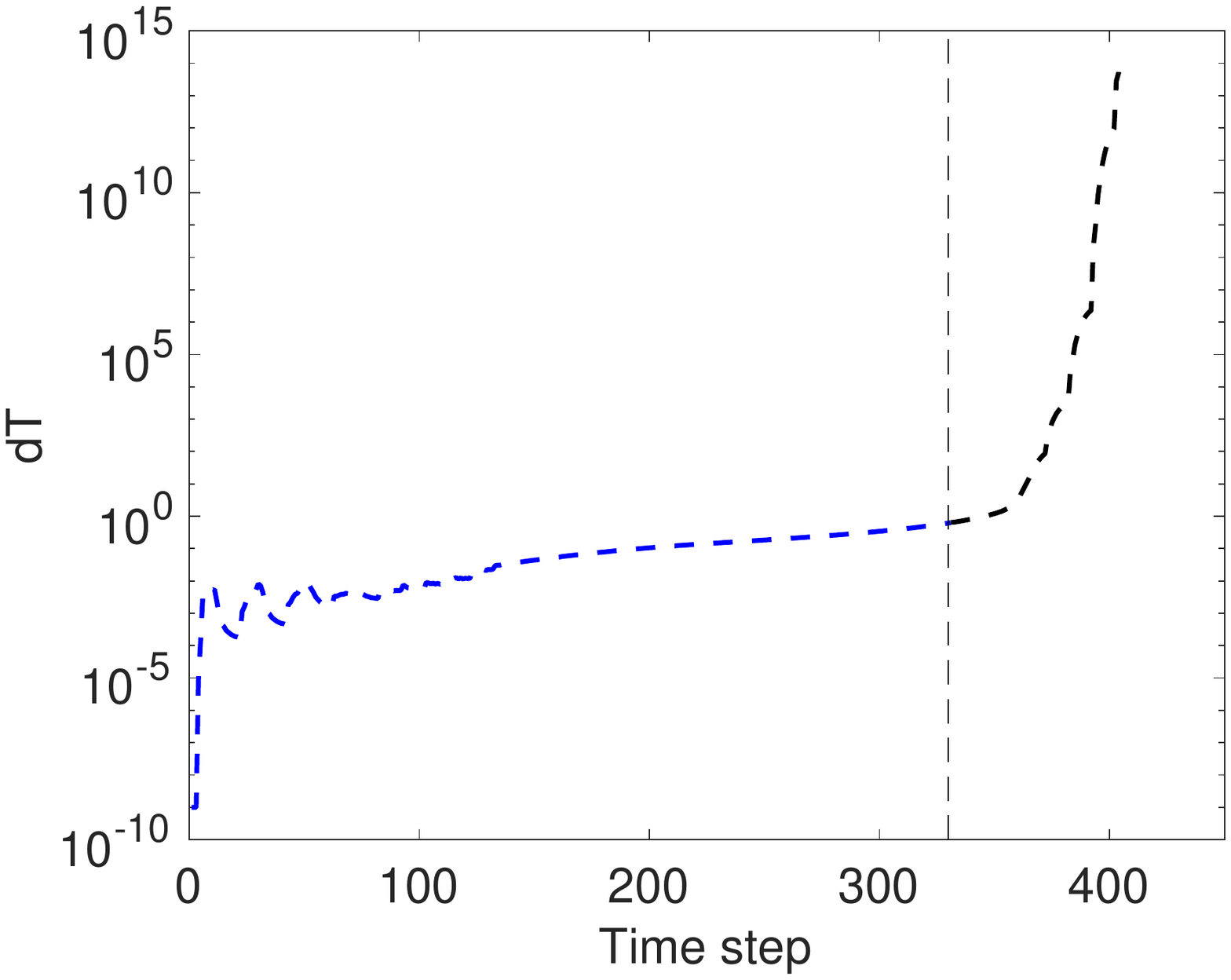}}
\put(165,240){\bf Benign perturbation}
\put(0,-30){\includegraphics[width=.6\textwidth, height=.90\textwidth]{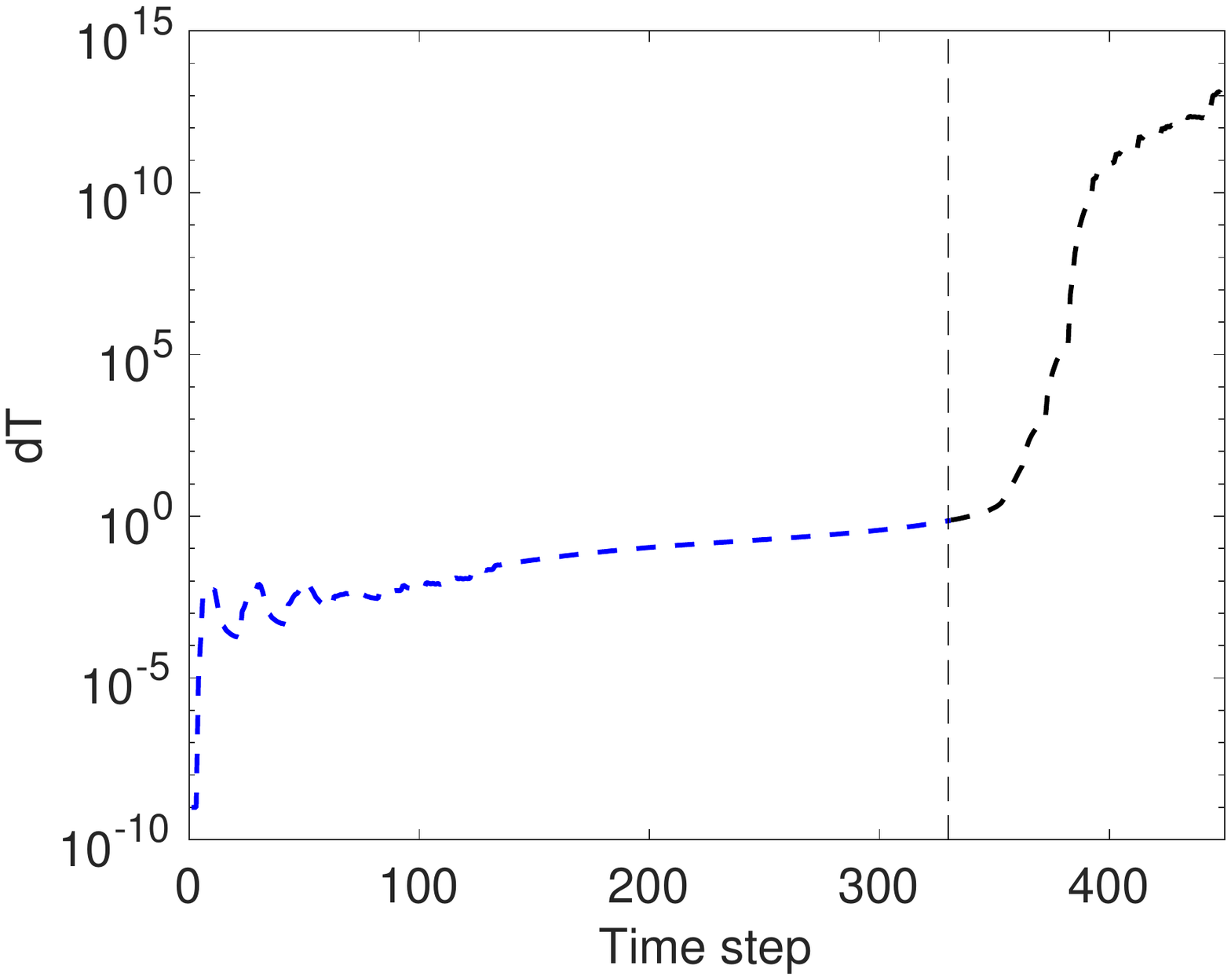}}
\put(190,-30){\includegraphics[width=.6\textwidth, height=.90\textwidth]{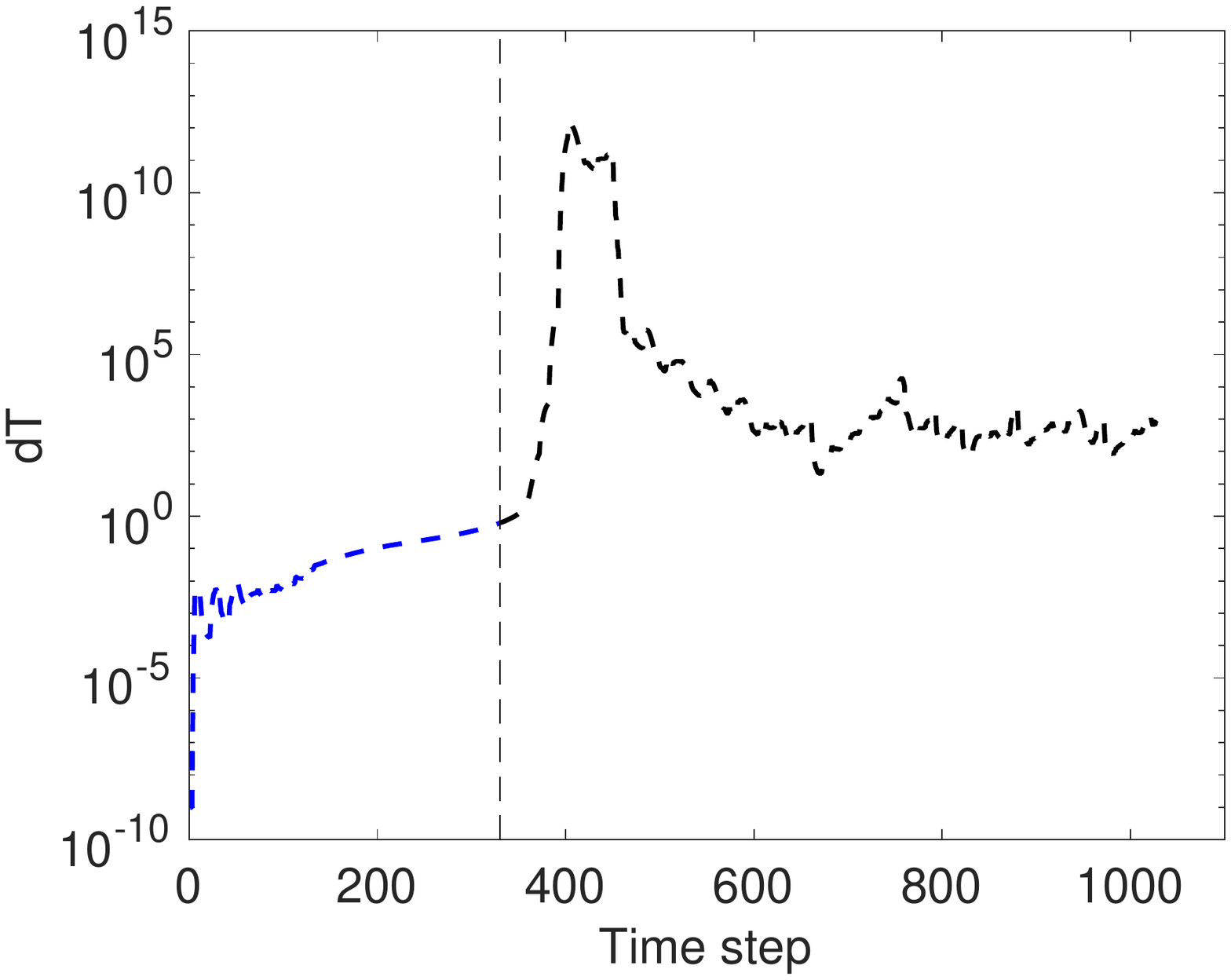}}
\put(165,40){\bf Lively perturbation}
\end{picture}
\end{center}
\caption{Time-step histories: time step size vs.\ time step count for various perturbations, 
symmetric step, $\nu=1/210$ (left), $\nu=1/220$ (right).} 
\label{timestep-histories-210-220}
\end{figure}

The  evolution of flow for the intermediate viscosity  parameter  
$\nu$ = 1/220 is shown in Figure \ref{second-phase-220}.
The unperturbed case is just about stable: 
the sTR1 integrator reaches the end time  75 steps after the restart. 
A virtually identical evolution is evident when the perturbation to the flow
is benign.
\rblack{(As in the previous example, roundoff effects lead to some growth in
the vorticity after a steady solution is obtained.)} 
The time evolution for the lively perturbation is noticeably different,
however. 
In this case, the sTR1 integrator rejects time step 415 (85 steps after the
restart) and the computational flow evolves to a 
numerically noisy solution where the magnitude of the 
oscillation is of the order of the time-stepping accuracy.

These observations are substantiated in Figure~\ref{timestep-histories-210-220}, which shows the 
history of the time step sizes chosen by the adaptive integrator.
For each of the plots in this figure, the switch from the first to the second stage
is identified by a vertical dotted line.
When either no perturbation or a benign perturbation is made, the time step sizes 
rapidly increase because the integration goes to a steady state for the subcritical 
values of $\nu$.
This behavior can also be seen for the lively perturbation and $\nu=1/210$.
In contrast,  the  integrator behaves differently  for $\nu=1/220$ --- here
the time step size is cut back at around 70 time steps after the 
perturbation is made in order  to resolve the 
nonstationary solution shown at the bottom of Figure~\ref{second-phase-220}.
Computing solutions  when so close to the stability limit is  a
delicate business.

\begin{figure}
\begin{center}
\begin{picture}(270,350) (42,0) 
\put(70,90){\includegraphics[width=.6\textwidth, height=.90\textwidth]{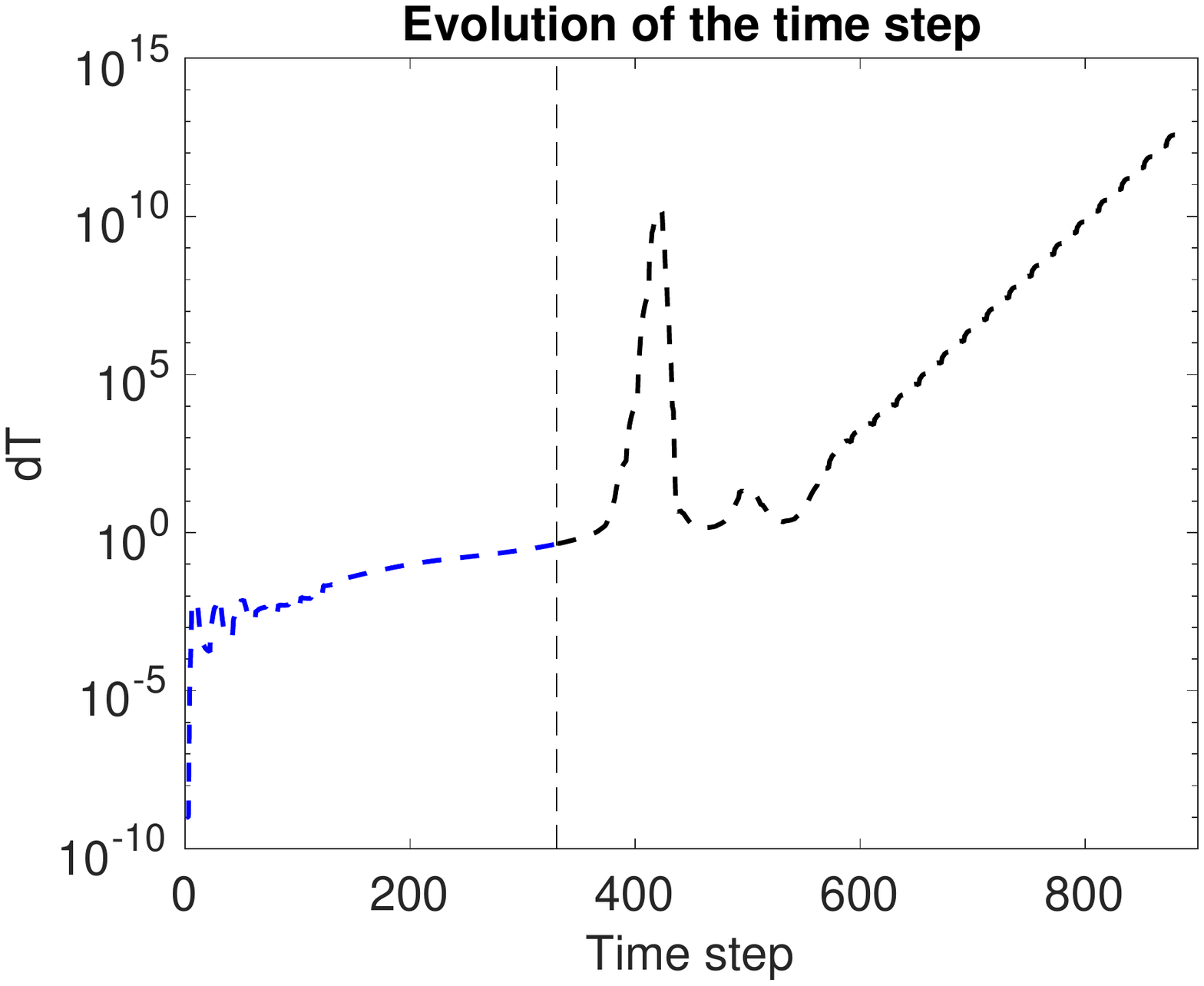}}
\put(0,-5){\includegraphics[width=.95\textwidth, height=.45\textwidth]{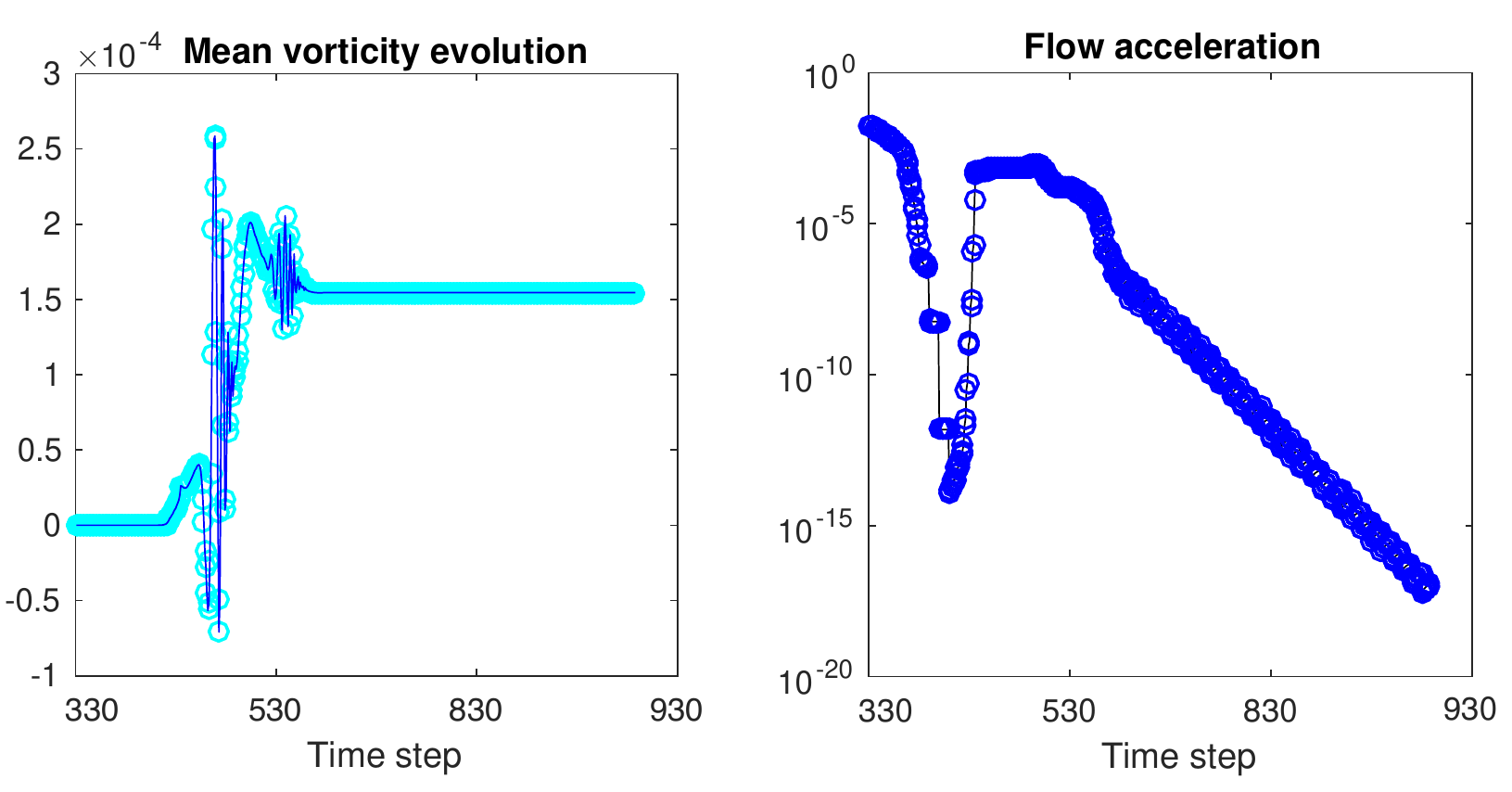}}
\end{picture}
\end{center}
\caption{Time-step history and long-term evolution for the unperturbed flow,
  symmetric step, $\nu=1/250$.}
\label{results-250}
\end{figure}

\begin{figure}
\begin{center}
\begin{picture}(270,340) (42,0) 
\put(0,0){\includegraphics[width=.95\textwidth, height=.95\textwidth]{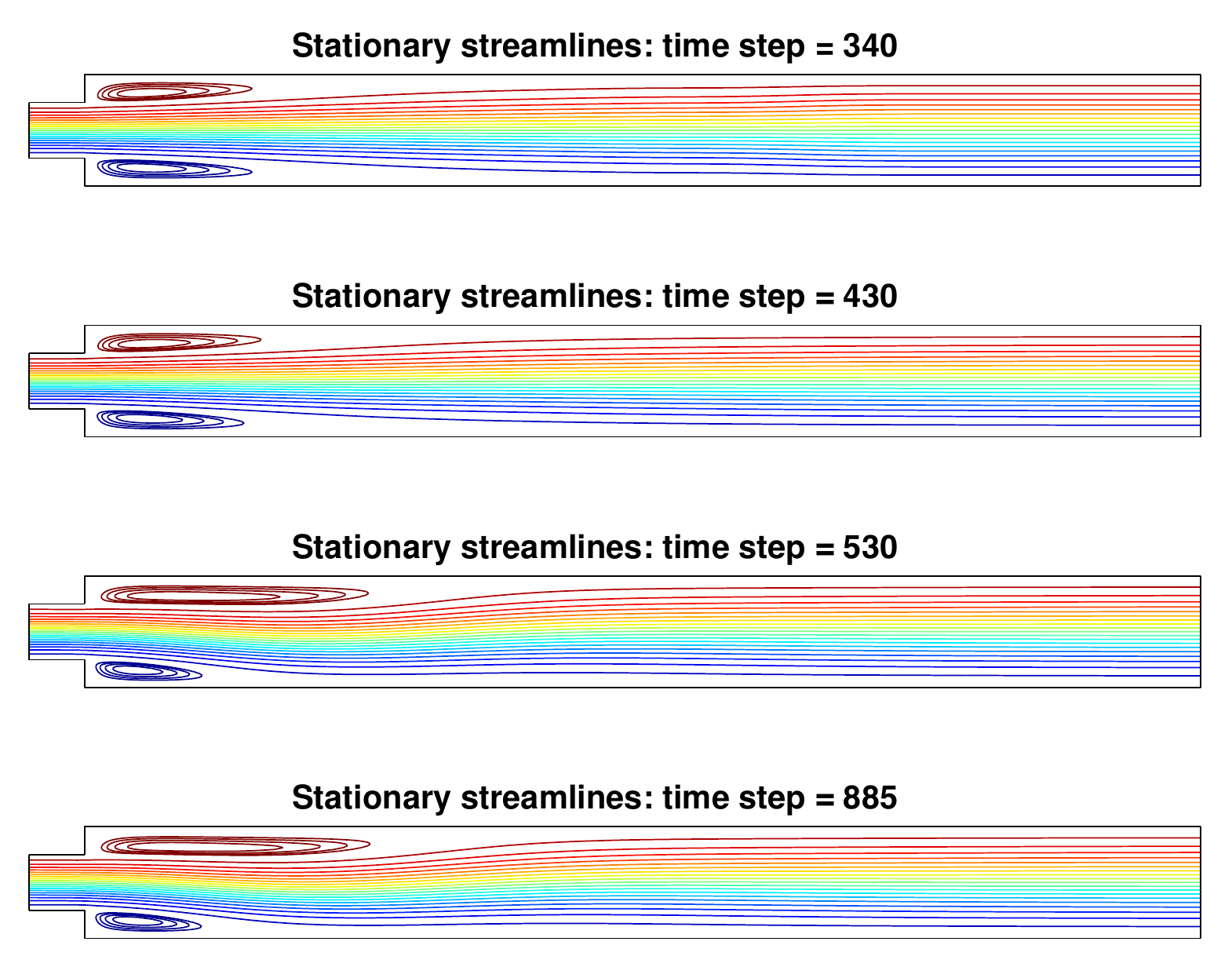}}
\end{picture}
\end{center}
\vspace{-.1in}
\caption{Snapshots of the (unperturbed) flow, symmetric step, $\nu=1/250$.}
\label{snaps-250}
\end{figure}

Results for the super-critical viscosity  parameter
$\nu$ = 1/250 are in Figure \ref{results-250}. 
In this instance, no perturbations are needed to excite instability.  
The time step history of the complete flow evolution
from $t$={\tt 0} to $t$={\tt 1e14} is presented  
at the top of the figure. 
Note the scale on the vertical axis --- this is  a pretty demanding
computational exercise!
The evolution  of mean vorticity and acceleration after the interrupt is shown in 
the two plots at the bottom of Figure \ref{results-250} and
should be contrasted with the results for the subcritical viscosity
shown in Figure \ref{second-phase-220}.
Just when the symmetric flow looks to be steady  (400 time steps;  70
after the restart) the time step is cut back to O(1) and after a transient
the flow goes to a computational steady state that does not have 
reflectional symmetry.  This is evident from the flow snapshots plotted at/after the
interrupt shown in Figure \ref{snaps-250}; 
the particular steady-state solution (top eddy longer than the bottom one) is solely 
determined by the build-up of roundoff error.
The two ``cups'' between 400 and 600 in the time step
history shown at the top of Figure \ref{results-250} suggest that the sTR1  
algorithm needed two
attempts to fix on the specific stationary solution --- it is instructive to contrast
this with the evolution that results  when vigorously perturbing the
flow close to the critical viscosity, which is shown at the bottom of 
Figure~\ref{timestep-histories-210-220}.

\subsection{Evolution of flow  around an obstacle} 
\label{evolution-obstacle}
Motivated by the  eigenvalue calculations  discussed in Section~\ref{benchmark-obstacle},
we now consider three distinct values of the viscosity parameter for the 
obstacle problem: 
$\nu=1/175$ (subcritical),
$\nu=1/185.6$ (close to critical) and  
$\nu=1/200$ (unstable). 
We consider $\nu=1/200$ first.
The same two-stage process described above gives the results shown in 
Figure~\ref{obstacle-results-200}.  
These results should be compared with those in Figure~\ref{results-250}. 
The difference is that instead of going to a nonsymmetric steady state solution,
 the computational flow evolves to a periodic (vortex-shedding) solution, at which point the 
 time step becomes essentially constant. 
The vortex-shedding solution is persistent --- it is unchanged when we
run the solver for another 10,000 time steps. 
The same long-term behavior is obtained if a perturbation is added  at the interrupt point. 
The different outcomes for the two benchmark problems are representative of the difference 
between a pitchfork bifurcation (for the step) and a Hopf bifurcation (for the obstacle) 
\cite{Cliffe-Spence-Taverner},\cite[p.\ 343]{ESW}.

\begin{figure}
\begin{center}
\begin{picture}(270,340) (45,10) 
\put(70,90){\includegraphics[width=.6\textwidth, height=.90\textwidth]{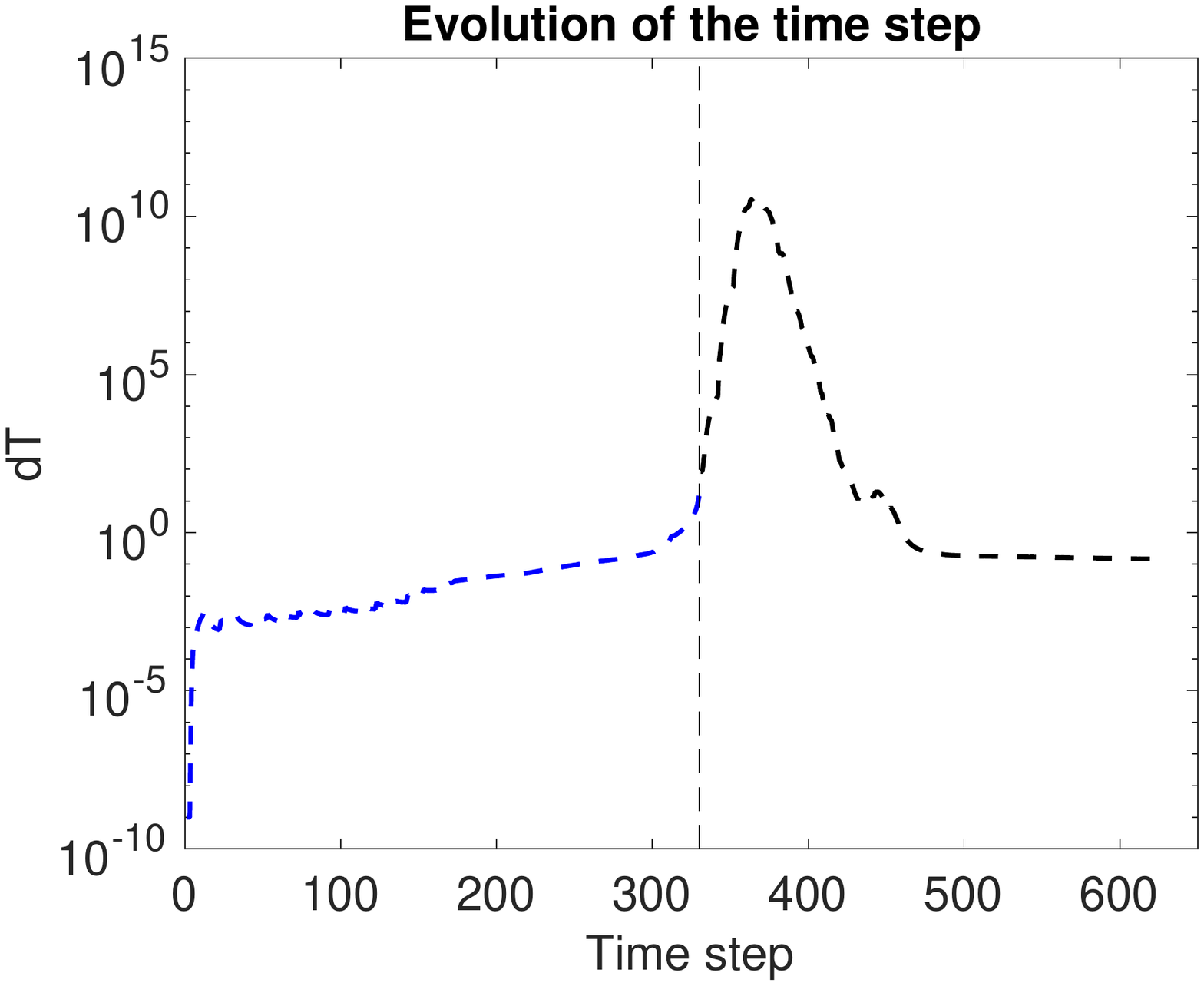}}
\put(0,-5){\includegraphics[width=.95\textwidth, height=.45\textwidth]{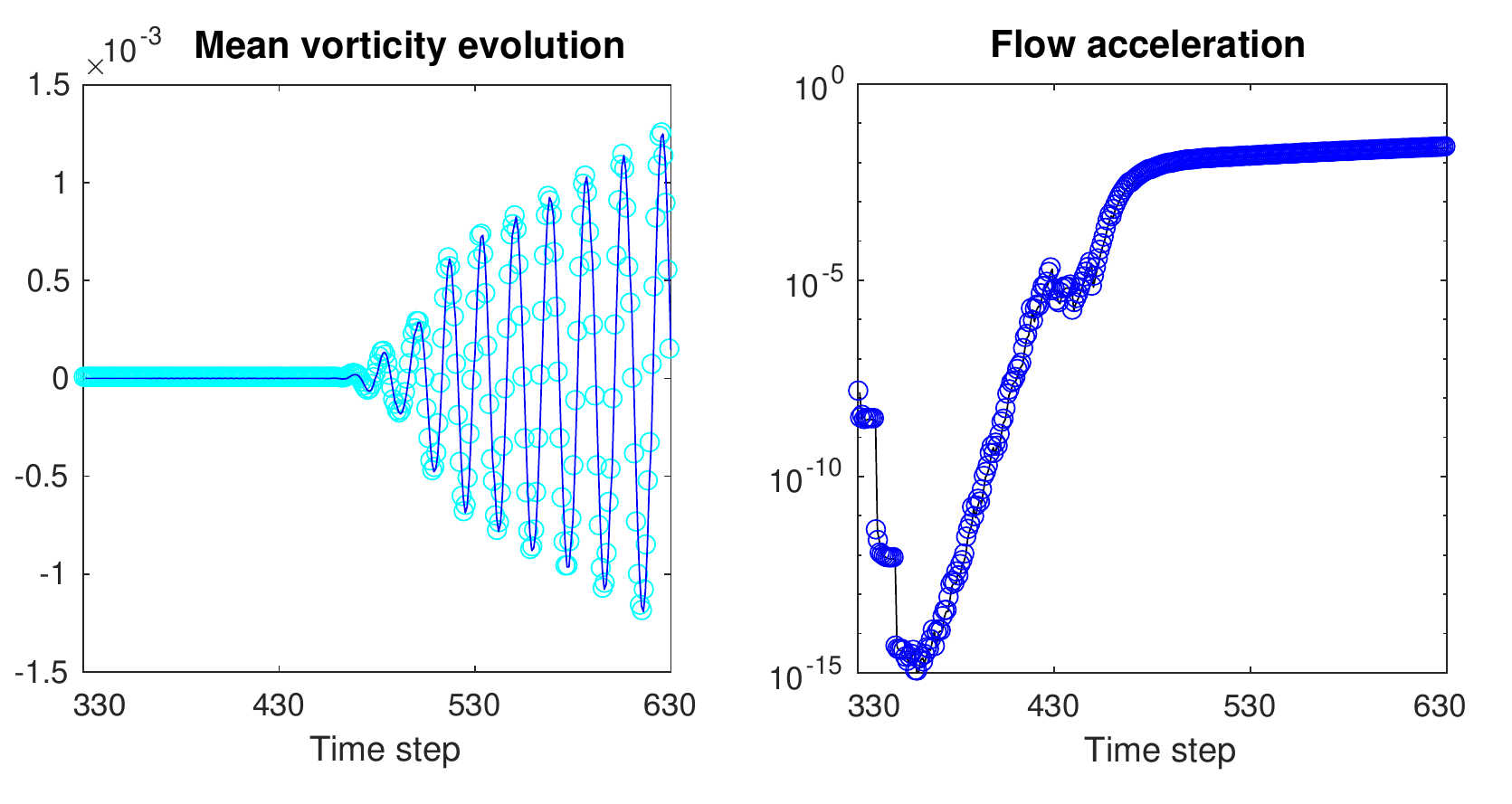}}
\end{picture}
\end{center}
\caption{Time-step history and long-term evolution for the unperturbed flow, obstacle, 
$\nu=1/200$.}
\label{obstacle-results-200}
\end{figure}

\begin{figure}
\begin{center}
\begin{picture}(270,570) (42,0) 
\put(7,390){\includegraphics[width=.93\textwidth, height=.50\textwidth]{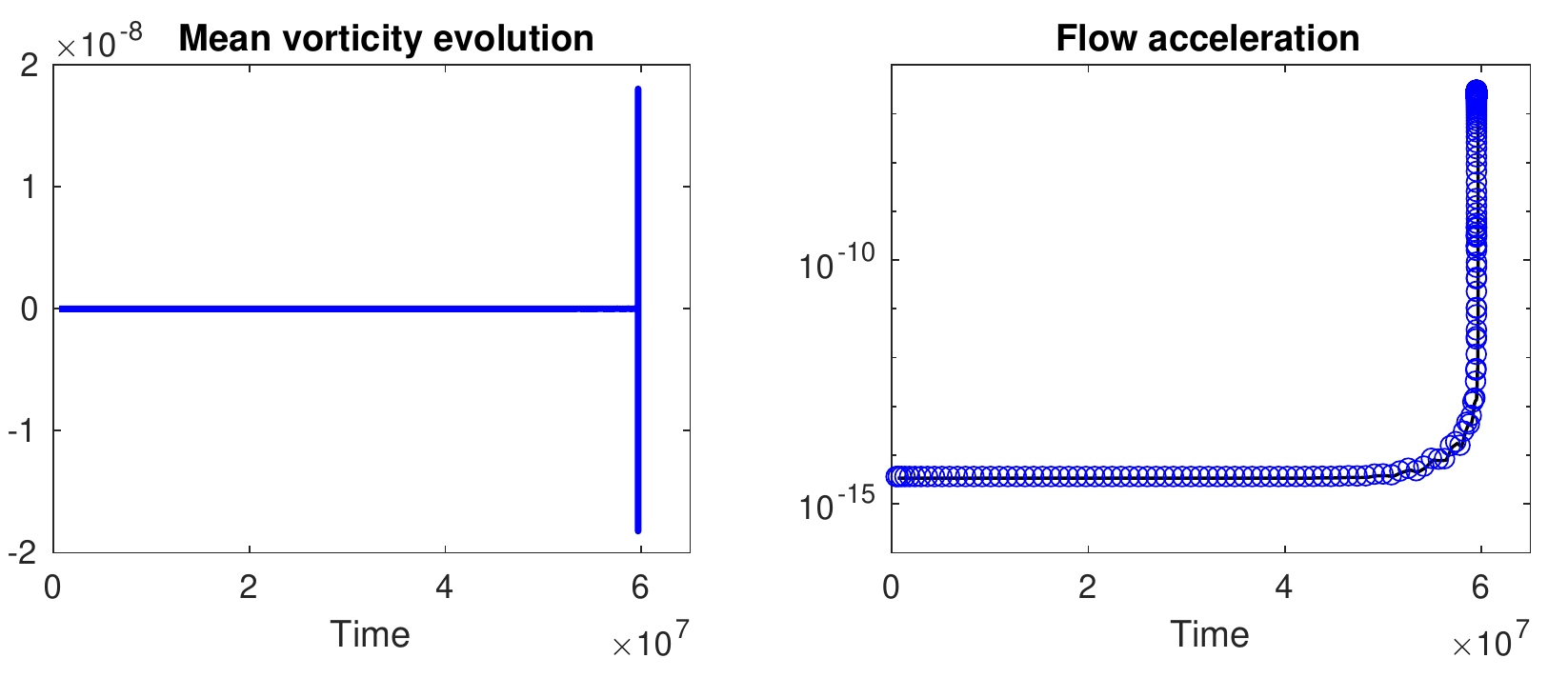}}
\put(110,385){\bf Flow evolution for $\nu=1/175$}
\put(2,200){\includegraphics[width=.95\textwidth, height=.45\textwidth]{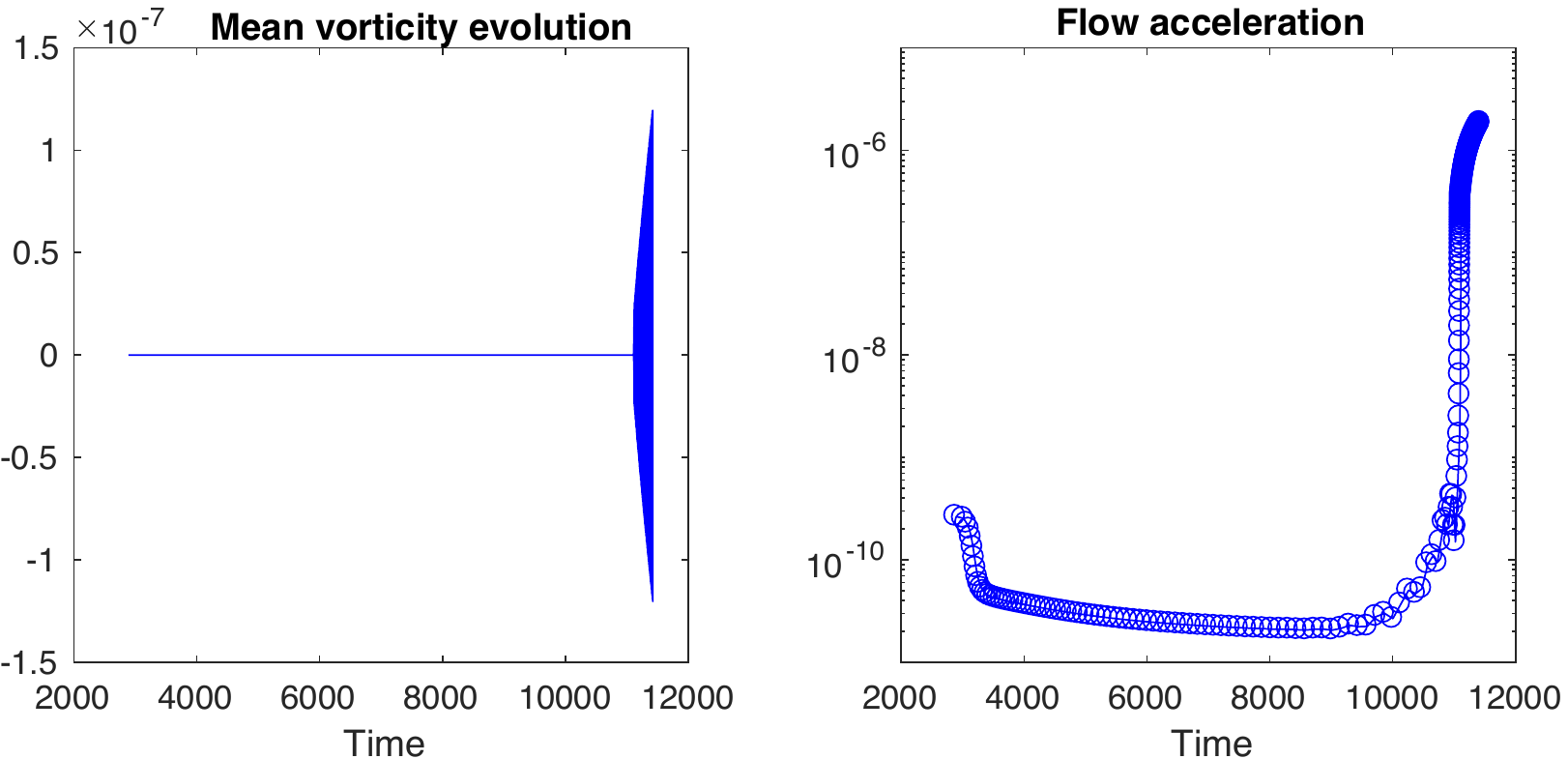}}
\put(110,188){\bf Flow evolution for $\nu=1/185.6$}
\put(5,8){\includegraphics[width=.95\textwidth, height=.45\textwidth]{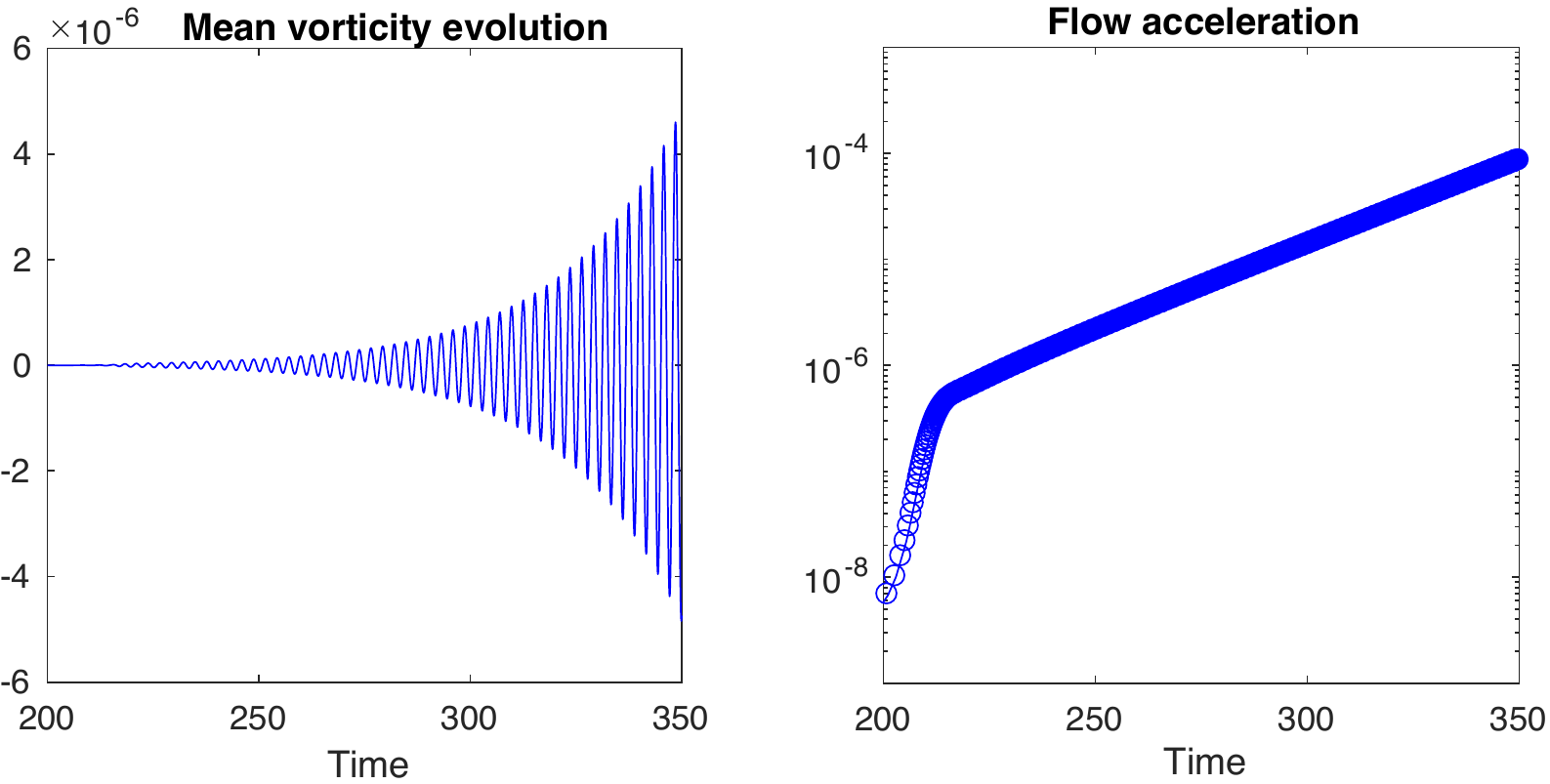}}
\put(110,0){\bf Flow evolution for  $\nu=1/200$}
\end{picture}
\end{center}
\caption{Evolution of mean vorticity and flow acceleration for three viscosity parameters, no 
perturbations, obstacle flow.} 
\label{fig-evolution-obstacle}
\end{figure}

To study the flow breakdown mechanism in detail,  the second phase of the time integration
is computed with a very small  accuracy tolerance  ({\tt 1e-9}) using the 
{\it unstabilised\/} TR1 integrator.\footnote{Stabilization of TR is not appropriate  when 
the accuracy tolerance is so small.}
In all cases discussed below the time integrator is run for 2500 
steps after the interrupt.
Figure~\ref{fig-evolution-obstacle} shows the evolution of the mean vorticity and the acceleration 
using this refined strategy for each value of the viscosity parameter, when no perturbation is done.
In the super-critical case of $\nu=1/200$ (bottom), there is a fast breakdown to the vortex-shedding
solution.
(Note that the evolution is plotted against physical time in this figure.)
For both the subcritical ($\nu=1/175$) and near-critical ($\nu=1/185.6$) cases, there are long delays
(until $t\approx 6e7$ and $t\approx 1.1e4$, respectively) after the interrupt, after which 
numerical instability kicks in and (as in the preceding section) generates a numerically 
noisy solution.
The onset of instability is dramatically later for the subcritical case.

We explore the breakdowns in more depth in Figure \ref{second-phase-obstacle-zoom}, 
which shows magnified images of the noisy solution measures at the time they become 
unsteady.
These images show that the magnitudes of the numerical oscillations (of order {\tt 1e-8} in the
sub-critical case and {\tt 1e-7} in the near-critical case) are comparable to the time-stepping accuracy.
Even when no explicit perturbation is done, time accuracy plays a role in long-term simulation to 
compute steady solutions in near-critical regimes.

\begin{figure}
\begin{center}
\begin{picture}(270,368) (42,0)
\put(6,193){\includegraphics[width=.95\textwidth, height=.47\textwidth]{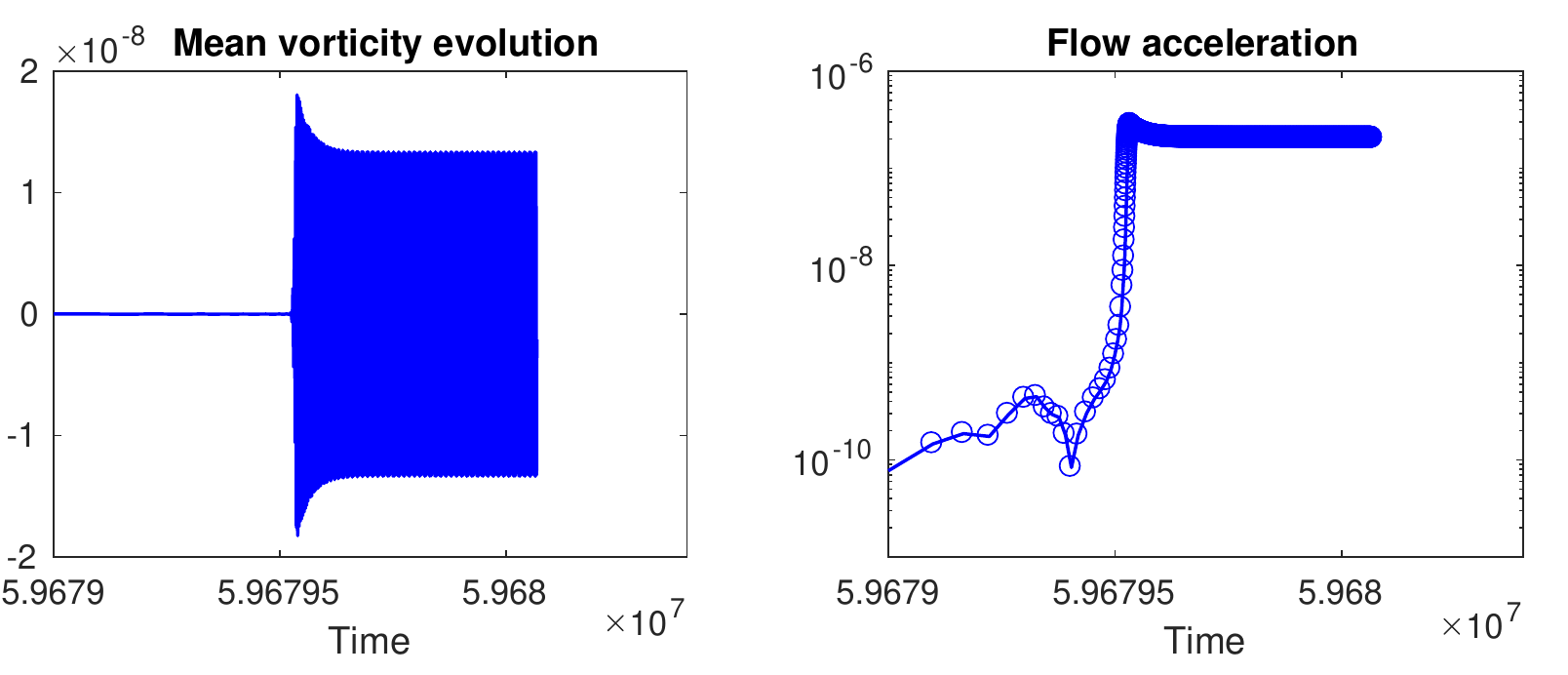}}
\put(160,192){$\nu=1/175$}
\put(2,10){\includegraphics[width=.96\textwidth, height=.45\textwidth]{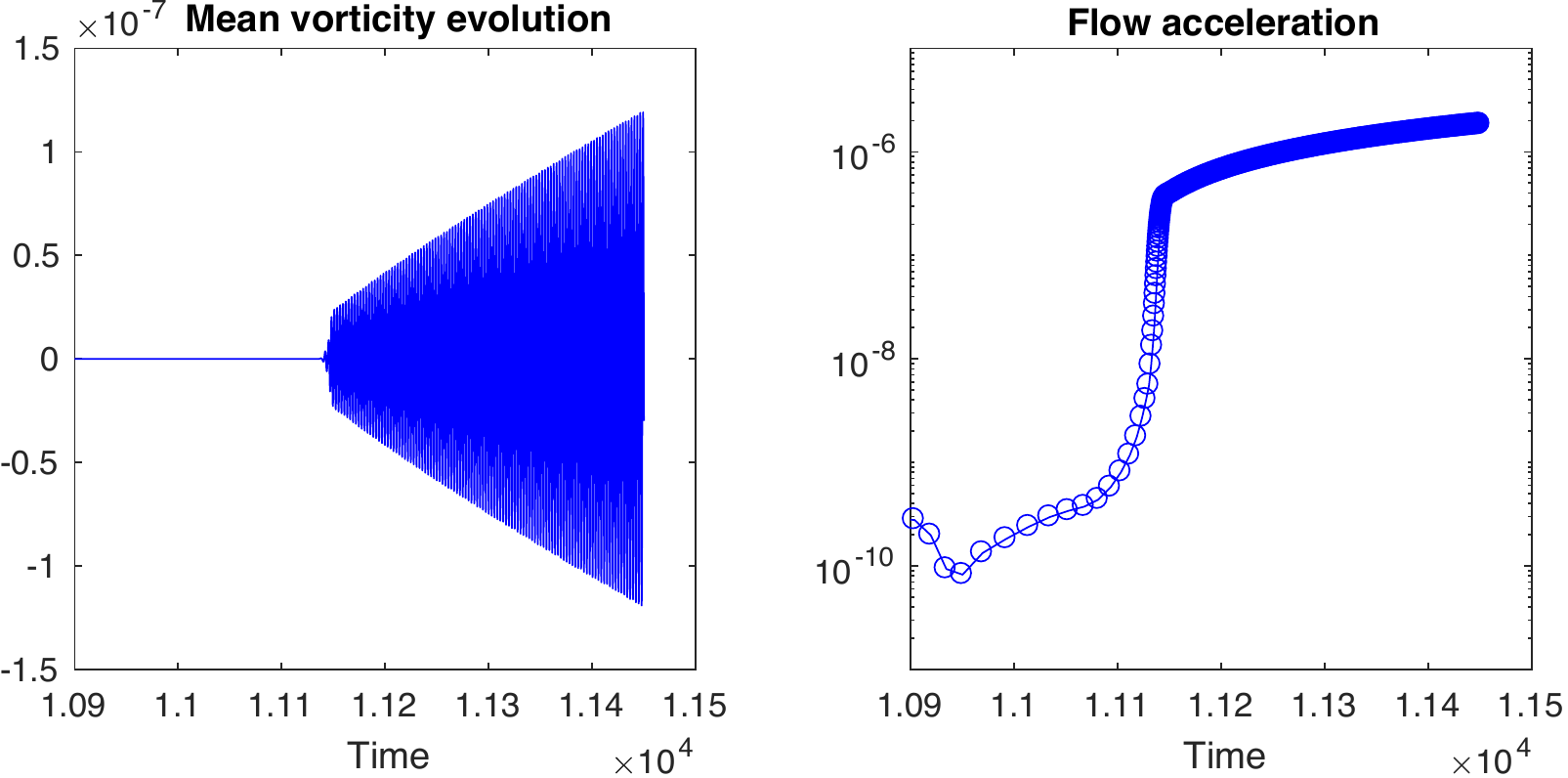}}
\put(160,2){$\nu=1/185.6$}
\end{picture}
\end{center}
\caption{Zoom of flow evolution for $\nu=1/175$ and $\nu=1/185.6$, obstacle flow.} 
\label{second-phase-obstacle-zoom}
\end{figure}

\begin{figure}
\begin{center}
\begin{picture}(270,370) (42,0) 
\put(2,194){\includegraphics[width=.95\textwidth, height=.49\textwidth]{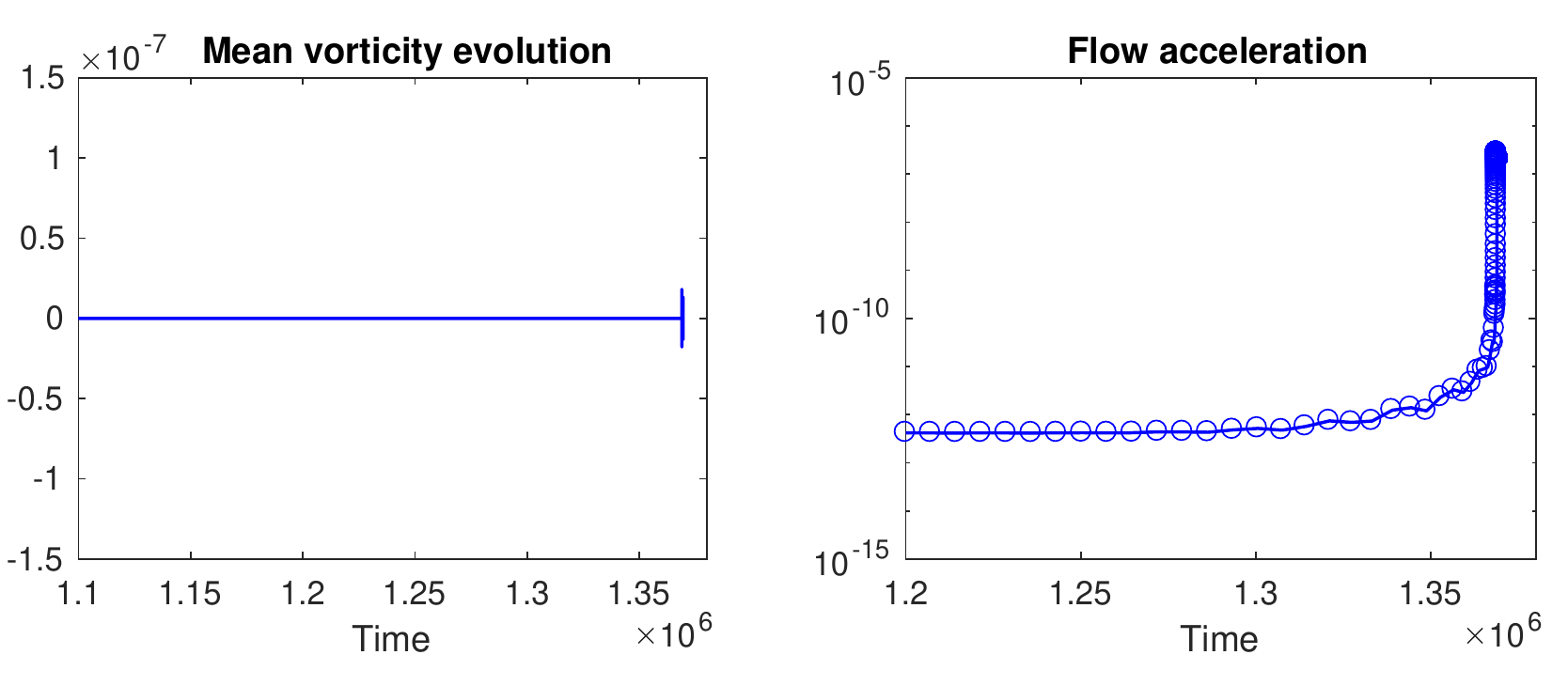}}
\put(135,189){\bf Benign perturbation}
\put(6,5){\includegraphics[width=.99\textwidth, height=.49\textwidth]{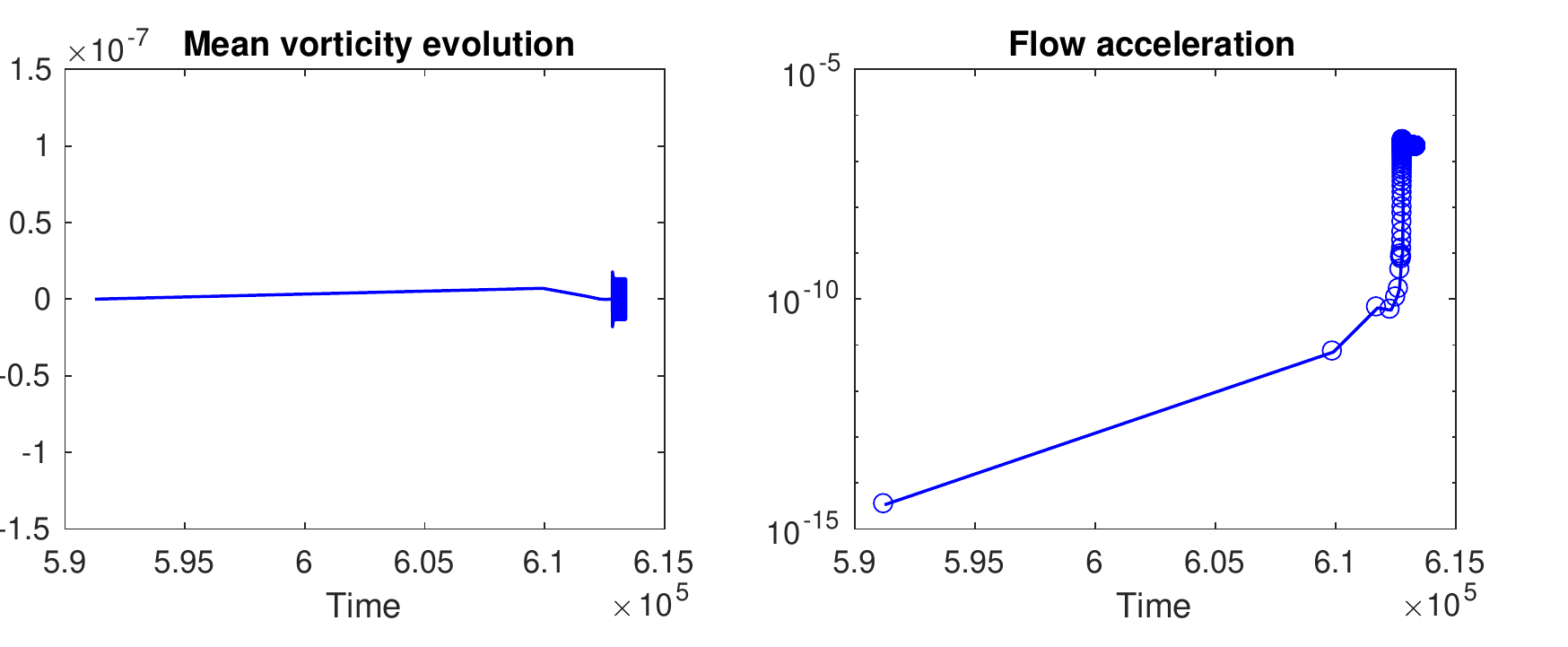}}
\put(81,22){\includegraphics[width=.17\textwidth,height=.22\textwidth]{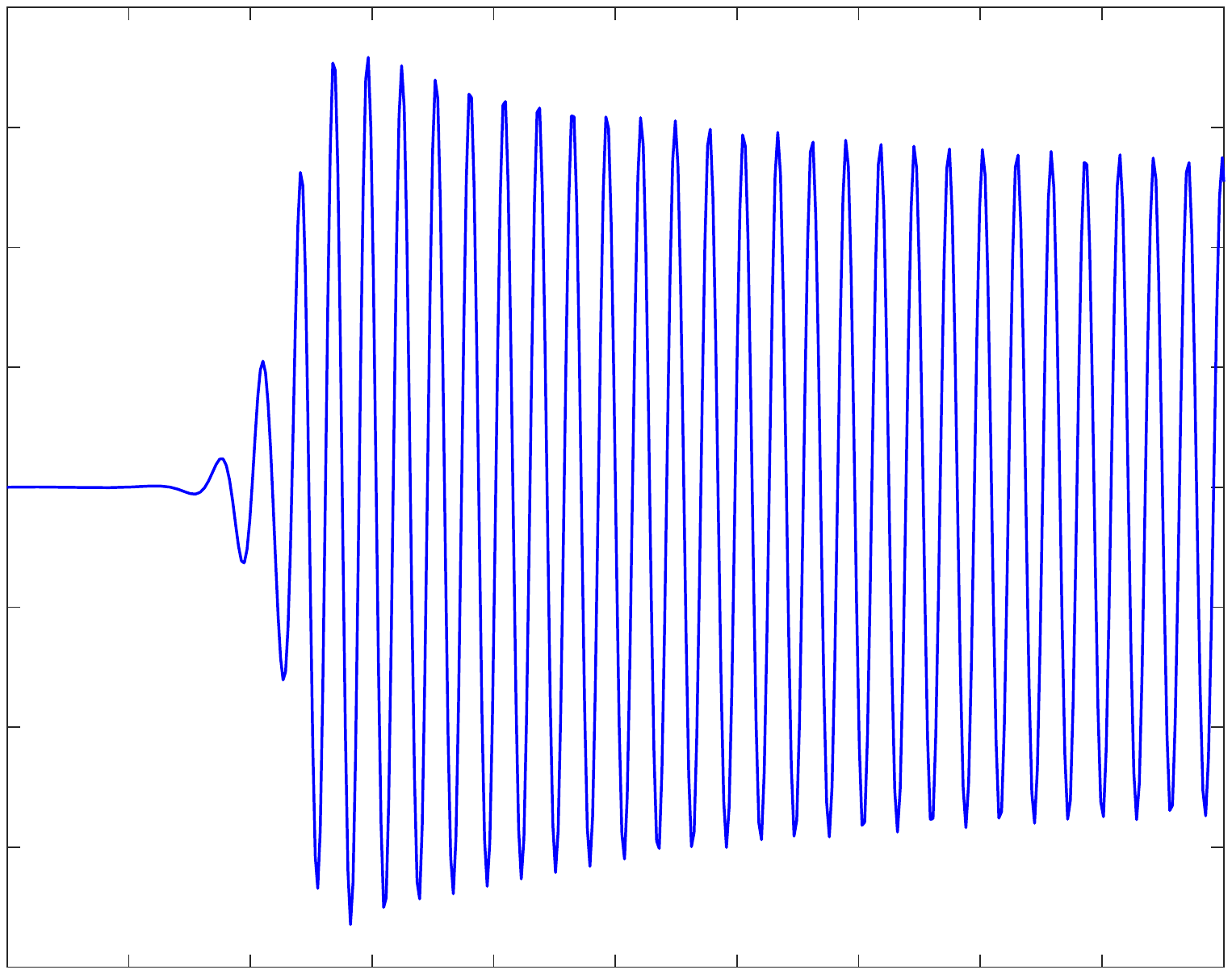}}
\put(135,0){\bf Lively perturbation}
\end{picture}
\end{center}
\caption{Short-term evolution for small perturbations, obstacle flow with $\nu=1/175$.
The inset shows a magnified image of the onset of periodic behavior.} 
\label{second-phase-175}
\end{figure}

Continuing this exploration of subcritical cases, we now consider the effects of
perturbations at the interrupt.
As in the previous section, we look at one perturbation that respects the 
reflectional symmetry of the flow  solution in Figure~\ref{obstacle-figure} and
is  expected  to be ``benign'' and  one  that breaks the reflectional symmetry and
so is expected to be ``lively''.
The results for $\nu=1/175$ are shown in Figure \ref{second-phase-175} and those
for $\nu=1/185.6$ are in Figure \ref{second-phase-1856}.
In these figures, the vertical scaling for the mean vorticities are now set to be
equal in order to discern differences for the two viscosity values.
These images should be compared with those corresponding to analogous experiments 
with no perturbation in
Figures~\ref{fig-evolution-obstacle}--\ref{second-phase-obstacle-zoom}.
In particular, for the sub-critical viscosity $\nu=1/175$ with either type of
perturbation, after a long delay, the solution moves away from a steady state.
This is not surprising, since the same phenomenon occurs when no perturbation is 
done. 
The onset of periodic behavior for the perturbed data is slightly earlier than
for no perturbation (and earlier still for the lively perturbation), but the magnitude
of the oscillations is small.
The results for $\nu=1/185.6$ (Figure \ref{second-phase-1856}) bear some similarity
to these --- most notably, the behavior for the  benign perturbation is virtually 
identical to that for no perturbation (middle of Figure \ref{fig-evolution-obstacle}).
\rblack{The structure of the oscillations for the near-critical viscosity is more
  like that for the super-critical viscosity (for both perturbations as well as
  without perturbation, compare the images for lively perturbation 
in Figure \ref{second-phase-1856} with the images in 
Figures~\ref{fig-evolution-obstacle}--\ref{second-phase-obstacle-zoom}).}
In contrast, for the sub-critical viscosity $\nu=1/175$, the structure of the
oscillations is is more like that  arising when no perturbation is done.
\rblack{But for the lively perturbation and $\nu=185.6$, the onset of unstable behavior 
is significantly earlier (bottom of Figure \ref{second-phase-1856}).}


\begin{figure}
\begin{center}
\begin{picture}(270,370) (42,0) 
\put(3,194){\includegraphics[width=.95\textwidth, height=.49\textwidth]{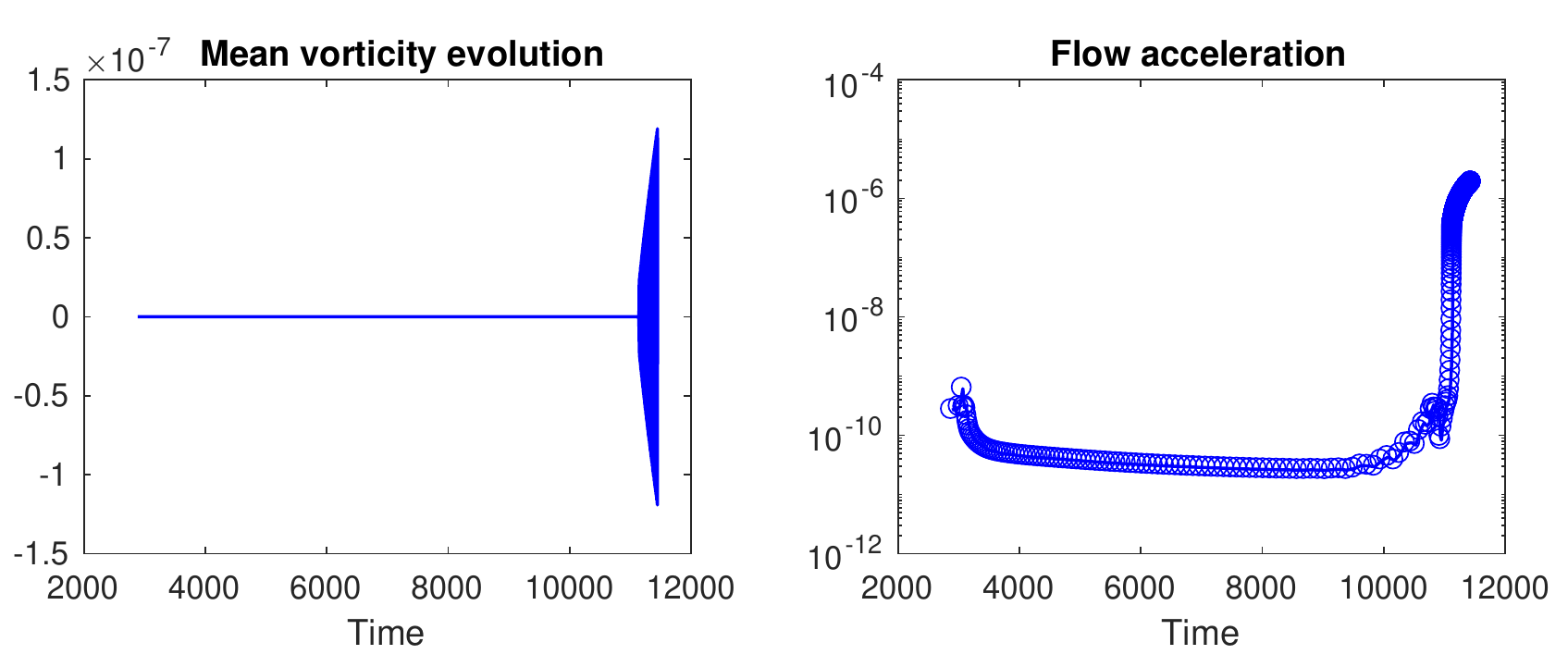}}
\put(135,189){\bf Benign perturbation}
\put(6,3){\includegraphics[width=.95\textwidth, height=.49\textwidth]{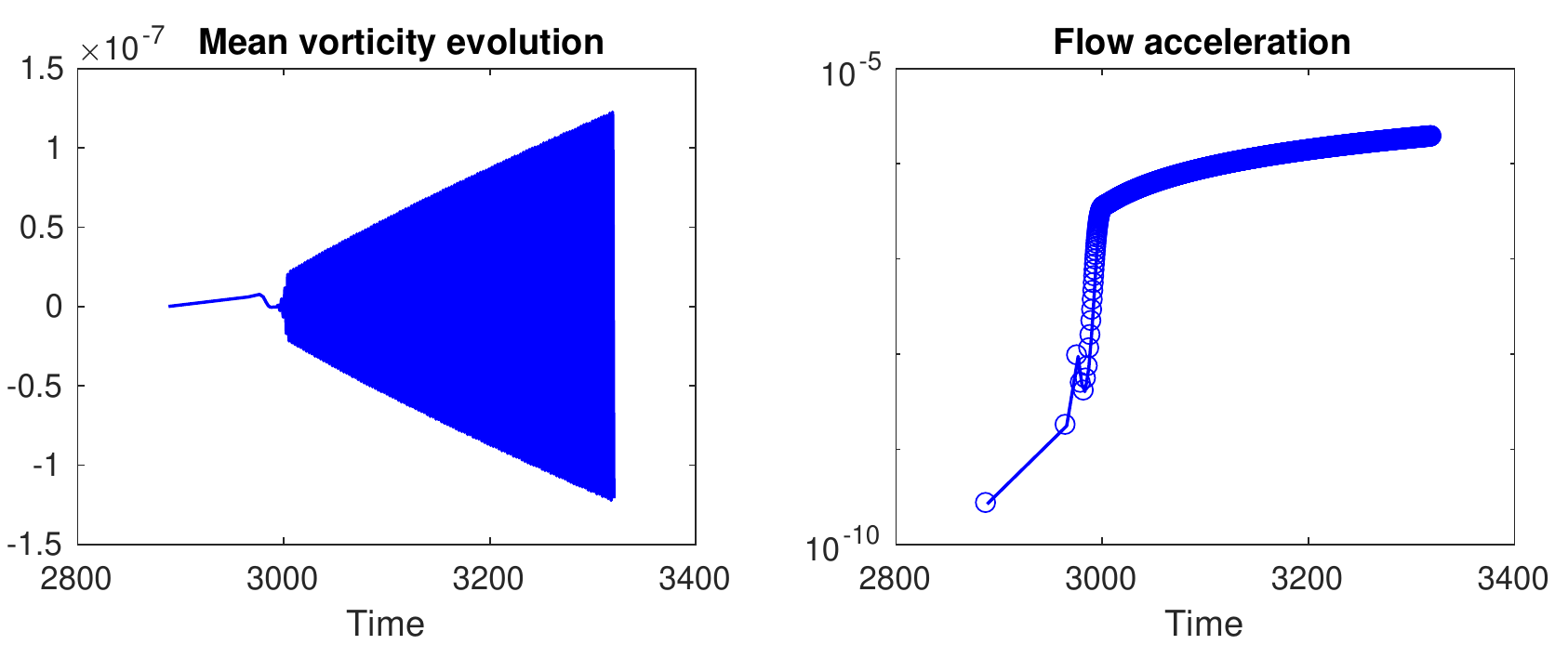}}
\put(21,16){\includegraphics[width=.17\textwidth,height=.22\textwidth]{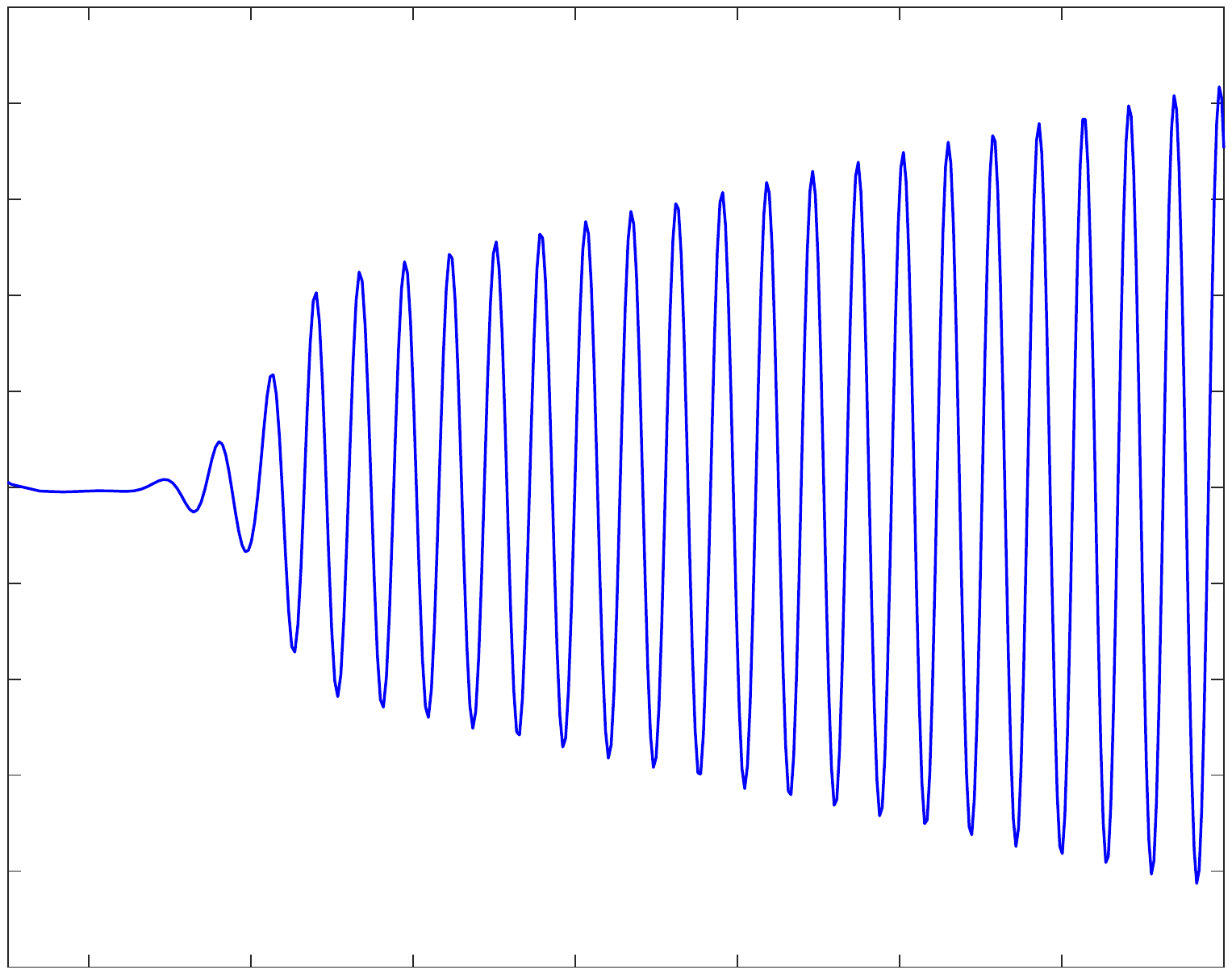}}
\put(135,0){\bf Lively perturbation}
\end{picture}
\end{center}
\caption{Short-term evolution for small perturbations, obstacle flow with $\nu=1/185.6$.
The inset shows a magnified image of the onset of periodic behavior.} 
\label{second-phase-1856}
\end{figure}

\begin{figure}
\begin{center}
\begin{picture}(270,375) (42,0) 
\put(2,194){\includegraphics[width=.95\textwidth, height=.45\textwidth]{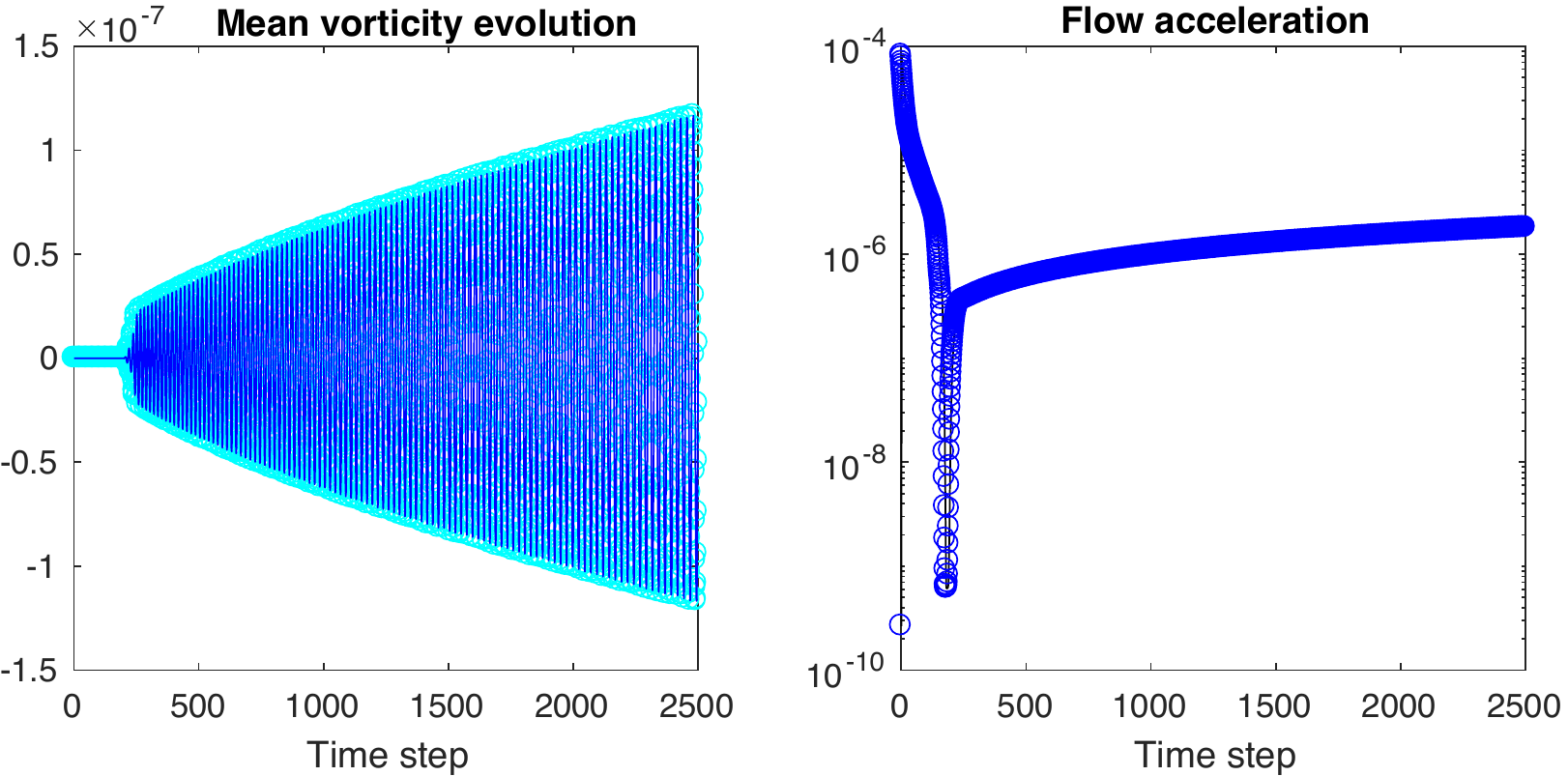}}
\put(135,189){\bf Benign perturbation}
\put(6,3){\includegraphics[width=.95\textwidth, height=.45\textwidth]{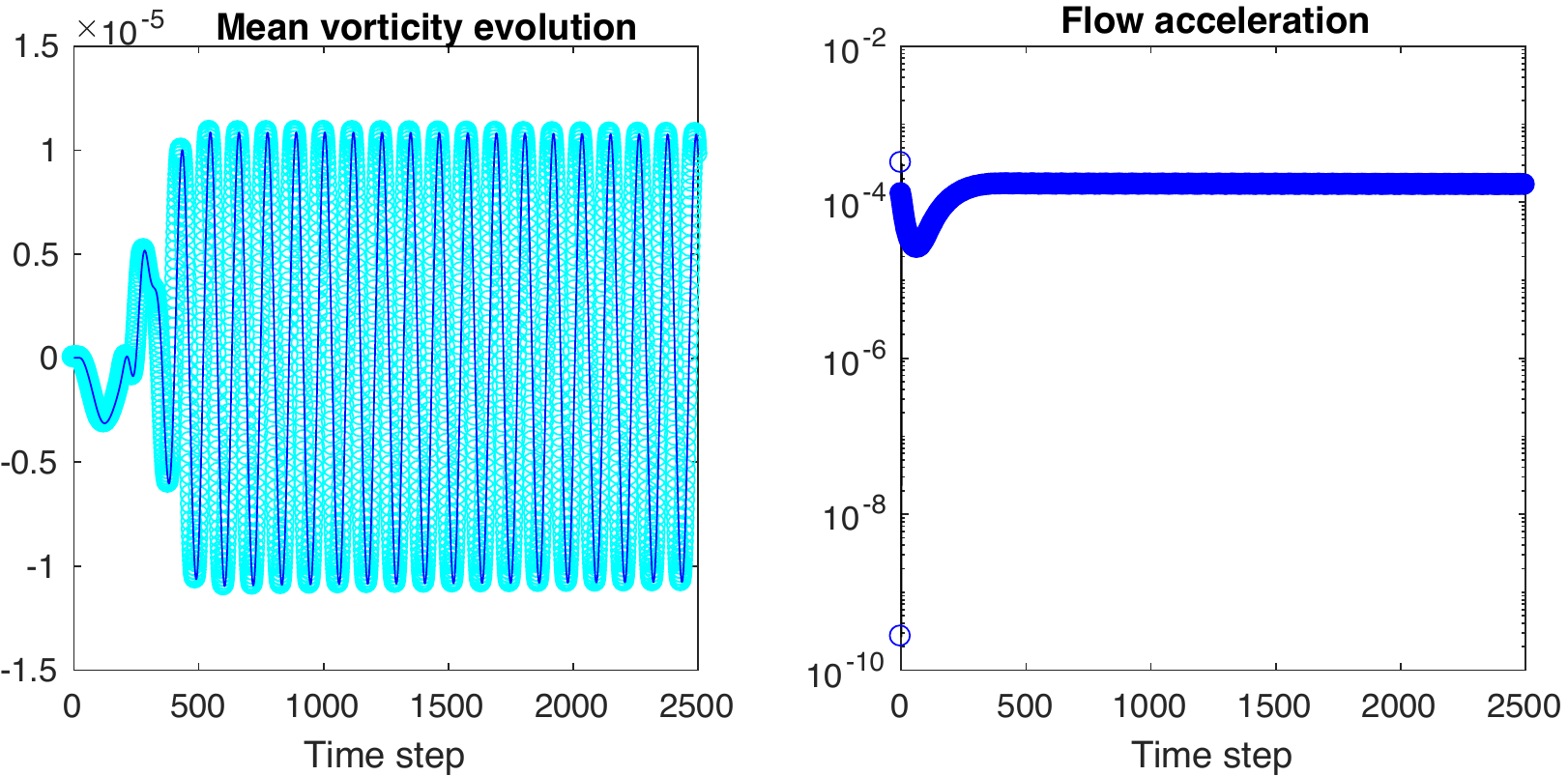}}
\put(135,0){\bf Lively perturbation}
\end{picture}
\end{center}
\caption{Short-term evolution for large perturbations, obstacle flow with $\nu=1/185.6$.} 
\label{second-phase-1856x}
\end{figure}

Finally, when we check to see  what happens   when the perturbation 
is significantly larger (of the order of the perturbation 
made in computing the pseudo-eigenvalues in Figure~\ref{eigs_obstacle}) we observe that
there is a big difference in the time-stepping behavior in any case where  the perturbation is not benign. 
This is illustrated by the results shown  in the bottom plot in Figure~\ref{second-phase-1856x}.
In this case the size of the lively perturbation is big enough to destabilize the integrator and
a noisy periodic solution is computed. This mirrors the vortex-shedding solution that
is computed in the unstable case but has an amplitude that is too small to be
seen when plotted.


\section{Concluding remarks} \label{sect-conclusions}
Our aims in this study were twofold.
First, we have developed a new approach to assess the stability of dynamical systems by constructing
perturbed systems based on collocation methods. 
This is reminiscent of methods for computing pseudospectra, but it has the advantage that  the process
of sampling (approximate) spectra is significantly less costly.
Second, we compared the results of such assessments with the performance of time-stepping 
computations for a nontrivial application, the incompressible Navier-Stokes equations.
In particular, for two benchmark problems, we examined the behavior of a stable integration 
scheme for simulating transient behavior for values of the viscosity in the system that 
are ``sub-critical'', nearly critical (very slightly smaller than the critical value), and super-critical.

In general, we found that the predictions of instability made by the collocation method were consistent
with the behavior of integrators:  in the nearly critical regime (of parameter values, viscosity in this case), 
there is more sensitivity to perturbation than in the sub-critical regime, and outcomes are qualitatively like 
those for super-critical parameters.
We also note that making such assessments is complicated somewhat by the delicate nature 
of computations in regimes at or near stability limits.
Eigenvalues and pseudoeigenvalues are not the sole determining factor affecting stability; the form of 
the perturbation also plays a significant role.

\vspace{.05in}

\vspace{.05in}
\rblack{
       {\bf Acknowledgements:} We thank Mark Embree and an anonymous referee for very
       constructive comments.
}

\bibliography{bib}

\end{document}

In addition, for both test examples, there appears to be some degree of instability caused by limitations
of time accuracy in the integrator; 
we believe this leads to the appearance of symmetry-breaking solutions for the step and 
of periodic solutions for the obstacle problem after long periods of numerical integration.